\documentclass[11pt,reqno]{amsart}

\usepackage{fixltx2e}                    
\usepackage{tgpagella} 
\usepackage[utf8]{inputenc}            
\usepackage{amsmath}                     
\usepackage{amssymb, latexsym, stmaryrd, amsthm, dsfont, amsfonts, amsbsy,amsthm, amsmath, mathrsfs}            
\usepackage{mathtools}                   
\usepackage{bm}                          
\usepackage{enumerate}                   
\usepackage{verbatim}                    
\usepackage{url}   
\usepackage{lscape}                      
\usepackage{soul}
\usepackage{centernot}

\usepackage{bussproofs} 

\usepackage[colorinlistoftodos]{todonotes} 

\usepackage[margin=1.1in]{geometry}

\usepackage{microtype,tikz,tikz-cd}  
\makeatletter                            
\def\MT@register@subst@font{\MT@exp@one@n\MT@in@clist\font@name\MT@font@list
 \ifMT@inlist@\else\xdef\MT@font@list{\MT@font@list\font@name,}\fi}
\makeatother

\usepackage[pdftex,bookmarks,bookmarksnumbered,linktocpage,   
         colorlinks,linkcolor=blue,citecolor=blue]{hyperref}

\newcommand{\bit}{\begin{itemize}}    
\newcommand{\eit}{\end{itemize}}
\newcommand{\ben}{\begin{enumerate}}
\newcommand{\een}{\end{enumerate}}
\newcommand{\benormal}{\ben[\normalfont 1.]}   

\let\enormal\een
\newcommand{\benroman}{\ben[\normalfont (i)]}  
\let\eroman\een

\newcommand{\bde}{\begin{description}}
\newcommand{\ede}{\end{description}}
\let\oper=\mathbb                               

\newcommand{\III}{\oper{I}}                     
\newcommand{\SSS}{\oper{S}}                     
\newcommand{\UUU}{\oper{U}}                     

\newcommand{\VVV}{\oper{V}}                     

\newcommand{\HHH}{\oper{H}}

\newcommand{\PPP}{\oper{P}}

\newcommand{\PPU}{\oper{P}_{\!\textsc{u}}^{}}

\theoremstyle{theorem}
\newtheorem{Theorem}{Theorem}[section]
\newtheorem{Theorem-n}{Theorem}

\newtheorem{Proposition}[Theorem]{Proposition}
\newtheorem{Modal Sahlqvist Theorem}[Theorem]{Modal Sahlqvist Theorem}
\newtheorem{Intuitionistic Sahlqvist Theorem}[Theorem]{Intuitionistic  Sahlqvist Theorem}
\newtheorem{Intuitionistic Sahlqvist Canonicity Theorem}[Theorem]{Intuitionistic Sahlqvist Canonicity Theorem}

\newtheorem{Transfer Lemma}[Theorem]{Transfer Lemma}
\newtheorem{Abstract Sahlqvist Theorem}[Theorem]{Abstract Sahlqvist Theorem}

\newtheorem{Lemma}[Theorem]{Lemma}
\newtheorem{Corollary}[Theorem]{Corollary}
\newtheorem{Claim}[Theorem]{Claim}

\theoremstyle{definition}
\newtheorem{Definition}[Theorem]{Definition}
\newtheorem{exa}[Theorem]{Example}

\theoremstyle{remark}

\newtheorem{Remark}[Theorem]{Remark}

\let\leq=\leqslant
\let\nleq=\nleqslant
\let\geq=\geqslant 

\usepackage[mathscr]{euscript}
 \let\mathscr\relax
\usepackage[scr]{rsfso}

\newcommand{\powerset}[1]{\mathscr{P}({#1})}
\newcommand{\powersetm}[1]{\mathscr{P}_{\mathsf{M}}({#1})}

\newcommand{\up}[1]{\mathsf{Up}({#1})}

\newcommand{\sub}{\subseteq}

\newcommand{\Fi}[1]{\mathsf{Fi}_\vdash({#1})}
\newcommand{\Fic}[1]{\mathsf{Fi}_\vdash^\omega({#1})}
\newcommand{\fg}[1]{\mathsf{Fg}^{\boldsymbol{A}}_\vdash({#1})}
\newcommand{\fgL}[1]{\mathsf{Fg}^{\boldsymbol{A}}_\mathsf{ILL}({#1})}

\newcommand{\fgt}[1]{\mathsf{Fg}_\vdash({#1})}

\newcommand{\Th}[1]{\mathsf{Th}({#1})}
\newcommand{\Thc}[1]{\mathsf{Th}^\omega({#1})}

\bmdefine{\A}{A} 
\bmdefine{\C}{C}                                
                              
\bmdefine{\2}{2}
\bmdefine{\B}{B}
\bmdefine{\D}{D}
\newcommand{\?}{\ensuremath{\mkern0.4\thinmuskip}}   

\DeclareMathOperator*{\foo}{\scalerel*{\&}{\sum}}
\usepackage{scalerel}



\begin{document}

\title[Intuitionistic Sahlqvist theory for deductive systems]{Intuitionistic Sahlqvist theory for deductive systems}


\author{Damiano Fornasiere and Tommaso Moraschini}

\address{Damiano Fornasiere: Departament de Filosofia, Facultat de Filosofia, Universitat de Barcelona (UB), Carrer Montalegre, $6$, $08001$ Barcelona, Spain}\email{damiano.fornasiere@ub.edu}
\address{Tommaso Moraschini: Departament de Filosofia, Facultat de Filosofia, Universitat de Barcelona (UB), Carrer Montalegre, $6$, $08001$ Barcelona, Spain}\email{tommaso.moraschini@ub.edu}

\date{\today}

\begin{abstract}
Sahlqvist theory is extended to the fragments of the intuitionistic propositional calculus that include the conjunction connective. This allows us to introduce a Sahlqvist theory of intuitionistic character amenable to arbitrary protoalgebraic deductive systems. As an application, we obtain a Sahlqvist theorem for the fragments of the intuitionistic propositional calculus that include the implication connective and for the extensions of the intuitionistic linear logic.
\end{abstract}

\maketitle

\section{Introduction}

Sahlqvist theorem is one of modal logic's crown jewels \cite{Sa75}. The theorem revolves around a family of syntactically defined modal formulas, known as \emph{Sahlqvist formulas}, and consists of two halves, related to the phenomena of \emph{canonicity} and \emph{correspondence}, respectively. The canonicity part of the theorem states that the validity of Sahlqvist formulas is preserved by canonical extensions of modal algebras \cite{JonTar51,JonTar52}, while the correspondence part states that the class of Kripke frames validating a Sahlqvist formula is strictly elementary. Taken together, the two imply that every normal modal logic axiomatized by Sahlqvist formulas is complete with respect to an elementary class of Kripke frames.

Whilst Sahlqvist theory has been at the center of many investigations in modal logic (see, e.g., \cite{Kr93,SaVa89}), it received less attention in the setting of the intuitionistic propositional calculus $\mathsf{IPC}$, with the notable exception of \cite{GhiMe97} (see also \cite{CoPaZa19}).\   This should not come as a surprise, as a version of Sahlqvist theory for $\mathsf{IPC}$ can easily be derived from the modal one, utilizing the \emph{G\"odel-McKinsey-Tarski translation}  of $\mathsf{IPC}$ into the modal system $\mathsf{S4}$ \cite{Go32,McKT48} and its semantic interpretation \cite{MakRyb74}, as explained in \cite{CoPaZa19}.  However, this method breaks down for fragments of $\mathsf{IPC}$, mostly because their duality theory is more opaque than that of $\mathsf{IPC}$. 

In this paper, we will fill this gap by extending Sahlqvist theory to fragments\footnote{The language of $\mathsf{IPC}$ is assumed to be	$\land, \lor, \to, \lnot , 0,  1$.} of $\mathsf{IPC}$ including the conjunction connective (Theorem \ref{thm : enhanced intuitionistic sahlqvist}). This will serve as the basis for our main contribution, which consists of a Sahlqvist theory of intuitionistic character amenable to arbitrary deductive systems (Theorem \ref{thm : correspondence}).

More precisely,  we say that a formula of $\mathsf{IPC}$ is a \emph{Sahlqvist antecedent} if it is constructed from negative formulas and the constants $0$ and $1$ using only $\land$ and $\lor$. It is a \emph{Sahlqvist implication} if it is either positive, or the negation of a Sahlqvist antecedent, or of the form $\varphi \to \psi$, where $\varphi$ is a Sahlqvist antecedent and $\psi$ a positive formula. Lastly, formulas obtained from Sahlqvist implications using only $\land$ and $\lor$ will be called \emph{Sahlqvist formulas}.

Instead of working with Sahlqvist formulas, we prefer to focus on \emph{Sahlqvist quasiequations}, i.e., expressions of the form 
	\[
	\varphi_1 \land y \leq z \, \& \dots \& \, \varphi_n \land y \leq z \Longrightarrow y \leq z,
\]
where $\varphi_1, \dots, \varphi_n$ are Sahlqvist formulas and $y$ and $z$ are distinct variables that do not occur in them. This is because, while in Heyting algebras (i.e., the algebraic models of $\mathsf{IPC}$) Sahlqvist quasiequations and formulas are equally expressive, in the case of fragments the expressive power of Sahlqvist quasiequations  is often greater. For instance, the so-called \color{black} \emph{bounded top width $n$} axiom \cite{Smory73}
\[
\mathsf{btw}_n \coloneqq \bigvee_{i = 1}^{n+1} \lnot (\lnot x_i \land \bigwedge_{0 <j< i}x_j)
\]
can be rendered as the Sahlqvist quasiequation
\begin{equation}\label{Eq:intro-example}
\Phi_n = \foo_{1 \leq i \leq n+1}\Big(\lnot(\lnot x_i \land \bigwedge_{0 < j< i}x_j) \land y \leq z\Big) \Longrightarrow y \leq z,
\end{equation}
which, contrarily to the formula $\mathsf{btw}_n$, does not contain the disjunction connective and, therefore, \color{black} makes the concept of ``bounded top width $n$'' amenable also to the algebraic models of the $\langle \land, \lnot \rangle$-fragment of $\mathsf{IPC}$, i.e., pseudocomplemented semilattices. Furthermore, the role of the quasiequation $\Phi_n$ cannot be taken over by any formula or equation in $\land$ and $\lnot$ only, because the expressive power of the latter is extremely limited \cite{Jon-PhD}.

Recall that the lattice $\mathsf{Up}(\mathbb{X})$ of upsets of a poset $\mathbb{X}$ can be endowed with the structure of a Heyting algebra (see, e.g., \cite{ChZa97,Esakia-book85}). Furthermore, a formula of $\mathsf{IPC}$ is valid   in a poset $\mathbb{X}$ viewed as a Kripke frame \color{black} iff it is valid in the Heyting algebra $\mathsf{Up}(\mathbb{X})$. Lastly, given a semilattice $\A = \langle A; \land \rangle$, we denote by $\A_\ast$ the poset of its meet irreducible filters. The Sahlqvist theorem for fragments of $\mathsf{IPC}$ with $\land$ takes the following form  (Theorems \ref{Thm:IPC-Sahlqvist} and \ref{thm : enhanced intuitionistic sahlqvist}):

\begin{Theorem-n}\label{Thm:1-intro}
The following conditions hold for a Sahlqvist quasiequation $\Phi$ in a language $\mathcal{L} \subseteq \{ \land, \lor, \to, \lnot, 0, 1 \}$  containing $\land$:
\benroman
\item \emph{Canonicity}: If an $\mathcal{L}$-subreduct $\A$ of a Heyting algebra validates $\Phi$, then also $\mathsf{Up}(\A_\ast)$ validates $\Phi$;
\item \emph{Correspondence}: There exists an effective computable first order sentence $\mathsf{tr}(\varphi)$ in the language of posets such that $\mathsf{Up}(\mathbb{X}) \vDash \Phi$ iff $\mathbb{X} \vDash \mathsf{tr}(\Phi)$, for every poset $\mathbb{X}$.
\eroman 
\end{Theorem-n}

The main obstacles in proving Theorem \ref{Thm:1-intro} can be summarized as follows. On the one hand, the method of \cite{CoPaZa19} is based on the observation that a Heyting algebra $\A$ validates a formula $\varphi$ iff the free Boolean extension of $\A$, viewed as a modal algebra, validates the G\"odel-McKinsey-Tarski translation of $\varphi$. This property, however, need not hold in subreducts of Heyting algebras, preventing the applicability of a similar method. On the other hand, the algebraic models of fragments of $\mathsf{IPC}$ with $\land$ admit the presence of operations that fail to be order preserving in every coordinate (such as the negation or the implication), contrarily to the case of \cite{KKTZ17}. Similarly, they need not have a lattice structure and, therefore, do not fall under the scope of \cite{CeJa99a,ACCNDL,GeNaVe05}. 

Lastly, even when these models have a lattice structure (as in the case of \textit{finite} pseudocomplemented semilattices), they need not be \textit{distributive} lattices. Because of this reason, their traditional canonical canonical extensions \cite{DuGePa05,GeHa01,GeJo94} may fail to be Heyting algebras. For instance, the smallest nonmodular lattice $\boldsymbol{N}_5$ can be viewed as a pseudocomplemented semilattice, whose canonical extension is order isomorphic to $\boldsymbol{N}_5$ itself and, therefore, is not a Heyting algebra. To overcome this problem, in the canonicity part of the theorem we work with completions of the form $\mathsf{Up}(\A_\ast)$ which, in turn, are always Heyting algebras.\footnote{Nevertheless, we kept the expression \textit{canonicity} for historical reasons.} As a consequence, the order theoretic properties typical of canonical extensions which serve as the basis of the approach of \cite{GhiMe97} and \cite{ACCNDL,ConTal20LMSC,DeRPa} need not hold in our setting: for instance, the completions of the form $\mathsf{Up}(\A_\ast)$ need not be \textit{dense} in the sense of \cite{DuGePa05} and are not induced by a polarity of filters and ideals in the sense of \cite{Delta1}.

Our main tools are \color{black} a model theoretic observation on universal classes (Theorem \ref{Thm:general-embedding}) and the correspondence in Section \ref{Sec:adjunction} between algebraic homomorphisms on the one hand and partial order preserving maps that generalize the notion of a p-morphism typical of Esakia duality \cite{Es74,Esakia-book85} on the other hand (see, e.g., \cite{BeBe09,Kh81,Zakha89,Za92,Za96}).\color{black}

As we mentioned, part of the interest of Theorem \ref{Thm:1-intro} is that it contains the germ of a Sahlqvist theory amenable to arbitrary \emph{logics} (a.k.a.\ \emph{deductive systems}), i.e., finitary and substitution invariant consequence relations on the set of formulas of some algebraic language. The price to pay in exchange for the great generality, however, is that the resulting Sahlqvist theory is of intuitionistic character and, therefore, does not apply to logics whose meet irreducible theories are maximally consistent (such as normal modal logics). 

This generalization is made possible by the methods of abstract algebraic logic \cite{AAL-AIT-f}, which allow to recognize that the pillars sustaining the intuitionistic Sahlqvist theory are certain metalogical properties that govern the behavior of the intuitionistic connectives $\lnot, \to$, and $\lor$. More precisely, a logic $\vdash$ is said to have
\benroman
\item the \emph{inconsistency lemma} \cite{JGR13} when for every $n \in \mathbb{Z}^+$ there exists a finite set of formulas ${\thicksim_{n}\!\!(x_1, \dots, x_n)}$ such that for every finite set of formulas $\Gamma \cup \{ \varphi_1, \dots, \varphi_n \}$,
\[
\Gamma \cup \{\varphi_1, \dots, \varphi_n \} \text{ is inconsistent}\, \, \text{ iff } \, \, \Gamma \vdash\, \sim_n\!\!(\varphi_1, \dots, \varphi_n);
\]
\item the \emph{deduction theorem} \cite{BP?} when there exists a finite set of formulas $x \Rightarrow y$ such that for every finite set of formulas $\Gamma \cup \{ \psi, \varphi \}$,
\[
\Gamma, \psi \vdash \varphi \, \, \text{ iff } \, \, \Gamma \vdash \psi \Rightarrow \varphi;
\]
\item the \emph{proof by cases} \cite{Cz84a,CzDz} when there exists a finite set $x \curlyvee y$ of formulas such that for every finite set of formulas $\Gamma \cup \{ \psi, \varphi, \gamma \}$,
\[
\Gamma, \psi \vdash \gamma \text{ and }\Gamma, \varphi \vdash \gamma \, \, \text{ iff } \, \, \Gamma, \psi\curlyvee \varphi \vdash \gamma.
\]
\eroman
 It is well known that $\mathsf{IPC}$ has the inconsistent lemma, the deduction theorem, and the proof by cases, as witnessed, respectively, by the sets
\[
\thicksim_{n}\!\!(x_1, \dots, x_n) \coloneqq \{ \lnot (x_1 \land \dots \land x_n) \} \qquad x \Rightarrow y \coloneqq \{ x \to y \} \qquad x \curlyvee y \coloneqq \{ x \lor y \}.
\]

Accordingly, when a logic $\vdash$ possesses the metalogical properties governing the behavior of the connectives among $ \lnot, \to$, and $\lor$ appearing in a formula $\varphi(x_1, \dots, x_n)$ of $\mathsf{IPC}$, we say that $\varphi$ is \emph{compatible} with $\vdash$. In this case, with every $k \in \mathbb{Z}^+$ we can associate a finite set of formulas
\[
\boldsymbol{\varphi}^k(x_1^1, \dots, x_1^k, \dots, x_n^1, \dots, x_n^k)
\]
of $\vdash$ which globally behaves as $\varphi$. For instance, suppose that $\varphi = \lnot \psi$ and that we already defined $\boldsymbol{\psi}^k$ to be $\{ \chi_1, \dots, \chi_n \}$. Since the connective $\lnot$ appears in $\varphi$, the assumption that $\varphi$ is compatible with $\vdash$ guarantees that the latter has the inconsistency lemma. Accordingly, we set
\[
\boldsymbol{\varphi}^k =  (\boldsymbol{\lnot \psi})^k \coloneqq \?\? \sim_n\!\!(\chi_1, \dots, \chi_n),
\]
thus ensuring that $\boldsymbol{\varphi}^k$ behaves as the negation of $\boldsymbol{\psi}^k = \{ \chi_1, \dots, \chi_n \}$ in $\vdash$.

Our main result applies to logics $\vdash$ that are \emph{protoalgebraic}, i.e., that possess a nonempty set of formulas $\Delta(x, y)$ which globally behaves as a weak implication, in the sense that $\emptyset \vdash \Delta(x, x)$ and modus ponens $x, \Delta(x, y) \vdash y$ hold \cite{Cz01}. It takes the form of a correspondence theorem connecting the validity of certain metarules in a logic $\vdash$ with the structure of the posets $\mathsf{Spec}_\vdash(\A)$ of meet irreducible deductive filters of $\vdash$ on arbitrary algebras $\A$ (Theorem \ref{thm : correspondence}).

\begin{Theorem-n}\label{Thm:2-intro}
Let $\varphi_1 \land y \leq z \, \& \dots \& \, \varphi_n \land y \leq z \Longrightarrow y \leq z
$ be a Sahlqvist quasiequation such that $\varphi_1, \dots, \varphi_n$ are compatible with a protoalgebraic logic $\vdash$. Then $\vdash$ validates all the metarules of the form
\begin{prooftree}
	\AxiomC{$\Gamma, \boldsymbol{\varphi_1}^k(\vec{\gamma}_1, \dots, \vec{\gamma}_n) \rhd \psi$} 
	\AxiomC{$\dots$}
	\AxiomC{$\Gamma, \boldsymbol{\varphi_m}^k(\vec{\gamma}_1, \dots, \vec{\gamma}_n) \rhd \psi$}
	\TrinaryInfC{$\Gamma \rhd \psi$}
\end{prooftree}
iff the poset $\mathsf{Spec}_{\vdash}(\A)$ validates $\mathsf{tr}(\Phi)$, for every algebra $\A$.
\end{Theorem-n}

For instance, a protoalgebraic logic with the inconsistency lemma validates the metarules corresponding to the bounded top width $n$ Sahlqvist quasiequation in Condition (\ref{Eq:intro-example}) iff the principal upsets in $\mathsf{Spec}_{\vdash}(\A)$ have at most $n$ maximal elements, for every algebra $\A$ (Theorem \ref{Thm:BTWL-correspondence-AAL}). In the case where $n=1$, this was first proved in \cite{LMR22} (see also \cite{PrenLav20}).

The connection between Theorems \ref{Thm:1-intro} and \ref{Thm:2-intro} is made possible by a series of bridge theorems that connect the validity of the inconsistency lemma, the deduction theorem, and the proof by cases in a protoalgebraic logic $\vdash$ with the demand that the semilattices $\Fic{\A}$ of compact deductive filters of $\vdash$ on algebras $\A$ are subreducts of Heyting algebras in a suitable language containing $\land$. For instance, a protoalgebraic logic $\vdash$ has the inconsistency lemma iff $\Fic{\A}$ is a pseudocomplemented semilattice, for every algebra $\A$ \cite{JGR13}. A similar result, where implicative semilattices and distributive lattices take over the role of pseudocomplemented semilattices, holds for the deduction theorem and the proof by cases	\cite{BP91, BP?, BP-AAL-DDT,CN13,Cz84a,CzDz}. This allows us to apply Theorem \ref{Thm:1-intro} to the semilattices of the form $\Fic{\A}$. Together with the observation that the poset $\Fic{\A}_\ast$ of meet irreducible filters of $\Fic{\A}$ is isomorphic to $\mathsf{Spec}_{\vdash}(\A)$, these are the keys for extending Theorem \ref{Thm:1-intro} to arbitrary protoalgebraic logics.
	
Lastly, we come full circle and derive a version of Sahlqvist theory for fragments of $\mathsf{IPC}$ including the implication connective $\to$ from Theorem \ref{Thm:2-intro}. As the correspondence part of the theory is already supplied by Theorem \ref{Thm:1-intro}, our result takes the form of a canonicity theorem. More precisely, given a language $\mathcal{L} \subseteq \{ \land, \lor, \to, \lnot, 0, 1 \}$ containing $\to$ and an $\mathcal{L}$-subreduct $\A$ of a Heyting algebra, we denote by $\A_\ast$ the poset of meet irreducible implicative filters of $\A$. Since every Sahlqvist quasiequation $\Phi$ in $\mathcal{L}$ can be rendered (up to equivalence in Heyting algebras) as a set $\mathsf{A}(\Phi)$ of formulas of $\mathcal{L}$, our canonicity theorem can be phrased as follows (Theorem \ref{Thm:canonicity-arrow-IPC}):

\begin{Theorem-n}
Let $\Phi$ be a Sahlqvist quasiequation in a language $\mathcal{L} \subseteq \{ \land, \lor, \to, \lnot, 0, 1 \}$ containing $\to$. If an $\mathcal{L}$-subreduct $\A$ of a Heyting algebra validates $\mathsf{A}(\Phi)$, then also $\mathsf{Up}(\A_\ast)$ validates $\mathsf{A}(\Phi)$.
\end{Theorem-n}

We close the paper with another application of Theorem \ref{Thm:2-intro}, namely a correspondence result for intuitionistic linear logic (Theorem \ref{Thm:linear-main}) which differs from the one in \cite{Suzukia,Suzuki} in that, while our theorem captures only the intuitionistic aspects of this logic, it extends naturally to its axiomatic extensions.

\section{Pseudocomplemented and implicative semilattices}

\begin{Definition} An algebra $\langle A; \land \rangle$ is said to be a \emph{semilattice} when $\land$ is an associative, commutative, and idempotent binary operation. 
\end{Definition}

With each semilattice $\langle A; \land \rangle$ we can associate a partial order $\leq$ on $A$ defined, for every $a, b \in A$, as follows:
\begin{equation}\label{Eq:order-semilattice-meet}
a \leq b \, \, \Longleftrightarrow\, \, a \land b = a.
\end{equation}
In this case, $\langle A; \leq \rangle$ is a poset in which the binary meet of every pair of elements $a, b$ exists and coincides with the element $a \land b$. Furthermore, given a poset $\langle A; \leq \rangle$ in which binary meets exist, the pair $\langle A; \land \rangle$, where $\land$ is the operation of taking binary meets, is a semilattice. These transformations are one inverse to the other.

For the present purpose, two kinds of semilattices are of special interest (see, e.g., \cite{Fr62,Ne65}):
\begin{Definition}
A semilattice $\langle A; \land \rangle$ is said to be:
\benroman
\item \emph{Pseudocomplemented} if it has a minimum element $0$ and for each $a \in A$ there exists an element $\lnot a \in A$ such that for every $c \in A$,
\begin{equation}\label{Eq:neg-res-law}
c \land a = 0 \, \, \Longleftrightarrow \, \, c \leq \lnot a;
\end{equation}
\item \emph{Implicative} if for each $a, b \in A$ there exists an element $a \to b \in A$ such that for every $c \in A$,
\begin{equation}\label{Eq:res-law}
c \land a \leq b \, \,  \Longleftrightarrow \, \, c \leq a \to b.
\end{equation}
\eroman
\end{Definition}

It follows that a semilattice $\langle A; \land \rangle$ is pseudocomplemented if it has a minimum element $0$ and for each $a \in A$ there exists the largest $c \in A$ such that $a \land c = 0$ (in which case, we take it to be $\lnot a$). As a consequence, every pseudocomplemented semilattice has a maximum element, namely $1 \coloneqq \lnot 0$.

Similarly, a semilattice $\langle A; \land \rangle$ is implicative if for each $a, b \in A$ there exists the largest $c \in A$ such that $c \land a \leq b$ (in which case, we take it to be $a \to b$). As a consequence, implicative semilattices $\langle A; \land \rangle$ have always a maximum, namely $a \to a$ for an arbitrary $a \in A$. Because of this, an implicative semilattice is said to be \emph{bounded} when it has a minimum element $0$. Notably, every bounded implicative semilattice $\langle A; \land \rangle$ is pseudocomplemented, since Condition (\ref{Eq:neg-res-law}) holds setting $\lnot a \coloneqq a \to 0$, for every $a \in A$.

Since for every lattice $\langle A; \land, \lor \rangle$ the pair $\langle A; \land \rangle$ is a semilattice, the above terminology extends naturally to lattices \cite{BaDw74,Esakia-book85,RaSi70}:

\begin{Definition}
A lattice $\langle A; \land, \lor \rangle$ is said to be \emph{pseudocomplemented} (resp.\ \emph{implicative}) if so is $\langle A; \land \rangle$. Bounded implicative lattices will be called \emph{Heyting algebras}. 
\end{Definition}

It is well known that the lattice reduct of an implicative lattice is always distributive.

\begin{Remark}\label{Rem:algebras}
Sometimes it will be convenient to treat pseudocomplemented and implicative semilattices as algebras whose basic operations include $\lnot$, $\to$, $0$, and $1$ (as opposed to $\land$ only). When this is the case, we will assume that pseudocomplemented semilattices are algebras $\langle A; \land, \lnot, 0, 1 \rangle$, where $\langle A; \land \rangle$ is a semilattice with minimum $0$ and maximum $1$ and $\lnot$ a unary operation on $A$ satisfying Condition (\ref{Eq:neg-res-law}). Similarly, we will treat implicative semilattices as algebras $\langle A; \land, \to, 1 \rangle$, where $\langle A; \land \rangle$ is a semilattice with maximum $1$ and $\to$ a binary operation on $A$ satisfying Condition (\ref{Eq:res-law}).\ Lastly, bounded implicative semilattices will be algebras $\langle A; \land, \to, 0, 1 \rangle$, where $\langle A; \land, \to, 1 \rangle$ is an implicative semilattice with minimum $0$. The analogous conventions apply to pseudocomplemented lattices, implicative lattices, and Heyting algebras with the only difference that the language of these structures will be assumed to contain the join operation $\lor$.
\qed
\end{Remark}

When $\langle A; \land \rangle$ is a finite semilattice with a maximum element, the partial order $\langle A; \leq \rangle$ is a lattice in which
\[
a \lor b = \bigwedge \{ c \in A : a, b \leq c \}, \text{ for every }a, b \in A.
\]
Because of this, Condition (\ref{item:finite-BSL2}) in the following result makes sense.

\begin{Proposition}\label{Prop:finite-BSL}
The following conditions hold:
\benroman
\item\label{item:finite-BSL1} If $\A = \langle A; \land, \lor, \lnot, 0, 1 \rangle$ is a finite pseudocomplemented distributive lattice, the structure $\langle A; \land, \lor, \to, 0, 1  \rangle$, where $\to$ is defined as $a \to b = \max \{ c \in A : c \land a \leq b \}$, is a Heyting algebra in which the term function $x \to 0$ coincides with the operation $\lnot$ of $\A$;
\item\label{item:finite-BSL2} If $\boldsymbol{A} = \langle A; \land, \to, 1 \rangle$ is a finite implicative semilattice, the structure $\langle A; \land, \lor, \to, 0, 1  \rangle$, where $\lor$ and $0$ are, respectively, the join operation and the minimum element of $\langle A; \leq \rangle$, is a Heyting algebra.
\eroman
\end{Proposition}

\begin{proof}
Condition (\ref{item:finite-BSL1}) holds because every finite distributive lattice is a Heyting algebra and the structure of a Heyting algebra is uniquely determined by its order reduct. For Condition (\ref{item:finite-BSL2}), see, e.g., \cite{Kh81}.
\end{proof}

\begin{Remark}
In contrast to this, finite pseudocomplemented semilattices cannot be given in general the structure of a Heyting algebra, because they need not be distributive, e.g., the nonmodular pentagon lattice $\boldsymbol{N}_5$ is a pseudocomplemented semilattice.
\qed
\end{Remark}

We denote the class operators of closure under isomorphic copies, homomorphic images, subalgebras, direct products, and ultraproducts, respectively, by $\III, \HHH, \SSS$, $\PPP$, and $\PPU$. A class of similar algebras is said to be a \emph{quasivariety} when it is closed under $\III, \SSS, \PPP$, and $\PPU$ or, equivalently, when it can be axiomatized by a set of \emph{quasiequations} (see, e.g., \cite[{Thm. V.2.25}]{BuSa00}), i.e., universal sentences of the form
\[
\forall \vec{x}((\varphi_1 \thickapprox \psi_1 \, \& \, \dots \, \& \, \varphi_n \thickapprox \psi_n)\Longrightarrow \epsilon \thickapprox \delta).
\]
We often drop the the universal quantifier at the beginning of a quasiequation and write simply
\[
\varphi_1 \thickapprox \psi_1 \, \& \, \dots \, \& \, \varphi_n \thickapprox \psi_n \Longrightarrow \epsilon \thickapprox \delta.
\]

Lastly, a class of similar algebras is said to be a \emph{variety} when it is closed under $\HHH$, $\SSS$, and $\PPP$ or, equivalently, when it can be axiomatized by a set of equations (see, e.g., \cite[Thm.\ II.11.9]{BuSa00}). When understood as in Remark \ref{Rem:algebras}, the following classes are examples of varieties:
\begin{align*}
\mathsf{PSL} &\coloneqq \text{ the class of pseudocomplemented semilattices};\\
\mathsf{ISL} &\coloneqq \text{ the class of implicative semilattices};\\
\mathsf{bISL} &\coloneqq \text{ the class of bounded implicative semilattices};\\
\mathsf{PDL} &\coloneqq \text{ the class of pseudocomplemented distributive lattices};\\
\mathsf{IL} &\coloneqq \text{ the class of implicative lattices};\\
\mathsf{HA} &\coloneqq \text{ the class of Heyting algebras}.
\end{align*}

From a logical standpoint, the interest of these varieties can be explained as follows. Let $\mathscr{L}$ be a sublanguage of the language of an algebra $\A$. The $\mathscr{L}$\emph{-reduct} of $\A$ is the $\mathscr{L}$-algebra $\A_\mathscr{L} \coloneqq \langle A; \{ f^\A : f \in \mathscr{L}\}\rangle$. Accordingly, the subalgebras of the $\mathscr{L}$-reduct of $\A$ are called $\mathscr{L}$\emph{-subreducts} of $\A$. The varieties in the above display consist of the subreducts of Heyting algebras in the appropriate signature. For instance, $\mathsf{PSL}$ is the class of $\langle \land, \lnot, 0, 1 \rangle$-subreducts of Heyting algebras, and similarly for the other cases.

A class of similar algebras $\mathsf{K}$ is said to be a \emph{universal class} if it can be axiomatized by a set of universal sentences or, equivalently, if it is closed under $\III, \SSS$, and $\PPU$ (see, e.g., \cite[Thm.\ V.2.10]{BuSa00}). The least universal class containing a class of similar algebras $\mathsf{K}$ coincides with $\III\SSS\PPU(\mathsf{K})$ and will be denoted by $\UUU(\mathsf{K})$. The rest of the section is devoted to proving the following:

\begin{Theorem}\label{Thm:general-embedding}
Let $\A$ be a semilattice in a variety between $\mathsf{(b)ISL}, \mathsf{PDL}$, $\mathsf{IL}$, and $\mathsf{HA}$. Then $\A$ embeds into the appropriate reduct $\B^-$ of a Heyting algebra $\B$ such that $\B^- \in \UUU(\A)$.
\end{Theorem}

To prove this, recall that a variety $\mathsf{K}$ is \emph{locally finite} if its finitely generated members are finite. We rely on the following observation:

\begin{Proposition}\label{Prop:locally-finite}
The varieties $\mathsf{PSL}, \mathsf{(b)ISL}$, and $\mathsf{PDL}$ are locally finite.
\end{Proposition}

\begin{proof}
For $\mathsf{(b)ISL}$ and $\mathsf{PSL}$ the result is essentially \cite[Cor.\ III.4.1]{Di65} (see also \cite{Di66,Jon-PhD}), while for $\mathsf{PDL}$ see, e.g., \cite[Cor.\ 4.55]{Be11g} and \cite{Le70}.
\end{proof}

Furthermore, we will make use of the following general embedding theorem.

\begin{Proposition}[\protect{\cite[Thm.\ V.2.14]{BuSa00}}]\label{Prop:embedding}
Let $\mathsf{K} \cup \{ \A \}$ be a class of similar algebras. If every finitely generated subalgebra of $\A$ belongs to $\III\SSS(\mathsf{K})$, then $\A \in \III\SSS\PPU(\mathsf{K})$.
\end{Proposition}

We are now ready to prove Theorem \ref{Thm:general-embedding}.

\begin{proof}
We begin by the case where $\A \in \mathsf{ISL}$. Since every finitely generated subalgebra of $\A$ embeds into itself, in view of Proposition \ref{Prop:embedding} there exist a family $\{ \A_i : i \in I \}$ of finitely generated subalgebras of $\A$ and an ultrafilter $U$ on $I$ with an embedding
\[
f \colon \A \to \prod_{i \in I}\A_i /U.
\]
Since universal classes are closed under $\SSS$ and $\PPU$, we have
\[
\prod_{i \in I}\A_i /U \in \PPU\SSS(\A) \subseteq \UUU(\A).
\]

In view of Proposition \ref{Prop:locally-finite}, each $\A_i$ is finite and, therefore, can be expanded to a Heyting algebra $\A_i^+$ by Proposition \ref{Prop:finite-BSL}(\ref{item:finite-BSL2}). Clearly,
\[
\B \coloneqq \prod_{i \in I}\A_i^+ / U
\]
is a Heyting algebra, whose $\langle \land, \to \rangle$-reduct $\B^-$ coincides with $\prod_{i \in I}\A_i / U$. Since $\A$ embeds into $\B^-$ and $\B^- \in \UUU(\A)$, we are done.

The same proof works for the case where $\A$ belongs to $\mathsf{bISL}$ or to $\mathsf{PDL}$ with the only difference that, when $\A \in \mathsf{PDL}$, we apply Condition (\ref{item:finite-BSL1}) of Proposition \ref{Prop:finite-BSL} instead of Condition (\ref{item:finite-BSL2}). Lastly, the case where $\A$ belongs to $\mathsf{HA}$ is straightforward, since can simply take $\B \coloneqq \A$.

It only remains to consider the case where $\A \in \mathsf{IL}$. Let $\A^+$ be the expansion of $\A$ with a new constant $c_a$ for each of its elements $a$. Let also $0$ be another new constant. Then, consider the set of sentences
	\[
	\Sigma \coloneqq \{ 0 \leq c_a \colon a \in A\} \cup \mathsf{Th}(\boldsymbol{A}^+),
	\]
	where $\mathsf{Th}(\A^+)$ is the elementary theory of $\A^+$ and $0 \leq c_a$ is a shorthand for $0 \land c_a \thickapprox 0$. Clearly, every finite part $\Gamma$ of $\Sigma$ is realizable in $\A^+$, for if $c_{a_1}, \dots, c_{a_n}, 0$ are the new constants appearing in $\Gamma$, we can interpret $0$ as $c_{a_1} \land \dots \land c_{a_n}$ in $\A^+$. Therefore, we can apply the Compactness Theorem of first order logic, obtaining that $\Sigma$ has a model $\boldsymbol{C}$.

Let $\B^+$ be the subalgebra of $\boldsymbol{C}$ generated by $\{ 0 \} \cup \{ c_a : a \in A \}$. As $\boldsymbol{C}$ is a model of $\mathsf{Th}(\A^+)$, the map $f \colon A \to C$ defined by the assignment $f(a) \coloneqq c_a$ is an elementary embedding of $\A$ into the $\langle \land, \lor, \to, 1 \rangle$-reduct $\C^-$ of $\boldsymbol{C}$. Consequently, $\A$ embeds also into the $\langle \land, \lor, \to, 1 \rangle$-reduct $\B^-$ of $\B^+$. Furthermore, as $\B^-$ embeds into the elementary extension $\C^-$ of $\A$ and the validity of universal sentences persists in elementary extensions and subalgebras, $\B^-$ satisfies all the universal sentences valid in $\A$ and, therefore, belongs to $\UUU(\A)$.

To conclude the proof, it only remains to show that $\B^-$ is the $\langle \land, \lor, \to, 1 \rangle$-reduct of a Heyting algebra $\B$. Since $\A \in \mathsf{IL}$ and $\mathsf{IL}$ is a universal class, from $\B^-\in \UUU(\A)$ it follows that $\B^- \in \mathsf{IL}$. \color{black} Therefore, it suffices to prove that $0$ is the minimum element of $\B^-$, as in this case we can let $\B$ be the expansion of $\B^-$ with $0$.

To prove this, recall that $\B^+$ is the subalgebra of $\boldsymbol{C}$ generated by of $\{ 0 \} \cup \{ c_a : a \in A \}$. Therefore, every element of $\B^-$ (equiv.\ of $\B^+$) \color{black} has the form $\varphi^{\C}(0, c_{a_1}, \dots, c_{a_n})$ for some $a_1, \dots, a_n \in A$ and $\langle \land, \lor, \to, 1 \rangle$-term $\varphi(y, x_1, \dots, x_n)$. We will prove by induction on the construction of $\varphi$ that
\[
0 \leq \varphi^{\C}(0, c_{a_1}, \dots, c_{a_n}),
\]
for every $a_1, \dots, a_n \in A$.

In the base case, $\varphi$ is either the constant $1$ or a variable. If it is the constant $1$, then $\varphi^\C(0, c_{a_1}, \dots, c_{a_n})$ is the maximum of $\C$, whence the above display holds. Then we consider the case where $\varphi$ is a variable. Since $\varphi \in \{ y, x_1, \dots, x_n \}$, we have $\varphi^\C(0, c_{a_1}, \dots, c_{a_n}) \in \{ 0, c_{a_1}, \dots, c_{a_n}\}$. If $\varphi^\C(0, c_{a_1}, \dots, c_{a_n}) = 0$, it is clear that the above display holds. Consider the case where $\varphi^\C(0, c_{a_1}, \dots, c_{a_n}) = c_{a_i}$ for some $i \leq n$. Since $\C$ is a model of the formula $0 \leq c_{a_i}$ (which belongs to $\Sigma$), we obtain $0 \leq c_{a_i} = \varphi^\C(0, c_{a_1}, \dots, c_{a_n})$ as desired.

In the induction step, $\varphi$ is a complex formula. The case where $\varphi$ is of the form $\psi_1 \land \psi_2$ or $\psi_1 \lor \psi_2$ is straightforward. Accordingly, we detail only the case where $\varphi$ is of the form $\psi_1 \to \psi_2$. By the inductive hypothesis, we have
\[
0 \leq \psi_2^{\C}(0, c_{a_1}, \dots, c_{a_n})
\]
and, therefore,
\[
0 \land \psi_1^{\C}(0, c_{a_1}, \dots, c_{a_n}) \leq \psi_2^{\C}(0, c_{a_1}, \dots, c_{a_n}).
\]
Since $\C$ is an elementary extension of $\A$, it satisfies Condition (\ref{Eq:res-law}). Consequently, from the above display it follows
\[
0 \leq \psi_1^{\C}(0, c_{a_1}, \dots, c_{a_n}) \to^{\C} \psi_2^{\C}(0, c_{a_1}, \dots, c_{a_n}) = \varphi^{\C}(0, c_{a_1}, \dots, c_{a_n}).
\]
Hence, we conclude that $0$ is the minimum of $\B^-$ as desired.
\end{proof}

\section{Posets and partial functions}\label{Sec:adjunction}

\color{black} In this section, we will individuate a correspondence between homomorphisms in the varieties $\mathsf{PSL}$, $\mathsf{(b)ISL}$, $\mathsf{PDL}$, $\mathsf{IL}$, and $\mathsf{HA}$ and appropriate partial functions between (possibly empty) posets that generalize the notion of a p-morphism typical of \emph{Esakia duality} for Heyting algebras \cite{Es74,Esakia-book85}.

\color{black} The idea of using partial functions to dualize varieties of subreducts of Heyting algebras can be traced back at least to \cite{Kh81,Vrancken86} and \cite{Zakha89,Za92,Za96} and was developed systematically in \cite{BeBe09,BeJa13acs,Ce03,CeMo12}.\ Our presentation is largely inspired by the approach of \cite{BeBe09}, which deals with categories of Heyting algebras with maps preserving the operations in some smaller signature.\ Since we work with semilattices (as opposed to Heyting algebras), some additional care will be needed, however.

For a poset $\mathbb{X}$ and $Y \subseteq X$, let
\begin{align*}
{\uparrow} Y &\coloneqq \{x\in X :  \text{there exists } y\in Y \text{ s.t. } y\leq x\}; \\
{\downarrow} Y &\coloneqq \{x\in X :  \text{there exists } y\in Y \text{ s.t. } x\leq y\}.
\end{align*}
We call the set $Y$ an {\em upset} if $Y={\uparrow}Y$ and a {\em downset} if $Y={\downarrow}Y$. If $Y = \{ y \}$, we simply write ${\uparrow} y$ and ${\downarrow} y$ instead of ${\uparrow} \{ y \}$ and ${\downarrow} \{ y\}$. When the poset $\mathbb{X}$ is not clear from the context, we will write ${\uparrow}^{\mathbb{X}} Y$ and ${\downarrow}^{\mathbb{X}} Y$, respectively, instead of ${\uparrow} Y$ and ${\downarrow} Y$.

A \emph{partial function} $p$ from a set $X$ to a set $Y$ is a function from a subset $Z$ of $X$ to $Y$. In this case, $Z$ is said to be the \emph{domain} of $p$ and will be denoted by $\mathsf{dom}(p)$. We will write $p \colon X \rightharpoonup Y$ to indicate that $p$ is a partial function from $X$ to $Y$. A partial function $p \colon \mathbb{X} \rightharpoonup \mathbb{Y}$ between posets is \emph{order preserving} when, for every $x, z \in \mathsf{dom}(p)$,
\[
\text{if }x \leq^{\mathbb{X}} z, \text{ then }p(x) \leq^{\mathbb{Y}} p(z).
\]

\begin{Definition}\label{Def:partial-maps}
An order preserving partial function $p \colon \mathbb{X}  \rightharpoonup \mathbb{Y}$ between posets is
\benroman
\item\label{item:weak-partial-p-morphism} a \emph{partial negative p-morphism} if
\[
X = {\downarrow}^{\mathbb{X}} \{ x \in X : {\uparrow}^{\mathbb{X}}x \subseteq \mathsf{dom}(p) \}
\]
and for every $x \in \mathsf{dom}(p)$ and $y \in Y$,
\[
\text{if }p(x) \leq^{\mathbb{Y}} y, \text{ there exists }z \in \mathsf{dom}(p)\text{ s.t. }x \leq^{\mathbb{X}}z \text{ and }y \leq^{\mathbb{Y}} p(z);
\]
\item\label{item:implicative-partial-p-morphism} a \emph{partial positive p-morphism} if for every $x \in \mathsf{dom}(p)$ and $y \in Y$,
\[
\text{if }p(x) \leq^{\mathbb{Y}} y, \text{ there exists }z \in \mathsf{dom}(p)\text{ s.t. }x \leq^{\mathbb{X}}z \text{ and }y = p(z);
\]
\item\label{item:partial-p-morphism} a \emph{partial p-morphism} if it is both a partial negative p-morphism and a partial positive p-morphism.
\eroman
When $p$ is a total function, we drop the adjective \emph{partial} in the above definitions.
\end{Definition}

\begin{Remark}
Partial p-morphisms $p \colon \mathbb{X}  \rightharpoonup \mathbb{Y}$ coincide with partial positive p-morphism such that $X = {\downarrow}\mathsf{dom}(p)$.
\qed
\end{Remark}

We say that a partial function $p \colon \mathbb{X}  \rightharpoonup \mathbb{Y}$ between posets is \emph{almost total} when $\mathsf{dom}(p)$ is a downset of $\mathbb{X}$. Notice that almost total partial functions $p \colon \mathbb{X} \rightharpoonup \mathbb{Y}$ such that $X = {\downarrow}\mathsf{dom}(p)$ are indeed total. In particular, almost total 
 partial (negative) p-morphisms are total. On the other hand, almost total partial implicative p-morphism need not be total in general, e.g., if $U$ is a proper downset of $\mathbb{X}$ and $y$ a maximal element of $\mathbb{Y}$, the partial function $p \colon \mathbb{X} \rightharpoonup \mathbb{Y}$ such that $\mathsf{dom}(p) = U$ and $p[U] = \{ y \}$ is an almost total partial positive morphism that fails to be total.
 
\color{black}With every variety $\mathsf{K}$ among $\mathsf{PSL}, \mathsf{(b)ISL}, \mathsf{PDL}$, $\mathsf{IL}$, and $\mathsf{HA}$ we associate a collection $\mathsf{K}^\partial$ consisting of the class of all posets with suitable partial functions between them as follows:\footnote{Notice that $\mathsf{K}^\partial$ need not be a category in general.}
\begin{align*}
\mathsf{PSL}^{\partial} &\coloneqq \text{ the collection of posets with partial negative p-morphisms;}\\
\mathsf{ISL}^{\partial} &\coloneqq \text{ the collection of posets with partial positive p-morphisms;}\\
\mathsf{bISL}^{\partial} &\coloneqq \text{ the collection of posets with partial p-morphisms;}\\
\mathsf{PDL}^{\partial} &\coloneqq \text{ the collection of posets with negative p-morphisms;}\\
\mathsf{IL}^{\partial} &\coloneqq \text{ the collection of posets with almost total partial positive p-morphisms;}\\
\mathsf{HA}^{\partial} &\coloneqq \text{ the collection of posets with p-morphisms.}
\end{align*}
We will refer to the partial functions in $\mathsf{K}^\partial$ as to the \textit{arrows} of $\mathsf{K}^\partial$.

Every variety $\mathsf{K}$ among $\mathsf{PSL}, \mathsf{(b)ISL}, \mathsf{PDL}$, $\mathsf{IL}$, and $\mathsf{HA}$ is related to $\mathsf{K}^\partial$ as follows. \color{black} An element $a$ of a semilattice $\A$ is said to be \emph{meet irreducible} if it is not the maximum of $\A$ and, for every $b, c \in A$,
\[
\text{if }a = b \land c \text{, then either }a = b \text{ or }a = c.
\]
A \emph{filter} of $\A$ is a nonempty upset of $\langle A; \leq \rangle$ closed under binary meets. When ordered under the inclusion relation, the set of filters of $\A$ forms a lattice in which meets are intersections. Accordingly, a filter $F$ of $\A$ is said to be \emph{meet irreducible} when it is meet irreducible in the lattice of filters of $\A$, i.e., when $F$ is proper and either $F = G$ or $F = H$ for every pair $G, H$ of filters of $\A$ such that $F = G \cap H$. The poset of meet irreducible filters of $\A$ will be denoted by $\A_\ast$.

\begin{Remark}\label{Rem:prime-filters}
When the poset underlying a semilattice $\A$ is a distributive lattice, a filter $F$ of $\A$ is meet irreducible iff it is \emph{prime}, i.e., iff $A \smallsetminus F$ is nonempty and closed under binary joins \cite[Thm.\ 12]{BiFr48}. \qed 
\end{Remark}
Let $\mathsf{K}$ be a variety among $\mathsf{PSL}, \mathsf{(b)ISL}, \mathsf{PDL}$, $\mathsf{IL}$, and $\mathsf{HA}$. Given $\A, \B \in \mathsf{K}$ and a homomorphism $f \colon \A \to \B$, let $f_\ast \colon \B_\ast \rightharpoonup \A_\ast$ be the partial function with
\[
\mathsf{dom}(f_\ast) \coloneqq \{ F \in \B_\ast : f^{-1}[F] \in \A_\ast \}
\]
defined as $f_\ast(F) \coloneqq f^{-1}[F]$ for every $F \in \mathsf{dom}(f_\ast)$.

Conversely, given a poset $\mathbb{X}$, let $\mathsf{Up}_{\mathsf{K}}(\mathbb{X})$ be the reduct in the language of $\mathsf{K}$ of the Heyting algebra
\[
\langle \mathsf{Up}(\mathbb{X}); \cap, \cup, \to, \emptyset, X \rangle,
\]
where $\mathsf{Up}(\mathbb{X})$ is the set of upsets of $\mathbb{X}$ and $\to$ is defined by
\[
U \to V \coloneqq X \smallsetminus {\downarrow} (U \smallsetminus V).
\]
Notice that $\mathsf{Up}_{\mathsf{K}}(\mathbb{X}) \in \mathsf{K}$, because $\mathsf{K}$ is the class of subreducts of Heyting algebras in the language of $\mathsf{K}$. Lastly, given an arrow  $p \colon \mathbb{X} \rightharpoonup \mathbb{Y}$ in $\mathsf{K}^\partial$, let $\mathsf{Up}_{\mathsf{K}}(p) \colon \mathsf{Up}_{\mathsf{K}}(\mathbb{Y}) \to \mathsf{Up}_{\mathsf{K}}(\mathbb{X})$ be the map defined for every $U \in \mathsf{Up}_{\mathsf{K}}(\mathbb{Y})$ as $\mathsf{Up}_{\mathsf{K}}(p)(U) \coloneqq X \smallsetminus {\downarrow}^{\mathbb{X}} p^{-1}[Y \smallsetminus U]$. When $\mathsf{K} = \mathsf{HA}$, we often drop the subscript $\mathsf{K}$ from $\mathsf{Up}_{\mathsf{K}}(\mathbb{X})$ and $\mathsf{Up}_{\mathsf{K}}(p)$.

\begin{Remark}\label{Rem:Esakia duality}
\color{black}In the case of $\mathsf{HA}$, the applications $(-)_{\ast}$ and $\mathsf{Up}(-)$ are the contravariant functors underlying \emph{Esakia duality} \cite{Es74,Esakia-book85}. \color{black}
\qed
\end{Remark}

\color{black}
The rest of the section is devoted to the proof of the following result.
\color{black}

\begin{Proposition}\label{Prop:weak-adjunction}
Let $\mathsf{K}$ be a variety among $\mathsf{PSL}, \mathsf{(b)ISL}, \mathsf{PDL}$, $\mathsf{IL}$, and $\mathsf{HA}$. The following conditions hold for every $\A, \B \in \mathsf{K}$ and every pair $\mathbb{X}, \mathbb{Y}$ of posets:
\benroman
\item\label{item:weak-adjunction1} If $f \colon \A \to \B$ is a homomorphism, then $f_\ast \colon \B_\ast \rightharpoonup \A_\ast$ is an arrow in $\mathsf{K}^\partial$;
\item\label{item:weak-adjunction2} If $p \colon \mathbb{X}\rightharpoonup \mathbb{Y}$ is an arrow in $\mathsf{K}^\partial$, then $\mathsf{Up}_{\mathsf{K}}(p) \colon \mathsf{Up}_{\mathsf{K}}(\mathbb{Y}) \to \mathsf{Up}_{\mathsf{K}}(\mathbb{X})$ is a homomorphism.
\eroman
Furthermore, if $f$ is injective (resp.\ $p$ is surjective), then $f_\ast$ is surjective (resp.\ $\mathsf{Up}_{\mathsf{K}}(p)$ is injective).
\end{Proposition}
 
When ordered under the inclusion relation, the set $\mathsf{Fi}(\A)$ of filters of a semilattice $\A$ with maximum forms a lattice $\langle \mathsf{Fi}(\A); \cap, + \rangle$, where the join operation $+$ is defined as
\[
F + G \coloneqq \{ a \in A : \text{there are }b \in F \text{ and }c \in G \text{ such that }b \land c \leq a \}.
\]
We rely on the next observation.

\begin{Lemma}\label{Lem:genera-Zorn}
Let $f \colon \A \to \B$ be a homomorphism between two semilattices, $F$ a filter of $\B$, and $G \in \A_\ast$. If $f^{-1}[F] \subseteq G$ and the filter of $\B$ generated by $F \cup f[G]$ is disjoint from $f[A \smallsetminus G]$, then there exists $H \in \B_\ast$ such that $F \subseteq H$ and $G = f^{-1}[H]$.
\end{Lemma}

\begin{proof}
Consider the poset $\mathbb{X}$ whose universe is
\[
\{ P : P\text{ is a filter of }\B \text{ such that }F \cup f[G] \subseteq P \text{ and }P \cap f[A \smallsetminus G] = \emptyset \}
\]
and whose order is the inclusion relation. By assumption, $\mathbb{X}$ is nonempty because it contains the filter of $\B$ generated by $F \cup f[G]$. Since $\mathbb{X}$ is closed under unions of chains, we can apply Zorn's Lemma obtaining that $\mathbb{X}$ has a maximal element $H$.

\begin{Claim} The filter $H$ is meet irreducible. 
\end{Claim}

\begin{proof}[Proof of the Claim]
Suppose the contrary, with a view to contradiction. By assumption, $G$ is a meet irreducible filter of $\A$ and, therefore, proper. Consequently, $A \smallsetminus G$ is nonempty, whence so is $f[A \smallsetminus G]$. Since $H$ is disjoint from $f[A \smallsetminus G]$, we conclude that $H$ is proper. Since $H$ is not meet irreducible, this means that there are two filters $H_1, H_2$ of $\B$ other than $H$ such that $H = H_1 \cap H_2$. From the maximality of $H$ in $\mathbb{X}$ it follows that neither $H_1$ nor $H_2$ is disjoint from $f[A \smallsetminus G]$. Therefore, there are $a, b\in A \smallsetminus G$ such that 
\begin{equation}\label{Eq:general-ZornI}
f(a) \in H_1 \, \, \text{ and } \, \, f(b) \in H_2.
\end{equation}

Now, from $a, b \in A \smallsetminus G$ it follows that $G$ is properly contained in $G + {\uparrow}^\A a$ and $G + {\uparrow}^\A b$. Since $G$ a meet irreducible filter of $\A$, this guarantees that
\[
G \subsetneq (G + {\uparrow}^\A a) \cap (G + {\uparrow}^\B a).
\]
Accordingly, there are $c \in G$ and $d \in A \smallsetminus G$ such that
\[
a \land^{\A} c \leq d \, \, \text{ and } \, \, b \land^{\A} c \leq d.
\]
Since $f$ is a homomorphism, we obtain
\begin{equation}\label{Eq:general-ZornII}
f(a) \land^\B f(c) \leq f(d) \, \, \text{ and } \, \, f(b) \land^\B f(c) \leq f(d).
\end{equation}
Furthermore, from $c \in  G$ and the assumption that $f[G] \subseteq H = H_1 \cap H_2$ it follows that $f(c) \in H_1 \cap H_2$. Together with Conditions (\ref{Eq:general-ZornI}) and (\ref{Eq:general-ZornII}) and the fact that $H_1, H_2$ are filters of $\B$, this implies that
\[
f(d) \in H_1 \cap H_2 = H,
\]
 a contradiction with the assumption that $d \in A \smallsetminus G$ and $H \cap f[A \smallsetminus G] = \emptyset$. Hence, we conclude that $H$ is meet irreducible.
\end{proof}

By the Claim, $H \in \B_\ast$. Furthermore, since $H \in \mathbb{X}$, we know that $F \subseteq H$. Therefore, it only remains to prove that $G = f^{-1}[H]$. The fact that $H \in \mathbb{X}$ guarantees that $f[G] \subseteq H$ and $H \cap f[A \smallsetminus G] = \emptyset$. From $f[G] \subseteq H$ it follows $G \subseteq f^{-1}[f[G]] \subseteq f^{-1}[H]$. To prove the other inclusion, consider $a \in f^{-1}[H]$. Then $f(a) \in H$. Since $H \cap f[A \smallsetminus G] = \emptyset$ and $a \in A$, this implies that $a \in G$ as desired.
\end{proof}

Lastly, we rely on the following technical lemma.

\begin{Lemma}\label{Lemma:weak-adjunction_meets_and_joins}
	Let $p \colon \mathbb{X} \rightharpoonup \mathbb{Y}$ be a partial function between posets and $U,V \in \up{\mathbb{Y}}$. Then
	\[
	{\downarrow}^{\mathbb{X}} p^{-1}[Y \smallsetminus (U \cap V)] = {\downarrow}^{\mathbb{X}} p^{-1}[Y \smallsetminus U] \cup {\downarrow}^{\mathbb{X}} p^{-1}[Y \smallsetminus V].
	\]
	Moreover, if $p$ is total and order preserving, it also holds
	\[
	{\downarrow}^{\mathbb{X}} p^{-1}[Y \smallsetminus (U \cup V)] = {\downarrow}^{\mathbb{X}} p^{-1}[Y \smallsetminus U] \cap {\downarrow}^{\mathbb{X}} p^{-1}[Y \smallsetminus V].
	\]
\end{Lemma}
\begin{proof}
The proof of the first part of the statement is straightforward. As for the second part, suppose that $p$ is total and order preserving. 
	The inclusion from left to right follows from the fact that the function ${\downarrow}^{\mathbb{X}}p^{-1}[-] \colon \powerset{Y} \to \powerset{X}$ is order preserving.
	To prove the other inclusion, consider $x \in {\downarrow}^{\mathbb{X}} p^{-1}[Y \smallsetminus U] \cap {\downarrow}^{\mathbb{X}} p^{-1}[Y \smallsetminus V]$. Then there are $u, v \in X$ such that $p(u) \in Y \smallsetminus U$, $p(v) \in Y \smallsetminus V$, and $x \leq^{\mathbb{X}}u, v$. Since $p$ is a total function, $x \in \mathsf{dom}(p)$. Furthermore, as $p$ is order preserving and $x \leq^{\mathbb{X}}u, v$, we obtain $p(x) \leq^{\mathbb{Y}}p(u), p(v)$. Together with the assumption that $U$ and $V$ are upsets and that $p(u) \notin U$ and $p(v) \notin V$, this yields $p(x) \in Y \smallsetminus (U \cup V)$. Hence, we conclude that $x \in {\downarrow}^{\mathbb{X}} p^{-1}[Y \smallsetminus (U \cup V)]$. 
\end{proof}

In order to prove Proposition \ref{Prop:weak-adjunction}, it will be convenient to consider the cases of $\mathsf{PSL}$ and $\mathsf{ISL}$ separately. We begin by the case of $\mathsf{PSL}$.

\begin{Lemma}\label{Lem:PSL-proper-filter}
Let $\A, \B \in \mathsf{PSL}$ with a homomorphism $f \colon \A \to \B$, let $F \in \mathsf{Fi}(\B)$, and let $G \in \mathsf{Fi}(\A)$ be maximal and proper. If $f^{-1}[F] \subseteq G$, the filter of $\B$ generated by $F \cup f[G]$ is disjoint from $f[A \smallsetminus G]$.
\end{Lemma}

\begin{proof}
Suppose, with a view to contradiction, that $f^{-1}[F] \subseteq G$ and that the filter of $\B$ generated by $F \cup f[G]$ is not disjoint from $f[A \smallsetminus G]$. Then there exist $a_1, \dots, a_n \in G$, $b_1, \dots, b_m \in F$, and $c \in A \smallsetminus G$ such that
\[
f(a_1) \land^\B \dots \land^\B f(a_n) \land^\B b_1 \land^\B \dots \land^\B b_m \leq f(c).
\]
Since $G$ is maximal and proper, $G + {\uparrow}^\A c = A$. Therefore, $c \land^\A d = 0^\A$ for some $d \in G$. By Condition (\ref{Eq:neg-res-law}), this yields $d \leq \lnot^\A c$ and, since $f$ is a homomorphism, $f(d) \leq \lnot^\B f(c)$. 

As $d \in G$, we may assume, without loss of generality, that $d \in \{ a_1, \dots, a_n \}$. Therefore, from $f(d) \leq \lnot^\B f(c)$ and the above display it follows
\[
f(a_1) \land^\B \dots \land^\B f(a_n) \land^\B b_1 \land^\B \dots \land^\B b_m \leq f(c) \land^\B \lnot^\B f(c) = 0^\B.
\]
Since $f$ is a homomorphism, we can apply Condition (\ref{Eq:neg-res-law}) obtaining
\[
b_1 \land^\B \dots \land^\B b_m \leq \lnot^\B f(a_1 \land^\A \dots \land^\A a_n) = f(\lnot^\A(a_1 \land^\A \dots \land^\A a_n)).
\]
As $F$ is a filter of $\B$ containing $b_1, \dots, b_m$, the above display guarantees that $f(\lnot^\A(a_1 \land^\A \dots \land^\A a_n)) \in F$. As a consequence, $\lnot^\A(a_1 \land^\A \dots \land^\A a_n) \in f^{-1}[F] \subseteq G$. But together with the fact that $a_1, \dots, a_n \in G$ and that $G$ is a filter of $\A$, this implies 
\[
0^\A = (a_1 \land^\A \dots \land^\A a_n) \land \lnot^\A(a_1 \land^\A \dots \land^\A a_n) \in G,
\]
a contradiction with the assumption that $G$ is proper. 
\end{proof}

\color{black}The homomorphisms in $\mathsf{PSL}$ and the arrows in $\mathsf{PSL}^\partial$ are related as follows.\color{black}

\begin{Proposition}\label{Prop:PSL-duality}
Let $\A, \B \in \mathsf{PSL}$ and let $\mathbb{X},\mathbb{Y}$ be posets. The following conditions hold:
\benroman
\item\label{item:duality-PSL-1} If $f \colon \A \to \B$ is a homomorphism, then $f_{\ast} \colon \B_\ast \rightharpoonup \A_\ast$ is a partial negative p-morphism;
\item\label{item:duality-PSL-2} If $p \colon \mathbb{X} \rightharpoonup \mathbb{Y}$ is a partial negative p-morphism, then $\mathsf{Up}_{\mathsf{PSL}}(p) \colon \mathsf{Up}_{\mathsf{PSL}}(\mathbb{Y}) \to \mathsf{Up}_{\mathsf{PSL}}(\mathbb{X})$ is a homomorphism.
\eroman
\end{Proposition}

\begin{proof}
(\ref{item:duality-PSL-1}): The definition of $f_\ast$ guarantees that $f_\ast \colon \B_\ast \rightharpoonup \A_\ast$ is a well-defined partial order preserving map. Therefore, it suffices to prove that 
\benormal
\item\label{item:duality-PSL-1a} $\B_\ast = {\downarrow}^{\B_\ast}\{ F \in \B_\ast : {\uparrow}^{\B_\ast}F \subseteq \mathsf{dom}(f_\ast)\}$;
\item\label{item:duality-PSL-1b} for every $F \in \mathsf{dom}(f_\ast)$ and $G \in \A_\ast$,
\[
\text{if }f_\ast(F)  \subseteq G, \text{ there exists }H \in \mathsf{dom}(f_\ast)\text{ s.t. }F \subseteq H \text{ and }G \subseteq f_\ast(H).
\]
\enormal

Notice that inverse images under $f$ of proper filters of $\B$ are proper filters of $\A$, because $f$ preserves binary meets and minimum elements. We will use this observation repeatedly.

To prove Condition \ref{item:duality-PSL-1a}, it suffices to establish the inclusion $\B_\ast \subseteq {\downarrow}^{\B_\ast}\{ F \in \B_\ast : {\uparrow}^{\B_\ast}F \subseteq \mathsf{dom}(f_\ast)\}$ as the other one is obvious. Accordingly, consider $F \in \B_\ast$. Since $F$ is a proper filter of $\B$, the set $f^{-1}[F]$ is a proper filter of $\A$. By Zorn's Lemma, we can extend it to a maximal proper filter $G$ of $\A$. Being maximal and proper, $G$ is meet irreducible and, therefore, it belongs to $\A_\ast$. Furthermore, since $f^{-1}[F] \subseteq G$, we can apply Lemma \ref{Lem:PSL-proper-filter} obtaining that the filter of $\B$ generated by $F \cup f[G]$ is disjoint from $f[A \smallsetminus G]$. By Lemma \ref{Lem:genera-Zorn}, there exists $H \in \B_\ast$ such that $F \subseteq H$ and $G = f^{-1}[H]$. Thus, $H \in \mathsf{dom}(f_\ast)$ and $f_\ast(H) = G$. To conclude the proof, it only remains to show that ${\uparrow}^{\B_\ast} H \subseteq \mathsf{dom}(f_\ast)$. To this end, consider $H^+ \in {\uparrow}^{\B_\ast} H$. Since $H^+$ is a proper filter of $\B$, the set $f^{-1}[H^+]$ is a proper filter of $\A$. Furthermore, $G = f^{-1}[H] \subseteq f^{-1}[H^+]$. Since $G$ is a maximal proper filter of $\A$ and $f^{-1}[H^+]$ is a proper filter of $\A$, this implies $f^{-1}[H^+] = G \in \A_\ast$. Hence, we conclude that $H^+ \in \mathsf{dom}(f_\ast)$ as desired.

Then we turn to prove condition \ref{item:duality-PSL-1b}.\ Consider $F \in \mathsf{dom}(f_\ast)$ and $G \in \A_\ast$ such that $f_\ast(F)  \subseteq G$, that is, $f^{-1}[F] \subseteq G$. We may assume, without loss of generality, that $G$ is maximal and proper (otherwise we use Zorn's Lemma to extend it to such a filter). Therefore, we can apply Lemma \ref{Lem:PSL-proper-filter} obtaining that the filter of $\B$ generated by $F \cup f[G]$ is disjoint from $f[A \smallsetminus G]$. By Lemma \ref{Lem:genera-Zorn}, there exists $H \in \B_\ast$ such that $F \subseteq H$ and $G = f^{-1}[H]$. We conclude that $H \in \mathsf{dom}(f_\ast)$ and $G = f_\ast(H)$.

(\ref{item:duality-PSL-2}): To see that $\mathsf{Up}_{\mathsf{PSL}}(p) \colon \mathsf{Up}_{\mathsf{PSL}}(\mathbb{Y}) \to \mathsf{Up}_{\mathsf{PSL}}(\mathbb{X})$ preserves binary meets, recall from Lemma \ref{Lemma:weak-adjunction_meets_and_joins} that 
\[
{\downarrow}^{\mathbb{X}} p^{-1}[Y \smallsetminus (U \cap V)] = {\downarrow}^{\mathbb{X}} p^{-1}[Y \smallsetminus U] \cup {\downarrow}^{\mathbb{X}} p^{-1}[Y \smallsetminus V],
\]
for every $U,V \in \up{\mathbb{Y}}$. Therefore, we obtain
\[
X \smallsetminus {\downarrow}^{\mathbb{X}} p^{-1}[Y \smallsetminus (U \cap V)] = ( X \smallsetminus {\downarrow}^{\mathbb{X}} p^{-1}[Y \smallsetminus U]) \cap (X \smallsetminus {\downarrow}^{\mathbb{X}} p^{-1}[Y \smallsetminus V]),
\]
that is, $\mathsf{Up}_{\mathsf{PSL}}(p)(U \cap V) = \mathsf{Up}_{\mathsf{PSL}}(p)(U) \cap \mathsf{Up}_{\mathsf{PSL}}(p)(V)$, as desired. 

We now detail the proof that $\mathsf{Up}_{\mathsf{PSL}}(p)$ preserves the operation $\lnot$ (this, in turn, guarantees that $\mathsf{Up}_{\mathsf{PSL}}(p)$ preserves also the constant $0$, because the latter is term-definable as $x \land \lnot x$). Consider $U \in \mathsf{Up}_{\mathsf{PSL}}(\mathbb{Y})$. We need to prove that
\[
\mathsf{Up}_{\mathsf{PSL}}(p) (\lnot^{\mathsf{Up}_{\mathsf{PSL}}(\mathbb{Y})} (U)) =   \lnot^{\mathsf{Up}_{\mathsf{PSL}}(\mathbb{X})}\mathsf{Up}_{\mathsf{PSL}}(p)(U).
\]
Using the definitions of $\mathsf{Up}_{\mathsf{PSL}}(p)$ and of the operation $\lnot$ in $\mathsf{Up}_{\mathsf{PSL}}(\mathbb{Y})$ and $\mathsf{Up}_{\mathsf{PSL}}(\mathbb{X})$, this amounts to
\[
	X \smallsetminus {\downarrow}^{\mathbb{X}}(p^{-1}[{\downarrow}^{\mathbb{Y}}U]) =  X \smallsetminus {\downarrow}^{\mathbb{X}}(X \smallsetminus {\downarrow}^{\mathbb{X}}(p^{-1}[Y \smallsetminus U])).
\]
Clearly, it suffices to show that the complements of the sets in the above display coincide, namely,
\begin{equation}\label{Eq:PSL-poset-to-alg}
	{\downarrow}^{\mathbb{X}}(p^{-1}[{\downarrow}^{\mathbb{Y}}U]) =  {\downarrow}^{\mathbb{X}}(X \smallsetminus {\downarrow}^{\mathbb{X}}(p^{-1}[Y \smallsetminus U])).
\end{equation}

To prove the inclusion from left to right in Condition (\ref{Eq:PSL-poset-to-alg}), consider $x \in {\downarrow}^{\mathbb{X}}(p^{-1}[{\downarrow}^{\mathbb{Y}}U])$. Then there are $z \in X$ and $u \in U$ such that $z \in \mathsf{dom}(p)$ and 
\[
x \leq^{\mathbb{X}} z \, \, \text{ and } \, \, p(z) \leq^{\mathbb{Y}} u.
\]
Since $p(z) \leq^{\mathbb{Y}} u$ and $p \colon \mathbb{X} \rightharpoonup \mathbb{Y}$ is a partial negative p-morphism, there exists $w \in \mathsf{dom}(p)$ such that $z \leq^{\mathbb{X}} w$ and $u \leq^{\mathbb{Y}} p(w)$. Together with the above display, this yields $x \leq^{\mathbb{X}} w$. Therefore, to conclude that $x \in {\downarrow}^{\mathbb{X}}(X \smallsetminus {\downarrow}^{\mathbb{X}}(p^{-1}[Y \smallsetminus U]))$, it suffices to show that $w \in X \smallsetminus {\downarrow}^{\mathbb{X}}(p^{-1}[Y \smallsetminus U])$. Accordingly, consider $v \in \mathsf{dom}(p)$ such that $w \leq^{\mathbb{X}} v$. We need to show that $p(v) \in U$. Since $p$ is order preserving, from $w \leq^{\mathbb{X}} v$ it follows $p(w) \leq^{\mathbb{Y}} p(v)$. Together with $u \leq^{\mathbb{Y}} p(w)$, this yields $u \leq^{\mathbb{Y}}p(v)$. Since $U$ is an upset of $\mathbb{Y}$ and $u \in U$, we conclude that $p(v) \in U$ as desired.

Lastly, we prove the inclusion from right to left in Condition (\ref{Eq:PSL-poset-to-alg}). Consider $x \in {\downarrow}^{\mathbb{X}}(X \smallsetminus {\downarrow}^{\mathbb{X}}(p^{-1}[Y \smallsetminus U]))$. Then there exists $z \in X$ such that $x \leq^{\mathbb{X}} z$ and for every $w \in X$,
\[
\text{if }z \leq^{\mathbb{X}} w \text{ and }w \in \mathsf{dom}(p)\text{, then }p(w) \in U.
\]
Since $p$ is a partial negative p-morphism, $X = {\downarrow}^{\mathbb{X}}\mathsf{dom}(p)$. Thus, there exists $w \in \mathsf{dom}(p)$ with $z \leq^{\mathbb{X}} w$. In view of the above display, we obtain $p(w) \in U$. Since $x \leq^{\mathbb{X}} z\leq^{\mathbb{X}} w$, this yields $x \in {\downarrow}^{\mathbb{X}}(p^{-1}[U]) \subseteq {\downarrow}^{\mathbb{X}}(p^{-1}[{\downarrow}^{\mathbb{Y}}U])$.
\end{proof}

Now, we turn our attention to the case of implicative semilattices.

\begin{Lemma}\label{Lem:ISL-proper-filter}
Let $\A, \B \in \mathsf{ISL}$ with a homomorphism $f \colon \A \to \B$, let $F \in \mathsf{Fi}(\B)$, and let $G \in \mathsf{Fi}(\A)$. If $f^{-1}[F] \subseteq G$, the filter of $\B$ generated by $F \cup f[G]$ is disjoint from $f[A \smallsetminus G]$.
\end{Lemma}

\begin{proof}
Suppose, with a view to contradiction, that $f^{-1}[F] \subseteq G$ and that the filter of $\B$ generated by $F \cup f[G]$ contains an element $c \in f[A \smallsetminus G]$. Then there exist $a_1, \dots, a_n \in G$ and $b_1, \dots, b_m \in F$ such that
\[
f(a_1) \land^\B \dots \land^\B f(a_n) \land^\B b_1 \land^\B \dots \land^\B b_m \leq c.
\]
Since $f$ preserves binary meets, this yields
\[
f(a_1 \land^\A \dots \land^\A a_n) \land^\B b_1 \land^\B \dots \land^\B b_m \leq c.
\]
Together $c \in f[A \smallsetminus G]$, this implies that there exists $d \in A \smallsetminus G$ such that
\[
f(a_1 \land^\A \dots \land^\A a_n) \land^\B b_1 \land^\B \dots \land^\B b_m \leq f(d).
\] 
Applying Condition (\ref{Eq:res-law}) and the fact that $f$ preserves $\to$ to the above display, we obtain
\[
b_1 \land^\B \dots \land^\B b_m \leq  f(a_1 \land^\A \dots \land^\A a_n) \to^\B f(d) = f((a_1 \land^\A \dots \land^\A a_n) \to^\A d).
\]

Lastly, as $F$ is a filter of $\B$ containing $b_1, \dots, b_m$, the above display guarantees that $f((a_1 \land^\A \dots \land^\A a_n) \to^\A d) \in F$. As a consequence, $(a_1 \land^\A \dots \land^\A a_n) \to^\A d \in f^{-1}[F] \subseteq G$. Together with the fact that $a_1, \dots, a_n \in G$ and that $G$ is a filter of $\A$, this implies
\[
(a_1 \land^\A \dots \land^\A a_n) \land^\A ((a_1 \land^\A \dots \land^\A a_n) \to^\A d) \in G.
\]
Since $G$ is an upset and, by Condition (\ref{Eq:res-law}), we have
\[
(a_1 \land^\A \dots \land^\A a_n) \land^\A ((a_1 \land^\A \dots \land^\A a_n) \to^\A d) \leq d,
\]
this implies $d \in G$, a contradiction with the assumption that $d \in A \smallsetminus G$.
\end{proof}

\color{black} The homomorphisms in $\mathsf{ISL}$ and the arrows in $\mathsf{ISL}^\partial$ are related as follows.\color{black}

\begin{Proposition}\label{Prop:ISL-duality}
Let $\A, \B \in \mathsf{ISL}$ and let $\mathbb{X},\mathbb{Y}$ be posets. The following conditions hold:
\benroman
\item\label{item:duality-ISL-1} If $f \colon \A \to \B$ is a homomorphism, then $f_{\ast} \colon \B_\ast \rightharpoonup \A_\ast$ is a partial positive p-morphism;
\item\label{item:duality-ISL-2} If $p \colon \mathbb{X} \rightharpoonup \mathbb{Y}$ is a partial positive p-morphism, then $\mathsf{Up}_{\mathsf{ISL}}(p) \colon \mathsf{Up}_{\mathsf{ISL}}(\mathbb{Y}) \to \mathsf{Up}_{\mathsf{ISL}}(\mathbb{X})$ is a homomorphism.
\eroman
\end{Proposition}

\begin{proof}
(\ref{item:duality-ISL-1}): The definition of $f_\ast$ guarantees that $f_\ast \colon \B_\ast \rightharpoonup \A_\ast$ is a well-defined partial order preserving map. Therefore, it suffices to prove that for every $F \in \mathsf{dom}(f_\ast)$ and $G \in \A_\ast$,
\[
\text{if }f_\ast(F) \subseteq G, \text{ there exists }H \in \mathsf{dom}(f_\ast)\text{ s.t. }F \subseteq H  \text{ and }G = f_\ast(H).
\]

Accordingly, let $F \in \mathsf{dom}(f_\ast)$ and $G \in \A_\ast$ be such that $f_\ast(F) \subseteq G$, that is, $f^{-1}[F] \subseteq G$. By Lemma \ref{Lem:ISL-proper-filter}, the filter of $\B$ generated by $F \cup f[G]$ is disjoint from $f[A \smallsetminus G]$. Therefore, we can apply Lemma \ref{Lem:genera-Zorn} obtaining an $H \in \B_\ast$ that contains $F$ such that $G = f^{-1}[H]$. Since $G \in \A_\ast$, we conclude that $H \in \mathsf{dom}(f_\ast)$ and $G = f_{\ast}(H)$ as desired.

(\ref{item:duality-ISL-2}): The proof of this condition coincides with that of \cite[Thm.\ 3.15]{BeBe09} (although the respective statements are slightly different).
\end{proof}

We are now ready to prove Proposition \ref{Prop:weak-adjunction}.

\begin{proof}
The cases where $\mathsf{K}$ is $\mathsf{PSL}$, $\mathsf{ISL}$, or $\mathsf{bISL}$ follow from Propositions \ref{Prop:PSL-duality} and \ref{Prop:ISL-duality}, while the case where $\mathsf{K} = \mathsf{HA}$ is well known (see Remark \ref{Rem:Esakia duality}). Therefore it only remains to detail the cases of $\mathsf{PDL}$ and $\mathsf{IL}$. We detail the case of $\mathsf{PDL}$ only, as that of $\mathsf{IL}$ is analogous. 

To prove Condition (\ref{item:weak-adjunction1}), consider $\A, \B \in \mathsf{PDL}$ and let $f \colon \A \to \B$ be a homomorphism. Since $f$ is a homomorphism of pseudocomplemented semilattices, $f_\ast \colon \B_\ast \rightharpoonup \A_\ast$ is a partial negative p-morphism by Proposition \ref{Prop:PSL-duality}(\ref{item:duality-PSL-1}). To prove that $f_\ast$ is total, consider $F \in \B_\ast$. Since $\B$ is a distributive lattice, $F$ is a prime filter in view of Remark \ref{Rem:prime-filters}. As $f$ preserves binary joins, $f^{-1}[F]$ is a prime filter of $\A$, whence $f^{-1}[F] \in \A_\ast$ by Remark \ref{Rem:prime-filters}. Thus, we conclude that $F \in \mathsf{dom}(f_\ast)$ and, therefore, that $f_\ast$ is total. This shows that $f_\ast$ is a negative p-morphism. 

To prove Condition (\ref{item:weak-adjunction2}), let $p \colon \mathbb{X} \rightharpoonup \mathbb{Y}$ be a negative p-morphism. By Proposition \ref{Prop:PSL-duality}(\ref{item:duality-PSL-2}), $\mathsf{Up}_{\mathsf{PDL}}(p)$ is a homomorphism of pseudocomplemented semilattices.
Therefore, it only remains to prove that it preserves binary joins. To this end, consider $U, V \in \mathsf{Up}(\mathbb{Y})$. Since
\begin{align*}
	\mathsf{Up}_{\mathsf{PDL}}(p)(U \lor^{\mathsf{Up}_{\mathsf{PDL}}(\mathbb{Y})} V)&= X \smallsetminus {\downarrow}^{\mathbb{X}} p^{-1}[Y \smallsetminus (U \cup V)]\\
	\mathsf{Up}_{\mathsf{PDL}}(p)(U) \lor^{\mathsf{Up}_{\mathsf{PDL}}(\mathbb{X})} \mathsf{Up}_{\mathsf{PDL}}(p)(V)&= (X \smallsetminus {\downarrow}^{\mathbb{X}} p^{-1}[Y \smallsetminus U]) \cup (X \smallsetminus {\downarrow}^{\mathbb{X}} p^{-1}[Y \smallsetminus V]),
\end{align*}
it suffices to show that
\[
(X \smallsetminus {\downarrow}^{\mathbb{X}} p^{-1}[Y \smallsetminus U]) \cup (X \smallsetminus {\downarrow}^{\mathbb{X}} p^{-1}[Y \smallsetminus V]) = X \smallsetminus {\downarrow}^{\mathbb{X}} p^{-1}[Y \smallsetminus (U \cup V)].
\]
As $p$ is order preserving and total, we can apply the second part of Lemma \ref{Lemma:weak-adjunction_meets_and_joins} obtaining that the above display holds. As a consequence, $\mathsf{Up}_{\mathsf{PDL}}(p)$ is a homomorphism as desired.

It only remains to prove the last part of the statement, namely, that for every variety $\mathsf{K}$ among $\mathsf{PSL}, \mathsf{(b)ISL}, \mathsf{PDL}$, $\mathsf{IL}$, and $\mathsf{HA}$, every homomorphism $f \colon \A \to \B$ in $\mathsf{K}$, and every arrow $p \colon \mathbb{X} \rightharpoonup \mathbb{Y}$ in $\mathsf{K}^\partial$, if $f$ is injective (resp.\ $p$ is surjective), then $f_\ast$ is surjective (resp.\ $\mathsf{Up}_{\mathsf{K}}(p)$ is injective).

Suppose first that $f \colon \A \to \B$ is injective and consider $G \in \A_\ast$. Since $G$ is proper, we can choose an element $a \in G$ and consider the filter $F \coloneqq {\uparrow}^\B f(a)$ of $\B$. Since $f$ is order reflecting, $f^{-1}[F] \subseteq {\uparrow}^\A a$. Together with $a \in G$ and the fact that $G$ is an upset, this implies $f^{-1}[F] \subseteq G$. We will prove that the filter of $\B$ generated by $F \cup f[G]$ is disjoint from $f[A \smallsetminus G]$. Suppose, on the contrary, that there exists some $b \in B$ in the filter of $\B$ generated by $F \cup f[G]$ and in $f[A \smallsetminus G]$. Since $b \in f[A \smallsetminus G]$ there exists $c \in A \smallsetminus G$ such that $f(c) = b$. Since $F = {\uparrow}^\B f(a)$ and $b$ belongs to the filter of $\B$ generated by $F \cup f[G]$, there are $a_1, \dots, a_n \in G$ such that
\[
f(a \land^\A a_1 \land^\A \dots \land^\A a_n )= f(a) \land^\B f(a_1) \land^\B \dots \land^{\B} f(a_n) \leq b = f(c).
\]
As $f$ is order reflecting, this implies $a \land^\A a_1 \land^\A \dots \land^\A a_n \leq a$. As $G$ is a filter and $a, a_1, \dots, a_n \in G$, we obtain that $c \in G$, a contradiction with the assumption that $c \in A \smallsetminus G$. In sum, $f^{-1}[F] \subseteq G$ and the filter of $\B$ generated by $F \cup f[G]$ is disjoint from $f[A \smallsetminus G]$. Therefore, we can apply Lemma \ref{Lem:genera-Zorn} obtaining an $H \in \B_\ast$ such that $G = f^{-1}[H]$. Hence, $H \in \mathsf{dom}(f_\ast)$ and $f_\ast(H) = G$ and we conclude that $f_\ast$ is surjective.

Lastly, consider a surjective arrow $p \colon \mathbb{X} \rightharpoonup \mathbb{Y}$ in $\mathsf{K}^\partial$ and let $U, V$ be distinct upsets of $\mathbb{Y}$. By symmetry, we may assume that there exists $y \in U \smallsetminus V$. Since $p$ is surjective, there exists $x \in \mathsf{dom}(p)$ such that $p(x) = y$. Since $p$ is order preserving, $p(x) \in U$, and $U$ is an upset, we have $x \notin {\downarrow}^{\mathbb{X}}p^{-1}[Y \smallsetminus U]$, whence $x \in X \smallsetminus {\downarrow}^{\mathbb{X}}p^{-1}[Y \smallsetminus U] = \mathsf{Up}_{\mathsf{K}}(p)(U)$. On the other hand, as $p(x) \in U \smallsetminus V \subseteq Y \smallsetminus V$, we obtain $x \in {\downarrow}^{\mathbb{X}}p^{-1}[Y \smallsetminus V]$, whence $x \notin X \smallsetminus {\downarrow}^{\mathbb{X}}p^{-1}[Y \smallsetminus V] = \mathsf{Up}_{\mathsf{K}}(p)(V)$. Hence, we conclude that $\mathsf{Up}_{\mathsf{K}}(p)(U) \ne \mathsf{Up}_{\mathsf{K}}(p)(V)$ and, therefore, that $\mathsf{Up}_{\mathsf{K}}(p)$ is injective as desired.
\end{proof}

\section{Sahlqvist theory for $\mathsf{IPC}$}

Sahlqvist theory \cite{Sa75} is usually formulated in the setting of modal logic (see, e.g., \cite{BlRiVe01,SaVa89}). However, the G\"odel-McKinsey-Tarski translation of the intuitionistic propositional calculus $\mathsf{IPC}$ into the modal system $\mathsf{S4}$ allows to extend Sahlqvist theory to $\mathsf{IPC}$, as shown in \cite{CoPaZa19}\footnote{For a different approach to the canonicity part of Sahlqvist theorem for $\mathsf{IPC}$, see \cite{GhiMe97}.}. In this section, we will review this process in such a way that it will help set the stage for the study of fragments of $\mathsf{IPC}$.

 Consider the modal language
	\begin{align*}	
		\mathcal{L}_{\Box}\quad \!\!\!  ::=    & \quad x \mid  \varphi \land \psi  \mid  \varphi \lor \psi  \mid  \varphi \to \psi  \mid  \lnot \varphi  \mid \Box \varphi \mid \Diamond \varphi \mid  0  \mid  1.
	\end{align*}
Formulas of $\mathcal{L}_\Box$ will be assumed to have variables in a denumerable set $Var = \{ x_n : n \in \mathbb{Z}^+ \}$ and arbitrary elements of $Var$ will often be denoted by $x, y$, and $z$.

	\begin{Definition}
		Let $\varphi$ be a formula of $\mathcal{L}_\Box$ and $x$ a variable. An occurrence of $x$ in $\varphi$ is said to be \emph{positive} (resp.\ \emph{negative}) if the sum of negations and antecedents of implications within whose scopes it appears is even (resp.\ odd). Moreover, we say a $x$ is \emph{positive} (resp.\ \emph{negative}) in $\varphi$ if every occurrence of $x$ in $\varphi$ is positive (resp.\ negative). Lastly, $\varphi$ is said to be \emph{positive} (resp.\ \emph{negative}) if every variable is positive (resp.\ negative) in $\varphi$.
	\end{Definition}
	
Formulas of the form $\Box^n x$ with $x \in Var$ and $n \in \mathbb{N}$ will be called \emph{boxed atoms}. Notice that the elements of $Var$ are also boxed atoms, because $x = \Box^0 x$ for every $x \in Var$.
	
	\begin{Definition}\label{def : sahlqvist formulas}
		A formula of $\mathcal{L}_\Box$ is said to be
		\benroman
			\item a \emph{modal Sahlqvist antecedent} if it is constructed from boxed atoms, negative formulas, and the constants $0$ and $1$ using only $\land$, $\lor$, and $\Diamond$;
			\item a \emph{modal Sahlqvist implication} if it is positive, or it is of the form $\lnot \varphi$ for a modal Sahlqvist antecedent $\varphi$, or it is of the form $\varphi \to \psi$ for a modal Sahlqvist antecedent $\varphi$ and a positive formula $\psi$.
		\eroman
	\end{Definition}
\begin{Remark}
When applied to modal logic, our definition of a modal Sahlqvist implication is intentionally redundant. For if $\varphi$ is positive and $\psi$ a modal Sahlqvist antecedent, then $\varphi$ is equivalent to $1 \to \varphi$ and $\lnot \psi$ is equivalent to $\psi \to 0$. Accordingly, in modal logic the third possibility in the definition of a modal Sahlqvist implication subsumes (up lo logical equivalence) the first two.
\qed
\end{Remark}
In the next definition $x \leq y$ is a shorthand for the equation $x \land y \thickapprox x$.

\begin{Definition}
A \emph{modal Sahlqvist quasiequation} is an expression $\Phi$ of the form
\[
\varphi_1 \land y \leq z \, \& \dots \& \, \varphi_n \land y \leq z \Longrightarrow y \leq z,
\]
where $y$ and $z$ are distinct variables that do not occur in $\varphi_1, \dots, \varphi_n$ and each $\varphi_i$ is constructed from modal Sahlqvist implications using only $\land, \lor$, and $\Box$. If, in addition, $\Phi$ does not contain any occurrence of $\Box$ or $\Diamond$, we say that $\Phi$ is simply a \emph{Sahlqvist quasiequation}.
\end{Definition}

\begin{exa}\label{Exa:BTWn}
For every $n \in \mathbb{Z}^+$, the \emph{bounded top width $n$} axiom\footnote{The formulation of $\mathsf{btw}_n$ given in \cite{Smory73} is equivalent to the one we employ, in the sense that the two formulas axiomatize the same axiomatic extension of $\mathsf{IPC}$ (for a proof, see, e.g., \cite{FornasierePhD}). Our formulation has the advantage of making the connection with Sahlqvist quasiequations apparent.} \cite{Smory73} is the formula of $\mathsf{IPC}$
\[
\mathsf{btw}_n \coloneqq \bigvee_{i = 1}^{n+1} \lnot (\lnot x_i \land \bigwedge_{0 <j< i}x_j).
\]
When $n = 1$, the formula $\mathsf{btw}_n$ is equivalent over $\mathsf{IPC}$ to the \emph{weak excluded middle law} $\lnot x \lor \lnot \lnot x$ \cite{Jankov68a}. Notably, each $\mathsf{btw}_n$ can be rendered as the Sahlqvist quasiequation
\[
\Phi_n = \foo_{1 \leq i \leq n+1}\Big(\lnot(\lnot x_i \land \bigwedge_{0 < j< i}x_j) \land y \leq z\Big) \Longrightarrow y \leq z
\]
in the sense that a Heyting algebra validates $\mathsf{btw}_n$ iff it validates $\Phi_n$.

Similarly, the \emph{excluded middle} $x \lor \lnot x$ and the \emph{G\"odel-Dummett} axiom $(x_1 \to x_2) \lor (x_2 \to x_1)$ \cite{Dm59,Go32ae} can be rendered, respectively, as the Sahlqvist quasiequations
\[
\pushQED{\qed}x \land y \leq z \, \& \, \lnot x \land y \leq z \Longrightarrow y \leq z \, \, \text{ and } \, \, (x_1 \to x_2) \land y \leq z \, \& \, (x_2 \to x_1) \land y \leq z \Longrightarrow y \leq z.\qedhere \popQED
\]
\end{exa}

In the modal logic literature, the role of modal Sahlqvist quasiequations is played by the so-called \emph{modal Sahlqvist formulas}, i.e., formulas that can be constructed from modal Sahlqvist implications using only $\land, \lor$, and $\Box$.\footnote{It is common to define modal Sahlqvist formulas as the formulas that can be obtained from modal Sahlqvist implications using only $\land, \Box$, and disjunctions of formulas with no variable in common (see, e.g., \cite{BlRiVe01}), but our definition coincides (up to logical equivalence) with the standard one as shown in \cite[Rmk.\ 4.3]{BeBeHo12}.} When $\Box$ and $\Diamond$ do not occur in a modal Sahlqvist formula $\varphi$, we will say that $\varphi$ is simply a \emph{Sahlqvist formula}.

In order to clarify the relation between (modal) Sahlqvist quasiequations and formulas, recall that a \emph{modal algebra} is a structure $\langle A; \land, \lor, \lnot, \Box, 0, 1 \rangle$ where $\langle A; \land,\lor, \lnot, 0, 1 \rangle$ is a Boolean algebra and for every $a, b \in A$,
\[
\Box (a \land b) = \Box a \land \Box b \, \, \text{ and } \, \, \Box 1 = 1.
\]

We say that a formula $\varphi$ is \emph{valid} in a modal (resp.\ Heyting) algebra $\A$, in symbols $\A \vDash \varphi$, when $\A$ satisfies the equation $\varphi \thickapprox 1$. 

\begin{Proposition}\label{Prop:meaning-Sahlqvist-modal}
A modal Sahlqvist quasiequation
\[
\Phi = \varphi_1 \land y \leq z \, \& \dots \& \, \varphi_n \land y \leq z \Longrightarrow y \leq z
\]
is valid in a modal algebra $\A$ iff $\A \vDash \varphi_1 \lor \dots \lor \varphi_n$.
\end{Proposition}

\begin{proof}
Suppose that $\A \vDash \Phi$ and consider $\vec{a} \in A$. For every $i \leq n$, we have
\[
\varphi_i(\vec{a}) \land 1 = \varphi_i(\vec{a}) \leq\varphi_1(\vec{a}) \lor \dots \lor \varphi_n(\vec{a}).
\]
Since $\A \vDash \Phi$, this implies $1 \leq \varphi_1(\vec{a}) \lor \dots \lor \varphi_n(\vec{a})$. As $1$ is the maximum of $\A$, we conclude that $\varphi_1(\vec{a}) \lor \dots \lor \varphi_n(\vec{a}) = 1$ as desired.

Conversely, suppose that $\A \vDash \varphi_1 \lor \dots \lor \varphi_n$ and consider $\vec{a}, b, c \in A$ such that $\varphi_i(\vec{a}) \land b \leq c$ for every $i \leq n$. Using the distributive laws, we obtain
\[
(\varphi_1(\vec{a}) \lor \dots \lor \varphi_n(\vec{a})) \land b = (\varphi_1(\vec{a}) \land b) \lor \dots \lor(\varphi_n(\vec{a}) \land b) \leq c.
\]
Since $\A \vDash \varphi_1 \lor \dots \lor \varphi_n$, this yields $b = 1 \land b \leq c$, whence $\A \vDash \Phi$.
\end{proof}

A similar argument yields the following:

\begin{Corollary}\label{Cor:meaning-Sahlqvist-IPC}
A Sahlqvist quasiequation
\[
\varphi_1 \land y \leq z \, \& \dots \& \, \varphi_n \land y \leq z \Longrightarrow y \leq z
\]
 is valid in a Heyting algebra $\A$ iff $\A \vDash \varphi_1 \lor \dots \lor \varphi_n$.
\end{Corollary}

The reason why (modal) Sahlqvist quasiequations and formulas are two faces of the same coin is that, in view of Proposition \ref{Prop:meaning-Sahlqvist-modal} and Corollary \ref{Cor:meaning-Sahlqvist-IPC}, a (resp.\ modal) Sahlqvist quasiequation
\[
\varphi_1 \land y \leq z \, \& \dots \& \, \varphi_n \land y \leq z \Longrightarrow y \leq z
\]
is valid in a Heyting (resp.\ modal) algebra $\A$ iff so is the (resp.\ modal) Sahlqvist formula $\varphi_1 \lor \dots \lor \varphi_n$. Conversely, a (resp.\ modal) Sahlqvist formula $\varphi$ is valid in $\A$ iff so is the (resp.\ modal) Sahlqvist quasiequation $\varphi \land y \leq z \Longrightarrow y \leq z$. 

\begin{figure}
\[
	\begin{tikzpicture}
		\tikzstyle{point} = [shape=circle, thin, draw=black, fill=black, scale=0.35]
		
		\node (0) at (0,0)[point] {};
		\node (a1) at (-1,1)[point] {};
		\node (a2) at (1,1)[point] {};
		\node[label=left:{$a$}] (b1) at (-2,2)[point] {};
		\node[label=left:{$c$}] (b2) at (0,2)[point] {};
		\node (b3) at (2,2)[point] {};
		\node[label=left:{$b$}] (c) at (0,3)[point] {};
		\node (1) at (0,4)[point] {};
		
		\draw (0)--(a1)--(b1)--(1)--(b3)--(0);
		\draw (0)--(a1)--(b2)--(a2);
		\draw (b2)--(c)--(1);
	\end{tikzpicture}
\]
\caption{A pseudocomplemented semilattice.}
\label{Fig:counterexample-PSL}
\end{figure}
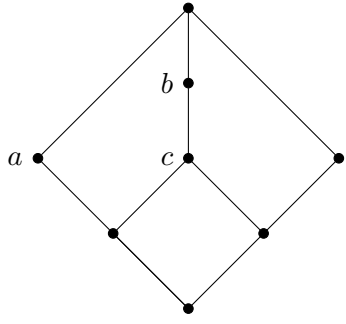

\begin{Remark}
The focus on Sahlqvist quasiequations (as opposed formulas or equations) is motivated by the fact that we deal with fragments of $\mathsf{IPC}$ where formulas have a very limited expressive power.\ For instance, in $\mathsf{PSL}$ there are only three nonequivalent equations \cite{Jon-PhD}, while there are infinitely many nonequivalent Sahlqvist quasiequations, as shown in Example \ref{Exa:BTWn}.

In addition, we cannot remove the ``context'' $y$ from Sahlqvist quasiequations. For instance, the Sahlqvist quasiequation
\[
\Phi = \lnot x \land y \leq z \, \& \, \lnot \lnot x \land y \leq z \Longrightarrow y \leq z
\]
corresponding to the weak excluded middle law (see Example \ref{Exa:BTWn}) is not equivalent to its context free version $\Psi = \lnot x \leq z \, \& \, \lnot \lnot x \leq z \Longrightarrow z \thickapprox 1$ over $\mathsf{PSL}$, for $\Psi$ holds in the pseudocomplemented semilattice $\A$ depicted in Figure \ref{Fig:counterexample-PSL}, while $\Phi$ fails in $\A$ as witnessed by the assignment
\[
\pushQED{\qed}x \longmapsto a \qquad y \longmapsto b \qquad z \longmapsto c.\qedhere \popQED
\]
\end{Remark}

With every Kripke frame $\mathbb{X} = \langle X, R \rangle$ we can associate a modal algebra
\[
\powersetm{\mathbb{X}} \coloneqq \langle \powerset{X}; \cap, \cup, \lnot, \Box, \emptyset, X \rangle,
\]
where $\lnot$ and $\Box$ are defined for every $Y \subseteq X$ as
\[
\lnot Y \coloneqq X \smallsetminus Y \, \, \text{ and } \, \, \Box Y \coloneqq \{ x \in X : \text{if }\langle x, y \rangle \in R \text{, then }y \in Y \}.
\]
Conversely, with a modal algebra $\A$ we can associate a Kripke frame $\A_+ \coloneqq \langle X, R \rangle$, where $X$ is the set of ultrafilters of $\A$ and
\[
R \coloneqq \{ \langle F, G \rangle \in X \times X : \text{for every $a \in A$, if }\Box^{\A} a \in F \text{, then }a \in G \}.
\]
Notably, $\A$ embeds into the algebra $\powersetm{\A_+}$, known as the \emph{canonical extension} of $\A$ \cite{JonTar51,JonTar52}.

Our aim is to extend the next classical version of Sahlqvist theorem to $\mathsf{IPC}$.

\begin{Modal Sahlqvist Theorem}[\protect{\cite[Thms.\ 3.54 and 5.91]{BlRiVe01}}]
The following conditions hold for a modal Sahlqvist quasiequation $\Phi$:
\benroman
\item \emph{Canonicity}: If a modal algebra $\A$ validates $\Phi$, then also $\powersetm{\A_+}$ validates $\Phi$;
\item \emph{Correspondence}: There is an effectively computable first order sentence $\mathsf{mtr}(\Phi)$\footnote{In the modal logic literature, $\mathsf{mtr}(\Phi)$ is the so-called \emph{standard translation} of the Sahlqvist formula $\varphi$ associated with $\Phi$. Furthermore, the demand that $\powersetm{\mathbb{X}} \vDash \Phi$ is equivalent to the requirement that $\varphi$ is valid in the Kripke frame $\mathbb{X}$.} in the language of Kripke frames such that $\powersetm{\mathbb{X}} \vDash \Phi$ iff $\mathbb{X} \vDash \mathsf{mtr}(\Phi)$, for every Kripke frame $\mathbb{X}$.
\eroman
\end{Modal Sahlqvist Theorem}

Let $\mathcal{L}$ be the language of $\mathsf{IPC}$, i.e., the language obtained from $\mathcal{L}_\Box$ by removing $\Box$ and $\Diamond$. Recall that the \emph{G\"odel-McKinsey-Tarski translation} \cite{Go32,McKT48} associates with every formula $\varphi$ of $\mathcal{L}$ a formula $\varphi_g$ of $\mathcal{L}_\Box$, defined recursively as follows: for every $x \in Var$,
\[
x_g \coloneqq \Box x \quad 0_g \coloneqq 0 \quad 1_g \coloneqq 1\quad (\varphi \land \psi)_g \coloneqq \varphi_g \land \psi_g
\]
\[
 (\varphi \lor \psi)_g \coloneqq \varphi_g \lor \psi_g\quad (\varphi \to \psi)_g \coloneqq \Box(\varphi_g \to \psi_g)\quad(\lnot \varphi)_g \coloneqq \Box \lnot\varphi_g.
\]

Given a Sahlqvist quasisequation
\[
\Phi = \varphi_1 \land y \leq z \, \& \dots \& \, \varphi_n \land y \leq z \Longrightarrow y \leq z,
\]
we set
\[
\Phi_g \coloneqq\varphi_{1g} \land y \leq z \, \& \dots \& \, \varphi_{ng} \land y \leq z \Longrightarrow y \leq z.
\]
The following observation is an immediate consequence of the definitions.

\begin{Lemma}\label{Lem:still-Sahlqvist}
If $\Phi$ is a Sahlqvist quasiequation, then $\Phi_g$ is a modal Sahlqvist quasiequation.
\end{Lemma}

The next result is instrumental to extend Sahlqvist theorem to $\mathsf{IPC}$.

\begin{Proposition}\label{Prop:free-extension}
The following conditions hold:
\benroman
\item\label{item:1:IPC-Sahlqvist} $\mathsf{Up}(\mathbb{X}) \vDash \Phi$ iff $\powersetm{\mathbb{X}} \vDash \Phi_g$, for every poset $\mathbb{X}$ and Sahlqvist quasiequation $\Phi$;
\item\label{item:2:IPC-Sahlqvist} For every Heyting algebra $\A$ there exists a modal algebra $\mathsf{f}(\A)$ with $\A_\ast \cong \mathsf{f}(\A)_+$ and such that $\A \vDash \Phi$ iff $\mathsf{f}(\A) \vDash \Phi_g$, for every Sahlqvist quasiequation $\Phi$.
\eroman
\end{Proposition}

\begin{proof}
(\ref{item:1:IPC-Sahlqvist}): It is well known that
\begin{equation}\label{Eq:Chag-Zakha}
\mathsf{Up}(\mathbb{X}) \vDash \varphi \iff \powersetm{\mathbb{X}} \vDash \varphi_g,
\end{equation}
 for every formula $\varphi$ of $\mathcal{L}$ and poset $\mathbb{X}$ (see, e.g., \cite[Cor.\ 3.82]{ChZa97}). Then for every poset $\mathbb{X}$ and Sahlqvist quasiequation $\Phi = \varphi_1 \land y \leq z \, \& \dots \& \, \varphi_n \land y \leq z \Longrightarrow y \leq z$, we have
\begin{align*}
\mathsf{Up}(\mathbb{X}) \vDash \Phi  \, \, &\Longleftrightarrow \, \, \mathsf{Up}(\mathbb{X}) \vDash \varphi_1 \lor \dots \lor \varphi_n\\
\, \, &\Longleftrightarrow \, \,\powersetm{\mathbb{X}} \vDash (\varphi_1 \lor \dots \lor \varphi_n)_g \\
\, \,&\Longleftrightarrow \, \,\powersetm{\mathbb{X}} \vDash \varphi_{1g} \lor \dots \lor \varphi_{ng}\\
\, \,&\Longleftrightarrow \, \,\powersetm{\mathbb{X}} \vDash \Phi_g.
\end{align*}
The equivalences above are justified as follows: the first and the last follow, respectively, from Corollary \ref{Cor:meaning-Sahlqvist-IPC} and Proposition \ref{Prop:meaning-Sahlqvist-modal}, the second holds by Condition (\ref{Eq:Chag-Zakha}), and the third by the definition of the G\"odel-McKinsey-Tarski translation.

(\ref{item:2:IPC-Sahlqvist}): Let $\mathsf{f}(\A)$ be the subalgebra of $\powersetm{\A_\ast}$ generated by the sets of the form
\[
\epsilon_{\A}(a) \coloneqq \{ F \in \A_\ast : a \in F \}
\]
for every $a \in A$. In view of \cite[Lem.\ 3.1 and 3.2]{MakRyb74}, for every formula $\varphi$ of $\mathcal{L}$ we have
\[
\A \vDash \varphi \, \,\Longleftrightarrow \, \,\mathsf{f}(\A)\vDash \varphi_g.
\]
As in the proof of Condition (\ref{item:1:IPC-Sahlqvist}), this implies that $\A \vDash \Phi$ iff $\mathsf{f}(\A) \vDash \Phi_g$, for every Sahlqvist quasiequation $\Phi$. For a proof that $\A_\ast \cong \mathsf{f}(\A)_+$, see, e.g., \cite[Construction 2.5.7 and Thm.\ 3.4.6(1)]{Esakia-book85}.
\end{proof}

As a consequence, we obtain a version of Sahlqvist theorem for $\mathsf{IPC}$:

\begin{Intuitionistic Sahlqvist Theorem}[\protect{\cite[Thms.\ 6.1 and 7.1]{CoPaZa19}}]\label{Thm:IPC-Sahlqvist}
The following conditions hold for a Sahlqvist quasiequation $\Phi$:
\benroman
\item\label{item:Sahlvist:IPC1} \emph{Canonicity}: If a Heyting algebra $\A$ validates $\Phi$, then also $\mathsf{Up}(\A_\ast)$ validates $\Phi$;
\item\label{item:Sahlvist:IPC2} \emph{Correspondence}: There is an effectively computable first order sentence $\mathsf{tr}(\Phi)$ in the language of posets such that $\mathsf{Up}(\mathbb{X}) \vDash \Phi$ iff $\mathbb{X} \vDash \mathsf{tr}(\Phi_g)$, for every poset $\mathbb{X}$.
\eroman
\end{Intuitionistic Sahlqvist Theorem}

\begin{proof}
(\ref{item:Sahlvist:IPC1}): Suppose that $\A \vDash \Phi$. In view of Proposition \ref{Prop:free-extension}(\ref{item:2:IPC-Sahlqvist}), we have $\mathsf{f}(\A) \vDash \Phi_g$. As $\Phi_g$ is a modal Sahlqvist quasiequation by Lemma \ref{Lem:still-Sahlqvist}, we can apply the canonicity part of the Modal Sahlqvist Theorem obtaining $\powersetm{\mathsf{f}(\A)_+} \vDash \Phi_g$. Since $\A_\ast \cong \mathsf{f}(\A)_+$ by Proposition \ref{Prop:free-extension}(\ref{item:2:IPC-Sahlqvist}), this amounts to $\powersetm{\A_\ast} \vDash \Phi_g$. Together with by Proposition \ref{Prop:free-extension}(\ref{item:1:IPC-Sahlqvist}), this implies that $\mathsf{Up}(\A_\ast) \vDash \Phi$ as desired.

(\ref{item:Sahlvist:IPC2}): From Proposition \ref{Prop:free-extension}(\ref{item:1:IPC-Sahlqvist}) it follows that $\mathsf{Up}(\mathbb{X}) \vDash \Phi$ iff $\powersetm{\mathbb{X}} \vDash \Phi_g$. Furthermore, as $\Phi_g$ is a modal Sahlqvist quasiequation by Lemma \ref{Lem:still-Sahlqvist}, we can apply the correspondence part of the Modal Sahlqvist Theorem obtaining that $\powersetm{\mathbb{X}} \vDash \Phi_g$ iff $\mathbb{X} \vDash \mathsf{mtr}(\Phi_g)$, where the first order sentence $\mathsf{mtr}(\Phi_g)$ is effectively computable. Therefore, setting $\mathsf{tr}(\Phi) \coloneqq \mathsf{mtr}(\Phi_g)$, we are done.
\end{proof}

\begin{exa}\label{Exa:BTWL-correspondence-easy}
Let $n \in \mathbb{Z}^+$ and consider the Sahlqvist quasiequation $\Phi_n$ associated with the formula $\mathsf{btw}_n$ defined in Example \ref{Exa:BTWn}. Since $\Phi_n$ and $\mathsf{btw}_n$ are equivalent over Heyting algebras, for every poset $\mathbb{X}$ we have
\[
\mathsf{Up}(\mathbb{X}) \vDash \Phi_n \, \,\Longleftrightarrow\, \, \mathsf{Up}(\mathbb{X}) \vDash \mathsf{btw}_n.
\]
On the other hand, it is known that $\mathsf{Up}(\mathbb{X})\vDash \mathsf{btw}_n$ iff for every $x \in X$ and $y_1, \dots, y_{n+1} \in {\uparrow}x$ there exist $z_1, \dots, z_n \in {\uparrow}x$ such that $y_1, \dots, y_{n+1} \in {\downarrow}\{z_1, \dots, z_n \}$ (see, e.g., \cite[Exercise 2.11]{ChZa97}). 

As the latter condition can be rendered as a first order sentence $\Psi_n$ in the language of posets, we obtain the following instance of the correspondence part of the Intuitionistic Sahlqvist Theorem: for every poset $\mathbb{X}$,
\[
\mathsf{Up}(\mathbb{X}) \vDash \Phi_n \, \,\Longleftrightarrow \, \,\mathbb{X} \vDash \Psi_n,
\]
whence $\mathsf{tr}(\Phi_{n})$ is logically equivalent to $\Psi_n$ over the class of posets. When $n = 1$, the condition $\Psi_n$ expresses the demand that the principal upsets of $\mathbb{X}$ are up-directed.

By the same token, when $\Phi$ is the Sahlqvist quasiequation associated with the excluded middle axiom, $\mathsf{tr}(\Phi)$ expresses the demand that the poset $\mathbb{X}$ is discrete. Lastly, when $\Phi$ is the Sahlqvist quasiequation associated with the G\"odel-Dummett axiom, $\mathsf{tr}(\Phi)$ is the sentence expressing the demand that $\mathbb{X}$ is a \emph{root system}, i.e., that ${\uparrow}x$ is a chain, for every $x \in X$ (see \cite{Ho69} or \cite[Prop.\ 2.36]{ChZa97}).
\qed
\end{exa}
\color{black}

\section{Sahlqvist theory for fragments of $\mathsf{IPC}$ with conjunction}

Recall that $\mathcal{L}$ is the algebraic language of $\mathsf{IPC}$, namely,
\begin{align*}	
\mathcal{L} = 		 x \mid  \varphi \land \psi  \mid  \varphi \lor \psi  \mid  \varphi \to \psi  \mid  \lnot \varphi  \mid 0  \mid  1.
\end{align*}
	\color{black}
The aim of this section is to extend Sahlqvist theory to fragments of $\mathsf{IPC}$ including the connective $\land$.\footnote{For a similar result for fragments $\mathsf{IPC}$ containing $\to$, see Theorem \ref{Thm:canonicity-arrow-IPC}.} As the correspondence part of Sahlqvist theorem is left unchanged by switching to fragments, the main result of this section takes the form of a canonicity result.

\begin{Theorem}\label{thm : enhanced intuitionistic sahlqvist}
Let $\Phi$ be a Sahlqvist quasiequation in a sublanguage $\mathcal{L}_\land$ of $\mathcal{L}$ containing $\land$. If an $\mathcal{L}_\land$-subreduct $\A$ of a Heyting algebra validates $\Phi$, then also $\mathsf{Up}(\A_\ast)$ validates $\Phi$.
\end{Theorem}

In order to prove the above result, we begin by ruling out some limit cases.

\begin{Proposition}
Let $\Phi$ be a Sahlqvist quasiequation in a language $\mathcal{L}_\land \subseteq \{ \land, \lor, 0, 1 \}$. If an $\mathcal{L}_\land$-subreduct $\A$ of a Heyting algebra validates $\Phi$, then also $\mathsf{Up}(\A_\ast)$ validates $\Phi$.
\end{Proposition} 

\begin{proof}
It is well known that, in view of the poorness of the language $\mathcal{L}_\land$, the class $\mathsf{K}$ of $\mathcal{L}_\land$-subreducts of Heyting algebras is a minimal quasivariety.\footnote{For instance, if $\mathcal{L}_\land = \{ \land, \lor, 0, 1 \}$, then $\mathsf{K}$ is the class of bounded distributive lattices.} This means that every quasiequation in $\mathcal{L}_\land$ is either true in $\mathsf{K}$ or false in all the nontrivial members of $\mathsf{K}$.

Suppose that $\Phi$ is valid in some $\A \in \mathsf{K}$. If $\mathsf{K} \vDash \Phi$, then $\Phi$ is also valid in the $\mathcal{L}_\land$-reduct of the Heyting algebra $\mathsf{Up}(\A_\ast)$. This, in turn, implies that $\mathsf{Up}(\A_\ast) \vDash \Phi$ as desired. Then we consider the case where $\Phi$ is false in all the nontrivial members of $\mathsf{K}$. In this case, the assumption that $\A \vDash \Phi$ forces $\A$ to be trivial. Therefore, $\A_\ast$ is the empty poset and the Heyting algebra $\mathsf{Up}(\A_\ast)$ is trivial. As a consequence, $\mathsf{Up}(\A_\ast)$ validates every quasiequation and, in particular, $\Phi$.
\end{proof}

In order to prove Theorem \ref{thm : enhanced intuitionistic sahlqvist}, it only remains to consider the cases where $\mathcal{L}_\land$ contains $\land$ and either $\lnot$ or $\to$. Up to term-equivalence, this amounts to proving that for every variety $\mathsf{K}$ among $\mathsf{PSL}, \mathsf{(b)ISL}, \mathsf{PDL}$, $\mathsf{IL}$, and $\mathsf{HA}$ the following holds: for every Sahlqvist quasiequation $\Phi$ in the language of $\mathsf{K}$ and every $\A \in \mathsf{K}$, if $\A$ validates $\Phi$, then also $\mathsf{Up}(\A_\ast)$ validates $\Phi$.

The next result does this for the case where $\mathsf{K}$ is any variety among $\mathsf{(b)ISL}, \mathsf{PDL}$, $\mathsf{IL}$, and $\mathsf{HA}$ (i.e., all cases except $\mathsf{K} = \mathsf{PSL})$.

\begin{Proposition}\label{Prop:Canonicity-all-but-PSL}
Let $\mathsf{K}$ be a variety among $\mathsf{(b)ISL}, \mathsf{PDL}$, $\mathsf{IL}$, and $\mathsf{HA}$ and $\Phi$ a Sahlqvist quasiequation in the language of $\mathsf{K}$. For every $\A \in \mathsf{K}$, if $\A$ validates $\Phi$, then also $\mathsf{Up}(\A_\ast)$ validates $\Phi$.
\end{Proposition}

\begin{proof}
Consider a variety $\mathsf{K}$ among $\mathsf{(b)ISL}, \mathsf{PDL}$, $\mathsf{IL}$, and $\mathsf{HA}$, a Sahlqvist quasiequation $\Phi$ in the language of $\mathsf{K}$, and an algebra $\A \in \mathsf{K}$ such that $\A \vDash \Phi$. By Theorem \ref{Thm:general-embedding}, $\A$ embeds into the appropriate reduct $\B^-$ of a Heyting algebra $\B$ such that $\B^- \in \UUU(\A)$. Since $\Phi$ is a universal sentence valid in $\A$, from $\B^- \in \UUU(\A)$ it follows $\B^- \vDash \Phi$. As $\B^-$ is the reduct of $\B$ in the language of $\mathsf{K}$, this guarantees that $\B \vDash \Phi$. 

Given that $\Phi$ is a Sahlqvist quasiequation, we can apply the canonicity part of the Intuitionistic Sahlqvist Theorem obtaining that $\mathsf{Up}(\B_\ast) \vDash \Phi$. Since $\B_\ast = \B^-_\ast$, the algebra $\mathsf{Up}_{\mathsf{K}}(\B^{-}_\ast)$ is the reduct of $\mathsf{Up}(\B_\ast)$ in the language of $\mathsf{K}$. Consequently, from $\mathsf{Up}(\B_\ast) \vDash \Phi$ it follows that $\mathsf{Up}_{\mathsf{K}}(\B_\ast) \vDash \Phi$.

Now, recall that there exists an embedding $f \colon \A \to \B^-$. By Conditions (\ref{item:weak-adjunction1}) 
and (\ref{item:weak-adjunction2}) of Proposition \ref{Prop:weak-adjunction}, the map $\mathsf{Up}_{\mathsf{K}}(f_\ast) \colon \mathsf{Up}_{\mathsf{K}}(\A_\ast) \to \mathsf{Up}_{\mathsf{K}}(\B^{-}_\ast)$ is a homomorphism 
between members of $\mathsf{K}$. Furthermore, applying the last part of the same proposition to the assumption 
that $f$ is injective, we obtain that $\mathsf{Up}_{\mathsf{K}}(f_\ast)$ is also injective, whence $\mathsf{Up}_{\mathsf{K}}(\A_\ast) \in \III\SSS(\mathsf{Up}_{\mathsf{K}}(\B^{-}_\ast))$. Since the validity of universal 
sentences persists under the formation of subalgebras and isomorphic copies, from $\mathsf{Up}_{\mathsf{K}}(\B^-_\ast) \vDash \Phi$ it follows that $\mathsf{Up}_{\mathsf{K}}(\A_\ast) \vDash \Phi$ and, therefore, $\mathsf{Up}(\A_\ast) \vDash \Phi$, thus concluding the proof.
\end{proof}

In order to complete the proof of Theorem \ref{thm : enhanced intuitionistic sahlqvist}, it only remains to prove the following:

\begin{Proposition}\label{Prop:Canonicity-PSL-main}
Let $\Phi$ be Sahlqvist quasiequation in the language of $\mathsf{PSL}$. For every $\A \in \mathsf{PSL}$, if $\A$ validates $\Phi$, then also $\mathsf{Up}(\A_\ast)$ validates $\Phi$.
\end{Proposition}

The proof result proceeds through a series of technical observations. An element $a$ of a semilattice $\A$ is said to be \emph{join irreducible} if it is not the minimum of $\A$ and for every pair of elements $b, c  \in A$ such that the join $b \lor c$ exists in $\A$, if $a = b \lor c$, then either $a = b$ or $a = c$. We denote by $\mathsf{J}(\A)$ the subposet of $\A$ whose universe is the set of join irreducible elements. 

\begin{Lemma}\label{Lem:J(A)-properties}
The following conditions hold for a finite semilattice $\A$:
\benroman
\item\label{item:J(A)-properties:1} If $a \nleq b$, there exists $c \in \mathsf{J}(\A)$ such that $c \leq a$ and $c \nleq b$;
\item\label{item:J(A)-properties:2} An element $a \in A$ is the minimum of $\A$ iff there is no $c \in \mathsf{J}(\A)$ such that $c \leq a$.
\eroman
\end{Lemma}

\begin{proof}
This is an immediate consequence of the fact that every element of a finite semilattice $\A$ is the join of a subset of $\mathsf{J}(\A)$.
\end{proof}

Furthermore, we rely on the following properties of pseudocomplemented semilattices. 

\begin{Lemma}\label{Lem:PSL-properties}
The following conditions hold for every $\A \in \mathsf{PSL}$:
\benroman
\item\label{Eq:PSL-properties-1} If $\varphi(x_1, \dots, x_n)$ is a negative formula, the term function $\varphi^\A(x_1, \dots, x_n)$ is \emph{order reversing} in every argument, i.e., for every $a_1, \dots, a_n, b_1, \dots, b_n \in A$, 
\[
\text{if }a_i \leq b_i \text{ for every }i \leq n \text{, then }\varphi^{\A}(b_1, \dots, b_n) \leq \varphi^{\A}(a_1, \dots, a_n);
\]
\item\label{Eq:PSL-properties-2} If $\A$ is finite and $X \subseteq A$, the join $\bigvee X$ exists in $\A$ and $\bigwedge_{a \in X} \lnot a = \lnot \bigvee X$.
\eroman
\end{Lemma}

\begin{proof}
Condition (\ref{Eq:PSL-properties-1}) follows from the fact that $\lnot$ is order reversing in $\mathsf{PSL}$ \cite[Condition (9)]{Fr62}, while $\land$ is order preserving in both arguments. For Condition (\ref{Eq:PSL-properties-2}), see \cite[Condition (19)]{Fr62}.
\end{proof}

The following construction will be instrumental to deal with finite members of $\mathsf{PSL}$.

\begin{Definition} With every finite semilattice $\A$ we associate an algebra 
\[
\A^+ \coloneqq \langle \mathsf{Dw}(\mathsf{J}(\A)); \cap, \lnot, \emptyset, X \rangle,
\]
where $\mathsf{Dw}(\mathsf{J}(\A))$ is the set of downsets of $\mathsf{J}(\A)$ and $\lnot$ is defined by
\[
\lnot D \coloneqq \{ a \in \mathsf{J}(\A) : D \cap {\downarrow} a = \emptyset \}.
\]
Furthermore, let $\epsilon_\A \colon \A \to \A^+$ be the map defined by the rule
\[
\epsilon_\A(a) \coloneqq \mathsf{J}(\A) \cap {\downarrow}a.
\]
\end{Definition}

\begin{Lemma}\label{Lem:A+-is-in-PSL}
Let $\A \in \mathsf{PSL}$ be finite. Then $\A^+$ is the $\langle \land, \lnot, 0, 1 \rangle$-reduct of a Heyting algebra, it belongs to $\mathsf{PSL}$, and the map $\epsilon_\A \colon \A \to \A^+$ is an embedding.
\end{Lemma}

\begin{proof}
Notice that $\A^+$ coincides with the algebra $\mathsf{Up}_{\mathsf{PSL}}(\mathbb{X})$, where $\mathbb{X}$ is the order dual of $\mathsf{J}(\A)$. Since $\mathsf{Up}_{\mathsf{PSL}}(\mathbb{X})$ is a pseudocomplemented semilattice, we infer that so is $\A^+$. Furthermore, the definition of $\A^+$ guarantees that it is a distributive lattice (whose join operation is $\cup$). Lastly, since $\A$ is finite, so is $\A^+$. Therefore, $\A^+$ is a finite distributive pseudocomplemented lattice. By Proposition \ref{Prop:finite-BSL}(\ref{item:finite-BSL1}), we conclude that $\A^+$ is the $\langle \land, \lnot, 0, 1 \rangle$-reduct of a Heyting algebra.

Then we turn to prove that $\epsilon_\A \colon \A \to \A^+$ is an embedding. Clearly, it is well defined and preserves $\land, 0$, and $1$. Furthermore, it is injective by Lemma \ref{Lem:PSL-properties}(\ref{Eq:PSL-properties-1}). To prove that it also preserves $\lnot$, consider $a \in A$. We will prove that for every $b \in \mathsf{J}(\A)$,
\begin{align*}
b \in \lnot^{\A^+}\epsilon_{\A}(a) \, \,&\Longleftrightarrow \, \,\epsilon_{\A}(a) \cap {\downarrow}b = \emptyset \, \,\Longleftrightarrow \, \,c \nleq a \land^\A b, \text{ for every }c \in \mathsf{J}(\A)\\
\, \,& \Longleftrightarrow \, \, a \land^\A b = 0 \, \,\Longleftrightarrow \, \, b \leq \lnot^\A a \, \,\Longleftrightarrow\, \, b \in \epsilon_{\A}(\lnot^{\A} a).
\end{align*}
The first of the above equivalences holds by the definition of $\lnot$ in $\A^+$, the second and the last by the definition of $\epsilon_\A$, the third by Lemma \ref{Lem:PSL-properties}(\ref{Eq:PSL-properties-2}), and the fourth by Condition (\ref{Eq:neg-res-law}). This shows that $\lnot^{\A^+}\epsilon_{\A}(a) = \epsilon_{\A}(\lnot^{\A} a)$. Hence, we conclude that $\epsilon_\A \colon \A \to \A^+$ preserves $\lnot$ and, therefore, it is an embedding.
\end{proof}

\begin{Remark}
The embedding $\epsilon_\A \colon \A \to \A^+$ need not be an isomorphism, because $\A^+$ is always a distributive lattice, while the (semi)lattice $\A$ may fails to be distributive.
\qed
\end{Remark}

We rely on the following technical observation.

\begin{Lemma}\label{lemma : technical two}
Let $\A \in \mathsf{PSL}$ be finite and $\varphi (x_1,\dots,x_n)$ a formula in the language of $\mathsf{PSL}$. For every $D_1,\dots,D_n \in \mathsf{Dw}(\mathsf{J}(\A))$, we have
	\[
	\lnot^{\A^+} \varphi^{\A^+} (D_1,\dots,D_n) = \lnot^{\A^+} \varphi^{\A^+}(\epsilon_\A(\bigvee^{\A} D_1),\dots,\epsilon_{\A}(\bigvee^\A D_n)). 
	\] 
\end{Lemma}

\noindent \emph{Proof.}
We begin by proving the following:

\begin{Claim}\label{Claim : technical one}
For every $D,V \in \mathsf{Dw}(\mathsf{J}(\A))$, we have
\[
\lnot^{\A^+}(D \cap V) = \lnot^{\A^+} (\epsilon_{\A}(\bigvee^{\A} D) \cap V).
\] 
\end{Claim}

\begin{proof}[Proof of the Claim]
In order to prove the inclusion from right to left, observe that for every $a \in D$ we have $a \leq \bigvee^{\A}D$ and, therefore, $a \in \epsilon_{\A}(\bigvee^{\A} D)$. Consequently, $D \subseteq \epsilon_{\A}(\bigvee^{\A} D)$. This, in turn, implies $D \cap V \subseteq \epsilon_{\A}(\bigvee^{\A} D) \cap V$. Bearing in mind that $\A^+ \in \mathsf{PSL}$ (Lemma \ref{Lem:A+-is-in-PSL}), we can apply Lemma \ref{Lem:PSL-properties}(\ref{Eq:PSL-properties-1}) obtaining that the operation $\lnot^{\A^+}$ is order reversing. Thus, $\lnot^{\A^+} (\epsilon_{\A}(\bigvee^{\A} D) \cap V) \subseteq \lnot^{\A^+}(D \cap V)$ as desired.
	
In order to prove the inclusion from left to right, we reason by contraposition. Consider $a \in \mathsf{J}(\A) \smallsetminus \lnot^{\A^+} (\epsilon_{\A}(\bigvee^{\A} D) \cap V)$. By the definitions of $\lnot^{\A^+}$ and $\epsilon_\A$, there exists $b \in V$ such that $b \leq \bigvee^{\A} D, a$. We have two cases depending on whether or not there exists $d \in D$ such that $b  \nleq \lnot^{\A}d$. 

Suppose first that such a $d$ exists. In view of Condition (\ref{Eq:neg-res-law}) we get $0 < b \land^{\A}d$. Therefore, Lemma \ref{Lem:J(A)-properties}(\ref{item:J(A)-properties:2}) gives us some $c \in \mathsf{J}(\A)$ such that $c \leq b, d$. Since $b \in V$, $c \in \mathsf{J}(\A)$, and $V$ is a downset of $\mathsf{J}(\A)$, we have $c \in V$. Furthermore, from $c \leq d \leq \bigvee^{\A}D$ and $c \in \mathsf{J}(\A)$ it follows that $c \in \epsilon_{\A}(\bigvee^{\A}D)$. Lastly, since $c \leq b \leq a$, we have $c \in {\downarrow}a$. Thus, $c \in \epsilon_{\A}(\bigvee^{\A}D) \cap V \cap {\downarrow}a$. By the definition of $\lnot^{\A^+}$, this amounts to $a \notin \lnot^{\A^+}(\epsilon_{\A}(\bigvee^{\A}D) \cap V)$ and we are done. 

To conclude the proof, it only remains to show that the case where $b  \leq \lnot^{\A}d$ for every $d \in D$ never happens. Suppose the contrary. By Lemma \ref{Lem:PSL-properties}(\ref{Eq:PSL-properties-2}), we have
\[
b \leq \bigwedge^{\A}_{d \in D}\lnot d = \lnot^{\A} \bigvee^{\A}D.
\]
As we assumed that $b \leq \bigvee^{\A}D$, this yields $b \leq (\bigvee^{\A}D) \land^{\A}(\lnot^{\A} \bigvee^{\A}D)$ which, by Condition (\ref{Eq:neg-res-law}), amounts to $b = 0$. But this contradicts with the fact that $b \in \mathsf{J}(\A)$.
\end{proof}

To prove the main statement, we reason by induction on the construction of $\varphi$. In the base case, $\varphi$ is either a constant or a variable. The case where $\varphi$ is a constant is straightforward. If $\varphi$ is a variable $x_i$, by applying the Claim in the third equality below, we obtain
\begin{align*}
\lnot^{\A^+} \varphi^{\A^+} (D_1,\dots,D_n) &= \lnot^{\A^+}D_i = \lnot^{\A^+}(D_i \cap \mathsf{J}(\A)) = \lnot^{\A^+}(\epsilon_\A(\bigvee^{\A}D_i) \cap \mathsf{J}(\A)) \\
&=\lnot^{\A^+}\epsilon_\A(\bigvee^{\A}D_i) = \lnot^{\A^+} \varphi^{\A^+}(\epsilon_\A(\bigvee^{\A} D_1),\dots,\epsilon_\A(\bigvee^{\A} D_n)).
\end{align*}

In the step case, the principal connective of $\varphi$ is either $\lnot$ or $\land$. The case where it is $\lnot$ follows immediately from the inductive hypothesis. Therefore, we detail only the case where the principal connective of $\varphi$ is $\land$. Since the operation $\land$ is associative and commutative in $\mathsf{PSL}$, we may assume that $\varphi$ is of the form
\[
\lnot \alpha_1 \land \dots \land \lnot \alpha_{m} \land \beta_1 \land \dots \land \beta_k \land x_{i_1} \land \dots \land x_{i_t},
\]
where $i_1, \dots, i_t \leq n$ and each $\beta_j$ is a constant. Furthermore, $m, k$, or $t$ can be $0$. As the inductive hypothesis applies to each $\alpha_j$, we obtain
\[
\lnot^{\A^+}\alpha_j^{\A^+}(D_1, \dots, D_n) = \lnot^{\A^+}\alpha_j^{\A^+}(\epsilon_\A(\bigvee^{\A} D_1),\dots,\epsilon_\A(\bigvee^{\A} D_n)), \text{ for every }j \leq m.
\]
Furthermore, as the various $\beta_j$ are constants, we have
\[
\beta_j^{\A^+}(D_1, \dots, D_n) = \beta_j^{\A^+}(\epsilon_\A(\bigvee^{\A} D_1),\dots,\epsilon_\A(\bigvee^{\A} D_n)), \text{ for every }j \leq k.
\]
Therefore, setting
\[
V \coloneqq \bigcap_{j \leq m}\lnot^{\A^+}\alpha_j^{\A^+}(\epsilon_\A(\bigvee^{\A} D_1),\dots,\epsilon_\A(\bigvee^{\A} D_n)) \cap \bigcap_{j \leq k}\beta_j^{\A^+}(\epsilon_\A(\bigvee^{\A} D_1),\dots,\epsilon_\A(\bigvee^{\A} D_n)),
\]
we obtain
\[
\lnot^{\A^+} \varphi^{\A^+} (D_1,\dots,D_n) = \lnot^{\A^+}(D_{i_1} \cap \dots \cap D_{i_t} \cap V).
\]
Lastly, applying $t$ times the Claim to the above display, we get
\[
\lnot^{\A^+} \varphi^{\A^+} (D_1,\dots,D_n) = \lnot^{\A^+}(\epsilon_{\A}(\bigvee^{\A}D_{i_1}) \cap \dots \cap \epsilon_{\A}(\bigvee^{\A}D_{i_t}) \cap V)
\]
which, by the definition of $V$, amounts to
\[
\pushQED{\qed}	\lnot^{\A^+} \varphi^{\A^+} (D_1,\dots,D_n) = \lnot^{\A^+} \varphi^{\A^+}(\epsilon_\A(\bigvee^{\A} D_1),\dots,\epsilon_{\A}(\bigvee^\A D_n)).  \qedhere \popQED
\]

 The next result is the hearth of the proof of Proposition \ref{Prop:Canonicity-PSL-main}.

\begin{Proposition}\label{Prop:PSL-canonicty-hard-part}
Let $\A \in \mathsf{PSL}$ be finite and $\Phi$ a Sahlqvist quasiequation in the language of $\mathsf{PSL}$. If $\A$ validates $\Phi$, then also $\A^+$ validates $\Phi$.
\end{Proposition}

\begin{proof}
We will reason by contraposition. Consider a Sahlqvist quasiequation
\[
\Phi = \varphi_1(x_1, \dots, x_k) \land y \leq z \, \& \dots \& \, \varphi_n(x_1, \dots, x_k) \land y \leq z \Longrightarrow y \leq z
\]
in the language of $\mathsf{PSL}$ such that $\A^+ \nvDash \Phi$. We need to prove that $\A \nvDash \Phi$. 

Since $\A^+$ is the $\langle \land, \lnot, 0, 1 \rangle
$-reduct of a Heyting algebra (Lemma \ref{Lem:A+-is-in-PSL}), we can apply Corollary \ref{Cor:meaning-Sahlqvist-IPC} to the assumption that $\A^+ \nvDash \Phi$, obtaining $D_1, \dots, D_k \in \mathsf{Dw}(\mathsf{J}(\A))$ such that
\[
\varphi_1^{\A^+}(D_1, \dots, D_k) \cup \dots \cup \varphi_n^{\A^+}(D_1, \dots, D_k) \ne \mathsf{J}(\A).
\]
Let then $a \in \mathsf{J}(\A)$ be such that
\begin{equation}\label{Eq:the-missing-step-PSL-7c}
a \notin \varphi_1^{\A^+}(D_1, \dots, D_k) \cup \dots \cup \varphi_n^{\A^+}(D_1, \dots, D_k).
\end{equation}
Recall that $\A$ is a finite semilattice with a maximum and, therefore, it is also a lattice. Thereby, for every $m \leq k$ we can define an element of $\A$ as follows:
\[
b_m \coloneqq \bigvee^{\A}_{d \in D_m}(d \land a).
\]
We will prove that
	\begin{equation}\label{Eq:the-missing-step-PSL-7}
	a \nleq (\varphi_1^\A(b_1, \dots, b_k) \land^\A a) \lor^{\A} \dots \lor^{\A} (\varphi_n^\A(b_1, \dots, b_k) \land^\A a).
	\end{equation}
This, in turn, implies that $\A \nvDash \Phi$, as witnessed by the assignment
\[
x_m \longmapsto b_m \quad  y \longmapsto a \quad z \longmapsto(\varphi_1^\A(b_1, \dots, b_k) \land^\A a) \lor^{\A} \dots \lor^{\A} (\varphi_n^\A(b_1, \dots, b_k) \land^\A a).
\]	
Therefore, to conclude the proof, it suffices to establish Condition (\ref{Eq:the-missing-step-PSL-7}).

Suppose, with a view to contradiction, that Condition (\ref{Eq:the-missing-step-PSL-7}) fails. Then 
\[
	a = (\varphi_1^\A(b_1, \dots, b_k) \land^\A a) \lor^{\A} \dots \lor^{\A} (\varphi_n^\A(b_1, \dots, b_k) \land^\A a).
\]
Since $a$ is join irreducible, by symmetry we may assume that $a = \varphi_1^\A(b_1, \dots, b_k) \land^\A a$, that is,
\begin{equation}\label{Eq:the-missing-step-PSL-7b}
a \leq \varphi_1^\A(b_1, \dots, b_k).
\end{equation}
Now, recall that $\varphi_1$ is obtained from Sahlqvist implications using only $\land, \lor$, and $\Box$. Since $\varphi_1$ is in the language of $\mathsf{PSL}$, this means that $\varphi_1$ is a conjunction of Sahlqvist implications. Consequently, we may assume that
	\begin{equation}\label{Eq:very-last-step-canonicity}
	\varphi_1 = \bigwedge_{i \leq p} \gamma_j \land \bigwedge_{j \leq q} \lnot \psi_i,
	\end{equation}
	where the various $\gamma_i$ and $\psi_j$ are, respectively, Sahlqvist antecedents and positive formulas, both in the language of $\mathsf{PSL}$. Furthermore, $p$ or $q$ can be $0$. Without loss of generality, we may assume each $\gamma_i$ is a variable. This is because if $\gamma_i = \lnot \alpha$ then $\alpha$ is a negative formula and, therefore, a Sahlqvist antecendent. Consequently, we may assume that $\gamma_i = \lnot \psi_j$ for some $j \leq q$ and remove $\gamma_i$ from the big conjunction on the left hand side of the above display. On the other hand, if $\gamma_i = \alpha \land \beta$, then both $\alpha$ and $\beta$ are positive formulas and, therefore, we may assume that there are $i_1, i_2 \leq p$ such that $\gamma_{i_1} = \alpha$ and $\gamma_{i_2} = \beta$ and remove $\gamma_i$ from the big conjunction on the left hand side of the above display. Iterating this process, we may assume that in the above display every $\gamma_i$ is either a constant or a variable, while the various $\psi_j$ are still Sahlqvist antecedents. Lastly, if some $\gamma_i$ is the constant $1$, we can remove it from the big conjunction on the left hand side of the above display, thereby producing a new formula that is still equivalent to $\varphi_1$ in $\mathsf{PSL}$. This is possible because $\varphi_1$ cannot simply be the constant $1$, otherwise Condition (\ref{Eq:the-missing-step-PSL-7}) would hold, contradicting the assumption. Moreover, no $\gamma_i$ is the constant $0$, otherwise Condition (\ref{Eq:the-missing-step-PSL-7b}) would imply that $a = 0$, contradicting the assumption that $a \in \mathsf{J}(\A)$. Therefore, we may assume that each $\gamma_i$ in Condition (\ref{Eq:very-last-step-canonicity}) is a variable. In addition, we may also assume that the various $\gamma_i$ are pairwise distinct and, renaming the variables when necessary, that each $\gamma_i$ is the variable $x_i$, thereby obtaining
		\[
\varphi_1 =	x_{1} \land \dots \land x_{p} \land  \lnot \psi_1 \land \dots \land \lnot \psi_q,
	\]
	where the various $\psi_j$ are Sahlqvist antecedents in the language of $\mathsf{PSL}$. 
	
	 In view of Condition (\ref{Eq:the-missing-step-PSL-7b}), this yields
\begin{equation}\label{Eq:the-missing-step-PSL-7d}
a \leq b_{1} \land^{\A} \dots \land^{\A} b_{p} \land^{\A} \lnot \psi_1^{\A}(b_1, \dots, b_k) \land^{\A} \dots \land^{\A} \lnot \psi_q^{\A}(b_1, \dots, b_k).
\end{equation}
On the other hand, from Condition (\ref{Eq:the-missing-step-PSL-7c}) it follows that
\[
a \notin D_{1} \cap \dots \cap D_{p} \cap \lnot^{\A^+} \psi_1^{\A^+}(D_1, \dots, D_k) \cap \dots \cap \lnot^{\A^+} \psi_q^{\A^+}(D_1, \dots, D_k).
\]
We have two cases depending on whether
\[
a \notin D_{1} \cap \dots \cap D_{p}\, \, \text{ or }\, \, a \notin \lnot^{\A^+} \psi_1^{\A^+}(D_1, \dots, D_k) \cap \dots \cap \lnot^{\A^+} \psi_q^{\A^+}(D_1, \dots, D_k).
\]

Suppose first that $a \notin D_{1} \cap \dots \cap D_{p}$. By symmetry, we may assume that $a \notin D_{1}$. From Condition (\ref{Eq:the-missing-step-PSL-7d}) and the definition of $b_{1}$ it follows that
\[
a \leq b_{1} = \bigvee_{d \in D_{1}}^{\A}(d \land^{\A} a).
\]
This amounts to $a = \bigvee_{d \in D_{1}}^{\A}(d \land^{\A} a)$ which, in turn, implies that $a \leq d$ for some $d \in D_{1}$ because $a \in \mathsf{J}(\A)$. Since $a \in \mathsf{J}(\A)$ and $D_{1}$ is a downset of $\mathsf{J}(\A)$, we conclude that $a \in D_{1}$, a contradiction.

Then we consider the case where $a \notin \lnot^{\A^+} \psi_1^{\A^+}(D_1, \dots, D_k) \cap \dots \cap \lnot^{\A^+} \psi_q^{\A^+}(D_1, \dots, D_k)$. By symmetry, we may assume that $a \notin \lnot^{\A^+} \psi_1^{\A^+}(D_1, \dots, D_k)$. Applying in sequence Lemma \ref{lemma : technical two} and the fact that $\epsilon_\A \colon \A \to \A^+$ is a homomorphism (Lemma \ref{Lem:A+-is-in-PSL}), we deduce 
\begin{align*}
a &\notin \lnot^{\A^+} \psi_1^{\A^+}(D_1, \dots, D_k) \\
&= \lnot^{\A^+} \psi_1^{\A^+}(\epsilon_\A(\bigvee^{\A} D_1),\dots,\epsilon_{\A}(\bigvee^\A D_k))\\
& = \epsilon_{\A}(\lnot^{\A} \psi_1^{\A}(\bigvee^{\A} D_1, \dots, \bigvee^{\A} D_k)). 
\end{align*}
Since $a \in \mathsf{J}(\A)$, by the definition of $\epsilon_\A$ this amounts to
\begin{equation}\label{Eq:the-missing-step-PSL-7dkss}
a \nleq \lnot^{\A} \psi_1^{\A}(\bigvee^{\A} D_1, \dots, \bigvee^{\A} D_k).
\end{equation}

Now, as $\psi_1$ is a Sahlqvist antecedent in the language of $\mathsf{PSL}$, it is a conjunction of variables, negative formulas, and constants.\ As before, we can remove the constants from this conjunction. Therefore, we may assume that $\psi_1$ is of the form
\begin{equation}\label{Eq:the-missing-step-PSL-7dk}
x_{1} \land \dots \land x_{p'} \land \alpha_1 \land \dots \land \alpha_{q'},
\end{equation}
where the various $\alpha_j$ are negative formulas. Furthermore, $p'$ or $q'$ can be $0$.

As the various $\alpha_j$ are negative formulas, the term function $\alpha_j^\A$ is order reversing in every argument by Lemma \ref{Lem:PSL-properties}(\ref{Eq:PSL-properties-1}). Bearing in mind that for every $m \leq k$ we have 
\[
b_m = \bigvee^{\A}_{d \in D_m} (d \land^{\A} a) \leq \bigvee^\A D_m,
\]
this implies that for every $j \leq q'$,
\[
\alpha_j^{\A}(\bigvee^{\A}D_1, \dots, \bigvee^{\A}D_k) \leq \alpha_j^{\A}(b_1, \dots, b_k).
\]
Since $\psi_1$ is the formula in Condition (\ref{Eq:the-missing-step-PSL-7dk}), we obtain
\[
\psi_1^{\A}(\bigvee^{\A}D_1, \dots, \bigvee^{\A}D_k) \leq \bigwedge^{\A}_{i \leq p'}\bigvee^{\A}D_{i}  \land^\A \bigwedge^{\A}_{j \leq q'} \alpha_j^{\A}(b_1, \dots, b_k).
\]
By applying the fact that the negation operation is order reversing in $\mathsf{PSL}$ to the above display and 
Condition (\ref{Eq:the-missing-step-PSL-7dkss}), we obtain 
\[
a \nleq \lnot^\A\Big(\bigwedge^{\A}_{i \leq p'}\bigvee^{\A}D_{i}  \land^\A \bigwedge^{\A}_{j \leq q'} \alpha_j^{\A}(b_1, \dots, b_k)\Big).
\]
In view of Condition (\ref{Eq:neg-res-law}), this amounts to
\begin{equation}\label{Eq:dwgijfdwkfmcpio}
0 < a \land^\A \bigwedge_{i \leq p'} \bigvee^{\A}D_{i} \land^\A \bigwedge^{\A}_{j \leq q'} \alpha_j^{\A}(b_1, \dots, b_k).
\end{equation}

We will prove that 
\begin{equation}\label{Eq:agvasqjhzvx}
0 < a \land^\A b_{1} \land^\A \bigvee^{\A}_{2 \leq i \leq p'}D_{i} \land^\A \bigwedge^{\A}_{j \leq q'} \alpha_j^{\A}(b_1, \dots, b_k).
\end{equation}

By applying Condition (\ref{Eq:neg-res-law}) to Condition (\ref{Eq:dwgijfdwkfmcpio}) and, subsequently, Lemma \ref{Lem:PSL-properties}(\ref{Eq:PSL-properties-2}), we obtain
\[
a \land^\A\bigvee^{\A}_{2 \leq i \leq p'}D_{i}  \land^\A \bigwedge^{\A}_{j \leq q'} \alpha_j^{\A}(b_1, \dots, b_k) \nleq \lnot \bigvee^{\A}D_{1} = \bigwedge^{\A}_{d \in D_{1}} \lnot^\A d.
\]
Consequently, there exists $d_1 \in D_{1}$ such that
\[
a \land^\A\bigvee^{\A}_{2 \leq i \leq p'}D_{i}  \land^\A \bigwedge^{\A}_{j \leq q'} \alpha_j^{\A}(b_1, \dots, b_k)  \nleq \lnot^\A d_1.
\]
By applying Condition (\ref{Eq:neg-res-law}) twice, this yields
\begin{equation}\label{Eq:PSL-missing-equation-bm}
a\land^\A\bigvee^{\A}_{2 \leq i \leq p'}D_{i}  \land^\A \bigwedge^{\A}_{j \leq q'} \alpha_j^{\A}(b_1, \dots, b_k)  \nleq \lnot^{\A}(d_1 \land^\A a).
\end{equation}
By the definition of $b_{1}$ and Lemma \ref{Lem:PSL-properties}(\ref{Eq:PSL-properties-2}) we have
\[
\lnot^\A b_{1} = \lnot^\A \bigvee^\A_{d \in D_{1}}(d \land^\A a) = \bigwedge^\A_{d \in D_{1}}(\lnot^{\A}(d \land^\A a)) \leq \lnot^{\A}(d_1 \land^\A a).
\]
Together with Condition (\ref{Eq:PSL-missing-equation-bm}), this yields
\[
a\land^\A\bigvee^{\A}_{2 \leq i \leq p'}D_{i}  \land^\A \bigwedge^{\A}_{j \leq q'} \alpha_j^{\A}(b_1, \dots, b_k)   \nleq \lnot^\A b_{1}.
\]
By Condition (\ref{Eq:neg-res-law}) this amounts to Condition (\ref{Eq:agvasqjhzvx}) as desired.

Iterating $p-1$ times the argument described for Condition (\ref{Eq:agvasqjhzvx}), where the role of $D_1$ is taken successively by $D_2, \dots, D_{p'}$, we obtain
\[
0 < a \land^\A  \bigwedge^\A_{i \leq p'} b_{i} \land^\A \bigwedge^{\A}_{j \leq q'} \alpha_j^{\A}(b_1, \dots, b_k).
\]
By Condition (\ref{Eq:neg-res-law}) and the fact that $\psi_1$ is the formula in Condition (\ref{Eq:the-missing-step-PSL-7dk}) this amounts to
\[
a \nleq \lnot^\A(\bigwedge^\A_{i \leq p'} b_{i} \land^\A \bigwedge^{\A}_{j \leq q'} \alpha_j^{\A}(b_1, \dots, b_k)) = \lnot^\A \psi_1^\A(b_1, \dots, b_k),
\]
a contradiction with Conditions (\ref{Eq:the-missing-step-PSL-7b}) and (\ref{Eq:very-last-step-canonicity}).
\end{proof}

We are now ready to conclude the proof of Proposition \ref{Prop:Canonicity-PSL-main}.

\begin{proof}
Suppose that $\A \vDash \Phi$. In view of Proposition \ref{Prop:locally-finite}, the finitely generated subalgebras of $\A$ are finite. Therefore, we can apply Lemma \ref{Lem:A+-is-in-PSL} obtaining that every finitely generated subalgebra $\C$ of $\A$ embeds into $\C^+$. Together with Proposition \ref{Prop:embedding}, this implies that there exist a family $\{ \A_i : i \in I \}$ of finitely generated subalgebras of $\A$ and an ultrafilter $U$ on $I$ with an embedding
\[
f \colon \A \to \prod_{i \in I}\A_i^+ /U.
\]

Consider $i \in I$. Since the validity of universal sentences persists in subalgebras, from $\A \vDash \Phi$ it follows $\A_i \vDash \Phi$. Therefore, Proposition \ref{Prop:PSL-canonicty-hard-part} guarantees that $\A_i^+\vDash\Phi$. As a consequence, all the factors of the ultraproduct in the above display validate $\Phi$. Since the validity of universal sentences persists in ultraproducts, we conclude that $\B^{-} \vDash \Phi$ for $\B^- \coloneqq \prod_{i \in I}\A_i^+ /U$.

Now, recall from Lemma \ref{Lem:A+-is-in-PSL} that each $\A_i^+$ is the $\langle \land, \lnot, 0, 1 \rangle$-reduct of a Heyting algebra $\B_i$. Therefore, $\B^-$ is the $\langle \land, \lnot, 0, 1 \rangle$-reduct of the ultraproduct $\B \coloneqq \prod_{i\in I}\B_i /U$. As $\Phi$ is in the language of $\mathsf{PSL}$ and $\B^- \vDash \Phi$, this implies that $\B \vDash \Phi$. Lastly, since $\mathsf{HA}$ is closed under $\PPU$, we have $\B \in \mathsf{HA}$. 

In sum, $\A$ embeds into some $\B^- \in \mathsf{PSL}$ that is the $\langle \land, \lnot, 0, 1 \rangle$-reduct of a Heyting algebra $\B$ such that $\B \vDash \Phi$. Because of this, we can repeat the argument detailed in the last two paragraphs of the proof of Proposition \ref{Prop:Canonicity-all-but-PSL}, thereby obtaining $\mathsf{Up}(\A_\ast) \vDash \Phi$.
\end{proof}

As a consequence of Theorem \ref{thm : enhanced intuitionistic sahlqvist}, we obtain the following:

\begin{Corollary}\label{corollary : corollary of enhanced sahlqvist}
Let $\Phi$ be a Sahlqvist quasiequation in a sublanguage $\mathcal{L}_\land$ of $\mathcal{L}$ containing $\land$. For every $\mathcal{L}_\land$-subreduct $\A$ of a Heyting algebra, it holds that $\boldsymbol{A} \vDash \Phi \text{ iff } \boldsymbol{A}_\ast \vDash \mathsf{tr}(\Phi)$.
\end{Corollary}
\begin{proof}
In view of the correspondence part of the Intuitionistic Sahlqvist Theorem, we have that $\mathsf{Up}(\A_\ast) \vDash \Phi$ iff $\boldsymbol{A}_\ast \vDash \mathsf{tr}(\Phi)$. Therefore, in order to complete the proof, it suffices to show that $\boldsymbol{A} \vDash \Phi$ iff $\mathsf{Up}(\A_\ast) \vDash \Phi$.

On the one hand, Theorem \ref{thm : enhanced intuitionistic sahlqvist} guarantees that $\boldsymbol{A} \vDash \Phi$ implies $\mathsf{Up}(\A_\ast) \vDash \Phi$. On the other hand, $\mathsf{Up}(\A_\ast) \vDash \Phi$ implies 
$\boldsymbol{A} \vDash \Phi$, because $\A$ embeds into the $\mathcal{L}_\land$-reduct of $\mathsf{Up}(\A_\ast)$ via the map defined by the rule
\[
a \longmapsto \{ F \in \A_\ast : a \in F \}\footnote{The proof that this map is a well-defined embedding of $\A$ into the $\mathcal{L}_\land$-reduct of $\mathsf{Up}(\A_\ast)$ is analogous to the proof that a Heyting algebra $\B$ embeds into $\mathsf{Up}(\B_\ast)$ typical of Esakia duality \cite{Es74,Esakia-book85}, the only difference being that, in our case, the role of the Prime Filter Theorem is played by the observation that for every $a, b \in A$ such that $a \nleq b$ there exisits $F \in \A_\ast$ with $a \in F$ and $b \notin F$.}
\]
and the validity of universal sentences persists in subalgebras.
\end{proof}

\begin{exa}
In view of Example \ref{Exa:BTWL-correspondence-easy} and Corollary \ref{corollary : corollary of enhanced sahlqvist}, a pseudocomplemented semilattice $\A$ validates the Sahlqvist quasiequation $\Phi_n$ associated with the bounded top width $n$ axiom $\mathsf{btw}_n$ iff every $(n+1)$-element antichain in a principal upset of $\A_\ast$ is below one that has at most $n$ elements.
\qed
\end{exa}

\section{Abstract algebraic logic}

In this section we will review the rudiments of abstract algebraic logic necessary to formulate a version of Sahlqvist theory amenable to arbitrary deductive systems \cite{CN20xx-book,Cz01,AAL-AIT-f,FJa09}.\ Recall that $Var$ is the set of variables $\{ x_n : n \in \mathbb{Z}^+ \}$.\ A \emph{logic} (or a \emph{deductive system}) $\vdash$ is a consequence relation on the set of formulas with variables in $Var$ of an algebraic language that, moreover, is \emph{substitution invariant}, in the sense that for every set $\Gamma \cup \{ \varphi\}$ of formulas and substitution $\sigma$,
\[
\text{if }\Gamma \vdash \varphi, \text{ then }\sigma[\Gamma] \vdash \sigma(\varphi).
\]
Furthermore, we will assume that logics $\vdash$ are \emph{finitary}, in the sense that 
\[
\text{if }\Gamma \vdash \varphi \text{, there exists a finite }\Sigma \subseteq \Gamma \text{ s.t. }\Sigma \vdash \varphi.
\]

We denote the set of formulas over which a logic $\vdash$ is defined and the corresponding algebra of formulas by $Fm(\vdash)$ and $\boldsymbol{Fm}(\vdash)$, respectively. In addition, given $\Gamma \cup \{ \varphi_1, \dots, \varphi_n, \psi_1, \dots, \psi_m \} \subseteq Fm(\vdash)$, we often write $\Gamma, \varphi_1, \dots, \varphi_n \vdash \psi_1, \dots, \psi_m$ to signify that $\Gamma \cup \{ \varphi_1, \dots, \varphi_n \} \vdash \psi_i$ for every $i \leq m$.

 Henceforth, we will assume that the algebras and logics under consideration are of the same similarity type. Let $\vdash$ be a logic and $\A$ an algebra. A subset $F$ of $A$ is said to be a \emph{deductive filter} of $\vdash$ on $\boldsymbol{A}$ when for every $\Gamma \cup \{ \varphi \} \subseteq Fm(\vdash)$ such that $\Gamma \vdash \varphi$ and every homomorphism $f \colon \boldsymbol{Fm}(\vdash) \to \A$,
\[
\text{if }f[\Gamma] \sub F, \text{ then } f(\varphi) \in F.
\]
\begin{exa}\label{exa-filers-ipc}
If $\A$ is a Heyting algebra, the deductive filters of $\mathsf{IPC}$ on $\A$ coincide with the lattice filters of $\A$.
\qed
\end{exa}

The set $\Fi{\A}$ of deductive filters of $\vdash$ on $\A$ is a closure system on $A$, whence for every $X \subseteq A$ there exists the least deductive filter of $\vdash$ on $\A$ containing $X$, in symbols $\fg{X}$.
Given $a_1,\dots,a_n \in A$, we often write $\fg{a_1,\dots,a_n}$ as a shorthand for $\fg{\{a_1,\dots,a_n\}}$.
 Furthermore, when ordered under the inclusion relation, $\Fi{\A}$ is a complete lattice in which meets are intersections and joins are defined for every $F, G \in \Fi{\A}$ as
\[
F +^{\A} G \coloneqq \fg{F \cup G}.
\]
When $\A = \boldsymbol{Fm}(\vdash)$, we will omit the superscript $\A$ from $+^\A$ and $\fg{-}$.

In order to describe the structure of the lattice $\Fi{\A}$, recall that an element $a$ of a complete lattice $\B$ is called \emph{compact} when for every $X \subseteq B$,
\[
\text{if }a \leq \bigvee X\text{, there exists a finite }Y \subseteq X \text{ s.t. }a \leq \bigvee Y.
\]
Notice that if two elements $a,b \in B$ are compact, then so is their join $a \lor b$. Consequently, when endowed with the restriction of the operation $\lor$, the set of the compact elements of $\B$ forms a semilattice ordered, in the sense of Condition (\ref{Eq:order-semilattice-meet}), under the restriction of the dual of the lattice order of $\B$. Lastly, a complete lattice $\B$ is said to be \emph{algebraic} when every element is the join of a set of compact ones. We will rely on the following classical representation theorem for algebraic lattices:

\begin{Theorem}[\protect{\cite[Thm.\ 42]{Gr11a}}]\label{thm : algebraic lattices and filters of compact}
Every algebraic lattice is isomorphic to the lattice of filters of the semilattice of its compact elements.
\end{Theorem}

We denote the set of compact elements of $\Fi{\A}$ by $\Fic{\A}$. 

\begin{Remark}\label{Rem:FiC-structure}
When endowed with the restriction of the operation $+^\A$, the set $\Fic{\A}$ forms a semilattice ordered under the superset relation. \qed
\end{Remark}

\noindent In order to describe the elements of $\Fic{\A}$, we say that a deductive filter $F \in \Fi{\A}$ is \emph{finitely generated} when there exists a finite set $X \subseteq A$ such that $F = \fg{X}$.

\begin{Proposition}\label{prop : finitarity implies algebraicity}
	Let $\vdash$ be a logic and $\boldsymbol{A}$ an algebra. Then $\langle \Fi{\boldsymbol{A}}; \cap, +^\A \rangle$ is an algebraic lattice whose compact elements are the finitely generated ones.
\end{Proposition}

\begin{proof}
This well-known fact is essentially \cite[Thm.\ 2.23.1]{AAL-AIT-f}.
\end{proof}

\begin{exa}\label{Exa:compact-filters-IPC}
In view of the above result and Example \ref{exa-filers-ipc}, the compact deductive filters of $\mathsf{IPC}$ on a Heyting algebra $\B$ coincide with the finitely generated lattice filters of $\B$. Since the latter are precisely the principal upsets of $\B$ and the order of the semilattice $\mathsf{Fi}_{\mathsf{IPC}}^\omega(\B)$ is the superset relation (see Remark \ref{Rem:FiC-structure}), we conclude that the poset associated with the semilattice $\mathsf{Fi}_{\mathsf{IPC}}^\omega(\B)$ is isomorphic to the lattice order of $\B$.
\qed
\end{exa}

Among the deductive filters of $\vdash$, those on $\boldsymbol{Fm}(\vdash)$ will be of special interest. They are called \emph{theories} and coincide with the sets $\Gamma \subseteq Fm(\vdash)$ such that $\varphi \in \Gamma$ for every $\varphi \in Fm(\vdash)$ with $\Gamma \vdash \varphi$.\ The algebraic lattice of theories of $\vdash$ will be denoted by $\Th{\vdash}$ and the semilattice of its compact elements by $\Thc{\vdash}$. 

The structure of compact deductive filters (resp.\ theories) can be used to capture the validity of various metalogical properties, as we proceed to explain. A finite set $\Gamma \subseteq Fm(\vdash)$ is said to be \emph{inconsistent} in a logic $\vdash$ if $\Gamma \vdash \varphi$ for every $\varphi \in Fm(\vdash)$. 

\begin{Definition}
A logic $\vdash$ is said to have:
\benroman
\item\label{item:def-metalogical-1} the \emph{inconsistency lemma} (IL, for short) when for every $n \in \mathbb{Z}^+$ there exists a finite set $\thicksim_n\!\!(x_1, \dots, x_n) \subseteq Fm(\vdash)$ such that
	\[
	\Gamma \cup \{\varphi_1,\dots,\varphi_n\} \text{ is inconsistent}\, \, \text{ iff } \, \, \Gamma \vdash\?\? \thicksim_n\!\!(\varphi_1,\dots,\varphi_n),
	\]
for every finite $\Gamma \cup \{\varphi_1,\dots,\varphi_1\} \subseteq Fm(\vdash)$;
\item\label{item:def-metalogical-2} the \emph{deduction theorem} (DT, for short) when for every $n, m \in \mathbb{Z}^+$ there exists a finite set $(x_1, \dots, x_n)\!\Rightarrow_{nm}\!\!(y_1, \dots, y_m) \subseteq Fm(\vdash)$ such that
\[
\Gamma, \varphi_1,\dots,\varphi_n \vdash \psi_1,\dots,\psi_m \, \, \text{ iff } \, \, \Gamma \vdash \?\?(\varphi_1, \dots, \varphi_n)\!\Rightarrow_{nm}\!\!(\psi_1, \dots, \psi_m),
\]
for every finite $\Gamma \cup \{\varphi_1,\dots,\varphi_n,\psi_1,\dots,\psi_m\} \sub Fm(\vdash)$;
\item\label{item:def-metalogical-3} the \emph{proof by cases} (PC, for short) when for every $n, m \in \mathbb{Z}^+$ there exists a finite set \linebreak $(x_1, \dots, x_n)\bigcurlyvee_{nm} (y_1, \dots, y_m) \subseteq Fm(\vdash)$ such that
\[
\Gamma, \varphi_1,\dots,\varphi_n \vdash \gamma \text{ and } \Gamma, \psi_1,\dots,\psi_m \vdash \gamma \, \, \text{ iff } \, \, \Gamma, (\varphi_1, \dots, \varphi_n)\!\bigcurlyvee_{nm} (\psi_1, \dots, \psi_m) \vdash \gamma,
\]
 for every finite $\Gamma \cup \{\varphi_1,\dots,\varphi_n,\psi_1,\dots,\psi_m, \gamma\} \sub Fm(\vdash)$.\footnote{For Conditions (\ref{item:def-metalogical-2}) and (\ref{item:def-metalogical-3}) to hold, the existence of a set with the desired property for $n = m = 1$ suffices. Our slightly redundant formulation, however, allows to simplify the presentation.}
\eroman
\end{Definition}

\begin{exa}\label{exa:IPC-metalogical} The logic $\mathsf{IPC}$ has the IL, the DT, and the PC witnessed, respectively, by the sets
\begin{align*}
\sim_n &\coloneqq \{ x_1 \to ( x_2 \to ( \dots ( x_n \to 0)\dots ) ) \};\\
\Rightarrow_{nm} &\coloneqq \{ x_1 \to ( x_2 \to ( \dots ( x_n \to y_k)\dots ) ) : k \leq m \};\\
\bigcurlyvee_{nm} &\coloneqq \{ x_i \lor y_j : i \leq n \text{ and } j \leq m \},
\end{align*}
for every $n,m \in \mathbb{Z}^+$.
\qed
\end{exa}

Henceforth, we will focus on the following class of logics \cite{BP86,Cz85,Cz86,Cz01}:

\begin{Definition}
A logic $\vdash$ is said to be \emph{protoalgebraic} if there exists a nonempty\footnote{The set $\Delta(x, y)$ is often allowed to be empty. However, the only protoalgebraic logic in a given algebraic language for which $\Delta(x, y)$ cannot be taken nonempty is the so-called \emph{almost inconsistent}, i.e., the logic $\vdash$ defined for every $\Gamma \cup \{ \varphi \} \subseteq Fm(\vdash)$ as follows: $\Gamma \vdash \varphi$ iff $\Gamma \ne \emptyset$ \cite[Prop.\ 6.11.4]{AAL-AIT-f}.} finite set $\Delta(x, y)$ of formulas such that
\[
\emptyset \vdash \Delta(x, x) \, \, \text{ and } \, \, x, \Delta(x, y) \vdash y.
\]
\end{Definition}

The class of protoalgebraic logics embraces most of the traditional logics. This is because if a logic $\vdash$ possesses a term-definable connective $\to$ such that $\emptyset \vdash x \to x$ and $x, x \to y \vdash y$ (as it is the case for $\mathsf{IPC}$), then it is protoalgebraic, as witnessed by the set $\Delta \coloneqq \{ x \to y \}$.

For protoalgebraic logics $\vdash$, the IL, the DT, and the PC admit a transparent description in terms of the structure of the semilattices $\Fic{\A}$ of compact deductive filters (see Remark \ref{Rem:FiC-structure}). More precisely, we have the following:

\begin{Theorem}[\protect{\cite[Thm.\ 3.7]{JGR13}}]\label{thm : il}
	Let $\vdash$ be a protoalgebraic logic. The following conditions are equivalent:
\benroman
		\item The logic $\vdash$ has the inconsistency lemma;
		\item The semilattice $\Thc{\vdash}$ is pseudocomplemented;
		\item \label{item:IL-3-proto} The semilattice $\Fic{\boldsymbol{A}}$ is pseudocomplemented, for every algebra $\A$. 
\eroman	Furthermore, if the inconsistency lemma for $\vdash$ is witnessed by $\{ \thicksim_n \? : n \in \mathbb{Z}^+\}$, then the operation $\lnot$ of the pseudocomplemented semilattice $\Fic{\boldsymbol{A}}$ is defined as follows: for every $a_1, \dots, a_n \in A$, 
	\[
	\lnot \fg{a_1, \dots, a_n} = \fg{\thicksim^{\boldsymbol{A}}_n \!\!(a_1, \dots, a_n)}.\footnote{Recall from Proposition \ref{prop : finitarity implies algebraicity} that the elements of $\Fic{\A}$ are precisely the finitely generated deductive filters of $\vdash$ on $\A$.}
	\]
\end{Theorem}

The next result is \cite[Thm.\ 2.11]{Cz85} (see also \cite{BP?,BP91, BP-AAL-DDT}).

\begin{Theorem}\label{thm : dt}
	Let $\vdash$ be a protoalgebraic logic. The following conditions are equivalent:
\benroman
		\item The logic $\vdash$ has the deduction theorem;
		\item The semilattice $\Thc{\vdash}$ is implicative;
		\item The semilattice $\Fic{\boldsymbol{A}}$ is implicative, for every algebra $\A$. 
\eroman	Furthermore, if the deduction theorem for $\vdash$ is witnessed by $\{ \Rightarrow_{nm} \? : n,m \in \mathbb{Z}^+\}$, then the operation $\to$ of the implicative semilattice $\Fic{\boldsymbol{A}}$ is defined as follows: for every $a_1, \dots, a_n, b_1, \dots, b_m \in A$, 
	\[
	\fg{a_1, \dots, a_n} \to \fg{b_1, \dots, b_m} = \fg{(a_1, \dots, a_n)\!\Rightarrow_{nm}^{\boldsymbol{A}} \!\!(b_1, \dots, b_m)}.
	\]
\end{Theorem}

Lastly, the next result originates in \cite{Cz84a} (see also \cite{CzDz} and \cite[Sec.\ 2.5]{Cz01}).

\begin{Theorem}\label{thm : pc}
	Let $\vdash$ be a protoalgebraic logic. The following conditions are equivalent: 
\benroman
		\item The logic $\vdash$ has the proof by cases;
		\item The semilattice $\Thc{\vdash}$ is a distributive lattice;
		\item The semilattice $\Fic{\boldsymbol{A}}$ is a distributive lattice, for every algebra $\A$. 
\eroman	
In this case, the lattice structure of the poset associated with the semilattice $\Fic{\boldsymbol{A}}$ is $\langle \Fic{\A}; +^\A, \cap \rangle$.\footnote{The meet operation of the lattice $\langle \Fic{\A}; +^\A, \cap \rangle$ is $+^\A$ and its join operation $\cap$. This is because the partial order associated with the semilattice $\Fic{\boldsymbol{A}}$ is the superset relation, as opposed to the inclusion relation (see Remark \ref{Rem:FiC-structure}, if necessary).} Furthermore, if the proof by cases for $\vdash$ is witnessed by $\{ \bigcurlyvee_{nm} \? : n,m \in \mathbb{Z}^+\}$, then for every $a_1, \dots, a_n, \linebreak b_1, \dots, b_m \in A$, 
	\[
	\fg{a_1, \dots, a_n} \cap \fg{b_1, \dots, b_m} = \fg{(a_1, \dots, a_n)\!\bigcurlyvee_{nm}^{\boldsymbol{A}} (b_1, \dots, b_m)}.
	\]
\end{Theorem}

\begin{Remark}
As we mentioned, $\mathsf{IPC}$ is protoalgebraic and it has the IL, the DT, and the PC. Therefore Condition (iii) of Theorems \ref{thm : il}, \ref{thm : dt}, and \ref{thm : pc} holds for $\mathsf{IPC}$. This should not come as a surprise, at least in the case where $\A$ is a Heyting algebra. This is because, in view of Example \ref{Exa:compact-filters-IPC}, the poset underlying $\mathsf{Fi}_{\mathsf{IPC}}^\omega(\A)$ is isomorphic to the lattice order of the Heyting algebra $\A$, which is obviously pseudocomplemented, implicative, and distributive.
\qed
\end{Remark}

\section{Sahlqvist theory for protoalgebraic logics}

The aim of this section is to extend Sahlqvist theory to protoalgebraic logics.\ To this end, it is convenient to introduce some terminology.\ A formula $\varphi$ is said to be a \emph{theorem} of a logic $\vdash$ when $\emptyset \vdash \varphi$. 

\begin{Remark}\label{Rem:protoalgebraic-have-thms}
Every protoalgebraic logic $\vdash$ has a theorem $\top(x)$. To prove this, let $\Delta(x, y)$ be the set witnessing the protoalgebraicity of $\vdash$. Since $\Delta(x, y)$ is nonempty, we can choose a formula $\varphi(x, y) \in \Delta(x, y)$. The definition of a protoalgebraic logic guarantees that $\emptyset \vdash \Delta(x, x)$ and, therefore, that $\emptyset \vdash \varphi(x, x)$. Thus, setting $\top(x) \coloneqq \varphi(x, x)$, we are done.
\qed
\end{Remark}	
	
Recall that $\mathcal{L}$ is the algebraic language of $\mathsf{IPC}$.

\begin{Definition}
	A formula $\varphi$ of $\mathcal{L}$ is \emph{compatible} with a protoalgebraic logic $\vdash$ when
	\benroman
	\item If $0$ or $\lnot$ occurs in $\varphi$, then $\vdash$ has the IL;
	\item If $\to$ occurs in $\varphi$, then $\vdash$ has the DT;
	\item If $\lor$ occurs in $\varphi$, then $\vdash$ has the PC.
	\eroman
\end{Definition}

\begin{Remark}\label{Rem:IPC-all-is-compatible}
Every formula of $\mathcal{L}$ is compatible with $\mathsf{IPC}$, because $\mathsf{IPC}$ has the IL, the DT, and the PC in view of Example \ref{exa:IPC-metalogical}. 
\qed
\end{Remark}

\begin{Remark}\label{Rem:induced-formula-term}
Let $\varphi(x_1, \dots, x_n)$ be a formula of $\mathcal{L}$ compatible with a protoalgebraic logic $\vdash$. In view of Condition (iii) of Theorems \ref{thm : il}, \ref{thm : dt}, and \ref{thm : pc}, the formula $\varphi$ can be interpreted in the semilattice $\Fic{\A}$, thus inducing a term-function
\[
\varphi^{\Fic{\A}} \colon \Fic{\A}^n \to \Fic{\A}
\]
on every algebra $\A$.
\qed
\end{Remark}

If the logic $\vdash$ in the next definition has the IL (resp.\ the DT or the PC), we denote the finite sets of formulas witnessing this property by $\thicksim_n\!\!(x_1, \dots, x_n)$ (resp.\ $(x_1, \dots, x_n)\!\Rightarrow_{nm}\!\!(y_1, \dots, y_m)$ or $(x_1, \dots, x_n)\bigcurlyvee_{nm} (y_1, \dots, y_m)$).

\begin{Definition}
With every $n \in \mathbb{Z}^+$ and formula $\varphi(x_1, \dots, x_n)$ of $\mathcal{L}$ that is compatible with a protoalgebraic logic $\vdash$ and every $k \in \mathbb{Z^+}$ we will associate a finite set
	\[
	\boldsymbol{\varphi}^k(x_1^1, \dots, x_1^k, \dots, x_n^1, \dots, x_n^k)
	\]
	of formulas of $\vdash$.
	 The case where $\varphi$ is a variable $x_m$ or a constant is handled as follows:
		\begin{align*}
	\boldsymbol{x_m}^k \coloneqq \{ x_m^1, \dots, x_m^k \}\qquad \boldsymbol{1}^k \coloneqq \{ \top(x_1^1) \} \qquad \boldsymbol{0}^k \coloneqq \{ x_1^1 \} \?\?\cup \sim_1\!\!(x_1^1),
	\end{align*}
where $\top(x)$ is a theorem of $\vdash$ (see Remark \ref{Rem:protoalgebraic-have-thms}).\footnote{Notice that the definition of $\boldsymbol{0}^k$ involves the set of formulas $\thicksim_1\!\!(x_1)$ typical of the IL. This makes sense because if the formula $0$ of $\mathcal{L}$ is compatible with $\vdash$, then the logic $\vdash$ has the IL. A similar remark applies to Conditions (\ref{ite:2:phi-set-bold}), (\ref{ite:3:phi-set-bold}), and (\ref{ite:4:phi-set-bold}) of this definition. Furthermore, notice that $0$ is viewed as a formula $0(x_1, \dots, x_n)$, where $n$ is a positive integer, and, therefore, $\boldsymbol{0}^k$ is allowed to contain formulas in the variable $x_1^1$. A similar remark applies to the definition of $\boldsymbol{1}^k$.} When $\varphi$ is a complex formula, we proceed as follows:
\benroman
		\item If $\varphi = \psi \land \chi$, we set
		\[
		\boldsymbol{\varphi}^k \coloneqq \boldsymbol{\psi}^k \cup \boldsymbol{\chi}^k;
		\]
		\item\label{ite:2:phi-set-bold} If $\varphi = \lnot \psi$ and $\boldsymbol{\psi}^k = \{ \psi_1, \dots, \psi_m \}$, we set
		\[
		\boldsymbol{\varphi}^k \coloneqq \,\, \thicksim_m\!\!(\psi_1, \dots, \psi_m);
		\]
				\item\label{ite:3:phi-set-bold} If $\varphi = \psi \to \chi$, $\boldsymbol{\psi}^k = \{\psi_1,\dots,\psi_m\}$, and $\boldsymbol{\chi}^k = \{\chi_1,\dots,\chi_t\}$, we set
\[
\boldsymbol{\varphi}^k \coloneqq  (\psi_1,\dots,\psi_m)\Rightarrow_{mt} (\chi_1,\dots,\chi_t);
\]
		\item\label{ite:4:phi-set-bold} If $\varphi = \psi \lor \chi$, $\boldsymbol{\psi}^k = \{\psi_1,\dots,\psi_m\}$, and $\boldsymbol{\chi}^k = \{\chi_1,\dots,\chi_t\}$, we set
\[
\boldsymbol{\varphi}^k \coloneqq (\psi_1,\dots,\psi_m)\!\bigcurlyvee_{mt} (\chi_1,\dots,\chi_t).
\]		
\eroman
\end{Definition}

\begin{exa}\label{Exa:IPC-implications-translated}
In view of Remark \ref{Rem:IPC-all-is-compatible}, the formula $\varphi = x_1 \to x_2$ is compatible with $\mathsf{IPC}$. Therefore, we can associate with $\varphi$ and every $k \in \mathbb{Z}^+$ a finite set $\boldsymbol{\varphi}^k$ of formulas of $\mathsf{IPC}$. The construction of $\boldsymbol{\varphi}^k$ depends on the sets of formulas witnessing the DT for $\mathsf{IPC}$, described in Example \ref{exa:IPC-metalogical}. As a result, we obtain
\[
\pushQED{\qed}\boldsymbol{\varphi}^k = \{ x_1^1 \to (x_1^2 \to ( \dots ( x_1^k \to x_2^i)\dots)) : i \leq k \}.\qedhere \popQED
\] 
\end{exa}

The connection between $\varphi(x_1, \dots, x_n)$ and $\boldsymbol{\varphi}^k(x_1^1, \dots, x_1^k, \dots, x_n^1, \dots, x_n^k)$ is made apparent by the following observation, where the function
\[
\varphi^{\Fic{\boldsymbol{A}}} \colon \Fic{\A}^n \to \Fic{\A}
\]
should be interpreted as in Remark \ref{Rem:induced-formula-term}.

\begin{Lemma}\label{lemma : filter generation}
Let $\varphi(x_1, \dots, x_n)$ be a formula of $\mathcal{L}$ compatible with a protoalgebraic logic $\vdash$. For every algebra $\boldsymbol{A}$, every $k \in \mathbb{Z}^+$, and every $a_1^1, \dots, a_1^k, \dots, a_n^1, \dots, a_n^k \in A$, it holds
\[
\fg{\boldsymbol{\varphi}^{k\A}(a_1^1, \dots, a_1^k, \dots, a_n^1, \dots, a_n^k)} = \varphi^{\Fic{\boldsymbol{A}}}(\fg{a_1^1, \dots, a_1^k},\dots,\fg{a_n^1, \dots, a_n^k}).
\] 
\end{Lemma}
\begin{proof}
The proof proceeds by induction on the construction of $\varphi$. In the base case, $\varphi$ is either a variable $x_m$ or one of the constants $1$ and $0$. The case of $x_m$ follows immediately from the definition of $\boldsymbol{x_m}^k$. Therefore, we only detail the cases of $1$ and $0$. 

On the one hand, we have that
\begin{align*}
\fg{\boldsymbol{1}^{k\A}(a_1^1, \dots, a_1^k, \dots, a_n^1, \dots, a_n^k)} &= \fg{\top(a_1^1)}\\
& = 1^{\Fic{\A}}\\
&= 1^{\Fic{\A}}(\fg{a_1^1, \dots, a_1^k},\dots,\fg{a_n^1, \dots, a_n^k}).
\end{align*}
The first equality above holds by the definition of $\boldsymbol{1}^k$ and the third is straightforward. To prove the second, recall that $\top(x_1^1)$ is a theorem of $\vdash$ and, therefore, $\fg{\top(a_1^1)}$ is the least compact deductive filter of $\vdash$ on $\A$. As $\Fic{\A}$ is ordered under the superset relation (see Remark \ref{Rem:FiC-structure}), this implies that $\fg{\top^\A(a_1^1)}$ is the top element of the semilattice $\Fic{\A}$, that is, $\fg{\top^\A(a_1^1)} = 1^{\Fic{\A}}$.

On the other hand, we have that
\begin{align*}
\fg{\boldsymbol{0}^{k\A}(a_1^1, \dots, a_1^k, \dots, a_n^1, \dots, a_n^k)} &= \fg{\{a_1^1 \} \? \? \cup \thicksim_1\!\!(a_1^1)}\\
& = \fg{a_1} +^\A \fg{\thicksim_1\!\!(a_1^1)}\\
& = \fg{a_1} +^\A \lnot^{\Fic{\A}}\fg{a_1^1}\\
&= 0^{\Fic{\A}}\\
&= 0^{\Fic{\A}}(\fg{a_1^1, \dots, a_1^k},\dots,\fg{a_n^1, \dots, a_n^k}).
\end{align*}
The equalities above are justified as follows: the first holds by the definition of $\boldsymbol{0}^{k}$, the second by the definition of $+^\A$ and $\fg{-}$, the third by the last part of Theorem \ref{thm : il}, the fourth by the fact that, in view of Theorem \ref{thm : il}(\ref{item:IL-3-proto}), $\Fic{\A}$ is a pseudocomplemented semilattice and, therefore, $0$ is term-definable as $x \land \lnot x$, and the last one is straightforward.

In the inductive step, $\varphi$ is a complex formula. If $\varphi = \psi \land \chi$, we have that
\begin{align*}
& \,\,\fg{\boldsymbol{\varphi}^{k\A}(a_1^1, \dots, a_1^k, \dots, a_n^1, \dots, a_n^k)} \\
= & \,\,\fg{(\boldsymbol{\psi \land \chi})^{k\A}(a_1^1, \dots, a_1^k, \dots, a_n^1, \dots, a_n^k)} \\
=& \,\,\fg{\boldsymbol{\psi}^{k\A}(a_1^1, \dots, a_1^k, \dots, a_n^1, \dots, a_n^k) \cup\boldsymbol{\chi}^{k\A}(a_1^1, \dots, a_1^k, \dots, a_n^1, \dots, a_n^k) }\\
=& \,\,\fg{\boldsymbol{\psi}^{k\A}(a_1^1, \dots, a_1^k, \dots, a_n^1, \dots, a_n^k)} +^\A \fg{\boldsymbol{\chi}^{k\A}(a_1^1, \dots, a_1^k, \dots, a_n^1, \dots, a_n^k) }\\
=& \,\,\psi^{\Fic{\boldsymbol{A}}}(\fg{a_1^1, \dots, a_1^k},\dots,\fg{a_n^1, \dots, a_n^k}) +^\A\chi^{\Fic{\boldsymbol{A}}}(\fg{a_1^1, \dots, a_1^k},\dots,\fg{a_n^1, \dots, a_n^k})\\
=& \,\,\varphi^{\Fic{\boldsymbol{A}}}(\fg{a_1^1, \dots, a_1^k},\dots,\fg{a_n^1, \dots, a_n^k}).
\end{align*}
The first equality above holds because $\varphi = \psi \land \chi$, the second by the definition of $(\boldsymbol{\psi \land \chi})^{k}$, the third by the definition of $+^\A$ and $\fg{-}$, the fourth by the inductive hypothesis, and the last one because $\varphi = \psi \land \chi$ and $+^\A$ is the operation of the semilattice $\Fic{\A}$.

It only remains to consider the cases where $\varphi$ is of the form $\lnot \psi$, $\psi \to \chi$, or $\psi \lor \chi$. Since they are handled essentially in the same way, we only detail the case where $\varphi = \lnot \psi$. Suppose that $\boldsymbol{\psi}^k = \{ \chi_1, \dots, \chi_m \}$. Then we have that
\begin{align*}
& \,\,\fg{\boldsymbol{\varphi}^{k\A}(a_1^1, \dots, a_1^k, \dots, a_n^1, \dots, a_n^k)} \\
=& \, \, \fg{(\boldsymbol{\lnot \psi})^{k\A}(a_1^1, \dots, a_1^k, \dots, a_n^1, \dots, a_n^k)}\\
= & \, \,  \fg{\sim_m\!\!(\chi_1^\A(a_1^1, \dots, a_1^k, \dots, a_n^1, \dots, a_n^k), \dots, \chi_m^\A(a_1^1, \dots, a_1^k, \dots, a_n^1, \dots, a_n^k))}
\\
= & \, \,  \lnot^{\Fic{\A}}\fg{\chi_1^\A(a_1^1, \dots, a_1^k, \dots, a_n^1, \dots, a_n^k), \dots, \chi_m^\A(a_1^1, \dots, a_1^k, \dots, a_n^1, \dots, a_n^k)}\\
= & \, \,  \lnot^{\Fic{\A}}\fg{\boldsymbol{\psi}^{k\A}(a_1^1, \dots, a_1^k, \dots, a_n^1, \dots, a_n^k)}\\
= & \, \, \lnot^{\Fic{\A}} \psi^{\Fic{\A}}(\fg{a_1^1, \dots, a_1^k}, \dots, \fg{a_n^1, \dots, a_n^k})\\
= & \, \, \varphi^{\Fic{\A}} (\fg{a_1^1, \dots, a_1^k}, \dots, \fg{a_n^1, \dots, a_n^k}).
\end{align*}
The first equality above holds because $\varphi = \lnot \psi$, the second by the definition of $(\boldsymbol{\lnot \psi})^k$ and the assumption that $\boldsymbol{\psi}^k = \{ \chi_1, \dots, \chi_m \}$, the third by the last part of Theorem \ref{thm : il}, the fourth by the assumption that $\boldsymbol{\psi}^k = \{ \chi_1, \dots, \chi_m \}$, and the fifth by the inductive hypothesis, and the last one because $\varphi = \lnot \psi$.
\end{proof}

The notion of compatibility can be extended to Sahlqvist quasiequations as follows:

\begin{Definition}
A Sahlqvist quasiequation $(\varphi_1 \land y \leq z \, \& \dots \& \, \varphi_m \land y \leq z) \Longrightarrow y \leq z$ is said to be \emph{compatible} with a protoalgebraic logic $\vdash$ if so are $\varphi_1, \dots, \varphi_m$.
\end{Definition}

In the case of protoalgebraic logics, the role of Sahlqvist quasiequations is played by the following metarules:

\begin{Definition}
Given a Sahlqvist quasiequation 
\[
\Phi = \varphi_1(x_1, \dots, x_n) \land y \leq z \, \& \dots \& \, \varphi_m(x_1, \dots, x_n) \land y \leq z \Longrightarrow y \leq z
\]
compatible with a protoalgebraic logic $\vdash$, let $\mathsf{R}_\vdash(\Phi)$ be the set of metarules of the form
\begin{prooftree}
	\AxiomC{$\Gamma, \boldsymbol{\varphi_1}^k(\gamma_1^1, \dots, \gamma_1^k, \dots, \gamma_n^1, \dots, \gamma_n^k) \rhd \psi$} 
	\AxiomC{$\dots$}
	\AxiomC{$\Gamma, \boldsymbol{\varphi_m}^k(\gamma_1^1, \dots, \gamma_1^k, \dots, \gamma_n^1, \dots, \gamma_n^k) \rhd \psi$}
	\TrinaryInfC{$\Gamma \rhd \psi$}
\end{prooftree}
where $k \in \mathbb{Z}^+$ and $\Gamma \cup \{ \psi \} \cup \{ \gamma_i^j : i \leq n, j \leq k \}$ is a finite subset of $Fm(\vdash)$.
\end{Definition}

We rely on the following observation, which generalizes \cite[Thm.\ 5.3]{LMR22}.

\begin{Proposition}\label{prop : transfer}
The following conditions are equivalent for a Sahlqvist quasiequation $\Phi$ compatible with a protoalgebraic logic $\vdash$:
\benroman
	\item\label{one1} The logic $\vdash$ validates the metarules in $\mathsf{R}_\vdash(\Phi)$;
	\item\label{two2} The semilattice $\Thc{\vdash}$ validates $\Phi$;
	\item\label{three3} The semilattice $\Fic{\boldsymbol{A}}$ validates $\Phi$, for every algebra $\boldsymbol{A}$.
\eroman
\end{Proposition}

The proof of Proposition \ref{prop : transfer} depends on the next well-known property of protoalgebraic logics.

\begin{Proposition}[\protect{\cite[Prop.\ 6.12]{AAL-AIT-f}}]\label{prop : protoalgebraic filters}
	Let $\vdash$ be a protoalgebraic logic, $\boldsymbol{A}$ an
	algebra, and $X \cup \{a\} \sub A$. Then $a \in \fg{X}$ iff there exist a finite $\Gamma \cup \{\varphi\} \sub Fm(\vdash)$ and a homomorphism
	$f \colon \boldsymbol{Fm}(\vdash) \to \boldsymbol{A}$ such that $\Gamma \vdash \varphi$, $f[\Gamma] \sub X \cup \fg{\emptyset}$, and $f(\varphi) = a$. 
\end{Proposition}

\begin{proof}[Proof of Proposition \ref{prop : transfer}]
Throughout the proof we will assume that
\[
\Phi = \varphi_1(x_1, \dots, x_n) \land y \leq z \, \& \dots \& \, \varphi_m(x_1, \dots, x_n) \land y \leq z \Longrightarrow y \leq z.
\]

(\ref{two2})$\Rightarrow$(\ref{one1}): Let $k \in \mathbb{Z}^+$ and let $\Gamma \cup \{ \psi \} \cup \{ \gamma_i^j : i \leq n, j \leq k \}$ be a finite subset of $Fm(\vdash)$ such that
\begin{equation}\label{Eq:qkfjwefojorpimpo}
\Gamma  \, {\cup}  \, \boldsymbol{\varphi_i}^k(\gamma_1^1, \dots, \gamma_1^k, \dots, \gamma_n^1, \dots, \gamma_n^k) \vdash \psi, \text{ for every } i \leq m.
\end{equation}
We want to prove that $\Gamma \vdash \psi$.

Consider $i \leq m$. We have that
\begin{align*}
\fgt{\psi} &\subseteq \fgt{\Gamma \cup \boldsymbol{\varphi_i}^k(\gamma_1^1, \dots, \gamma_1^k, \dots, \gamma_n^1, \dots, \gamma_n^k)}\\
 &=\fgt{\Gamma} + \fgt{\boldsymbol{\varphi_i}^k(\gamma_1^1, \dots, \gamma_1^k, \dots, \gamma_n^1, \dots, \gamma_n^k)}\\
&=\fgt{\Gamma} + \varphi_i^{\Thc{\vdash}}(\fgt{\gamma_1^1, \dots, \gamma_1^k}, \dots, \fgt{\gamma_n^1, \dots, \gamma_n^k}),
\end{align*}
where the first step follows from Condition (\ref{Eq:qkfjwefojorpimpo}), the second from the definition of $\fgt{-}$ and $+$, and the last one from Lemma \ref{lemma : filter generation}.

Since the semilattice $\Thc{\vdash}$ is ordered under the superset relation and its operation is $+$ (see Remark \ref{Rem:FiC-structure}), the above display yields
\[
\fgt{\Gamma} \land^{\Thc{\vdash}}\varphi_i^{\Thc{\vdash}}(\fgt{\gamma_1^1, \dots, \gamma_1^k}, \dots, \fgt{\gamma_n^1, \dots, \gamma_n^k}) \leq \fgt{\psi}.
\]
As this holds for every $i \leq m$, we can apply the assumption that $\Thc{\vdash}$ validates $\Phi$, obtaining $\fgt{\Gamma} \leq \fgt{\psi}$. But, since the order of $\Thc{\vdash}$ is the superset relation, this amounts to $\fgt{\psi} \subseteq \fgt{\Gamma}$, whence $\Gamma \vdash \psi$ as desired.

(\ref{three3})$\Rightarrow$(\ref{two2}): This implication is straightforward, since $\Thc{\vdash} = \Fic{\boldsymbol{Fm}(\vdash)}$.

(\ref{one1})$\Rightarrow$(\ref{three3}): Let $\A$ be an algebra and let $F, G_1,\dots,G_{n}, H\in \Fic{\A}$ be such that 
\begin{equation}\label{Eq:fdkjhfdpppp}
H \subseteq F +^\A \varphi_{i}^{\Fic{\boldsymbol{A}}}(G_1,\dots,G_{n}) \text{, for every }i \leq m.
\end{equation}
We want to prove that $H \subseteq F$.

Recall that $\vdash$ has theorems, because it is protoalgebraic (see Remark \ref{Rem:protoalgebraic-have-thms}). Consequently, $G_1, \dots, G_n$ are nonempty. Furthermore, they are finitely generated. This is because, in view of Proposition \ref{prop : finitarity implies algebraicity}, compact and finitely generated deductive filters coincide. Therefore, there are $k \in \mathbb{Z}^+$ and $a_1^1, \dots, a_1^k, \dots, a_n^1, \dots, a_n^k \in A$ such that
\[
G_1 = \fg{a_1^1, \dots, a_1^k}, \dots, G_n = \fg{a_n^1, \dots, a_n^k}.
\]
Together with Lemma \ref{lemma : filter generation}, this implies that
\begin{equation}\label{Eq:fdkjhfdpppp2}
\varphi_{i}^{\Fic{\boldsymbol{A}}}(G_1,\dots,G_{n}) = \fg{\boldsymbol{\varphi_i}^k(a_1^1, \dots, a_1^k, \dots, a_n^1, \dots, a_n^k)}, \text{ for every }i \leq m.
\end{equation}

In order to prove that $H \subseteq F$, let $a \in H$.\ From Conditions (\ref{Eq:fdkjhfdpppp}) and (\ref{Eq:fdkjhfdpppp2}) it follows that
\[
a \in F +^\A \fg{\boldsymbol{\varphi_i}^k(a_1^1, \dots, a_1^k, \dots, a_n^1, \dots, a_n^k)}, \text{ for every }i \leq m.
\]
Thus, from Proposition \ref{prop : protoalgebraic filters} we deduce that for each $i \leq m$ there exists a homomorphism $f_i \colon \boldsymbol{Fm}(\vdash) \to \boldsymbol{A}$ and a finite set $\Psi_i \cup \{\psi_i\} \sub Fm(\vdash)$ such that 
\begin{equation}\label{Eq:gjdpqlvplqeWW}
\Psi_i \vdash \psi_i, \, \, \,  f_i[\Psi_i] \sub F \cup \boldsymbol{\varphi_i}^{k\boldsymbol{A}}(a_1^1, \dots, a_1^k, \dots, a_n^1, \dots, a_n^k), \, \, \text{ and } \,\, f_i(\psi_i) = a.
\end{equation}

As the sets $\Psi_1, \dots, \Psi_m$ are finite and $\vdash$ is substitution invariant, we may assume, without loss of generality, that if $i < j \leq m$, then the set of variables occurring in the members of $\Psi_i \cup \{\psi_i\}$ is disjoint from the set of variables occurring in the members of $\Psi_j \cup \{\psi_j\}$. Consequently, we may also assume that $f_1 =  \dots = f_m$. Accordingly, from now on, we will denote these maps by $f$ and drop the subscripts. Lastly, we may assume that there exists a set of variables $\{ z_j^t : j \leq n, t \leq k \}$ not occurring in any $\Psi_i \cup \{ \psi_i \}$ such that $f(z_j^t) = a_j^t$ for every $j \leq n$ and $t \leq k$.

Now, in view of Condition (\ref{Eq:gjdpqlvplqeWW}), we can split each $\Psi_i$ into two subsets $\Psi_i^1$ and $\Psi_i^2$ such that 
\begin{equation}\label{Eq:gjdpqlvplqeWWsdwd}
f[\Psi_i^1] \sub F, \quad f[\Psi_i^2] \sub \boldsymbol{\varphi_i}^{k\boldsymbol{A}}(a_1^1, \dots, a_1^k, \dots, a_n^1, \dots, a_n^k), \,\, \text{ and } \, \,\Psi_i = \Psi_i^1 \cup \Psi_i^1.
\end{equation}
We will construct a finite set $\Gamma \sub f^{-1}[F]$ as follows. First, we stipulate that $\Psi_1^1 \cup \dots \cup \Psi_m^1  \sub \Gamma$, as the above display guarantees that $\Psi_1^1 \cup \dots \cup \Psi_m^1 \subseteq f^{-1}[F]$.

 Then let $\Delta(x, y)$ be the finite set of formulas witnessing the protoalgebraicity of $\vdash$. We will make extensive use of the observation that, since $\emptyset \vdash \Delta(x, x)$ and $F$ is a deductive filter, we have that $\Delta^{\A}(b, b) \subseteq F$ for every $b \in F$. 

Recall from Condition (\ref{Eq:gjdpqlvplqeWW}) that $f(\psi_1) = \dots = f(\psi_m) = a$. Therefore,
\[
f[\bigcup_{i, j \leq m}\Delta(\psi_i, \psi_j)]  = \Delta^\A(a, a) \subseteq F
\]
and so we may assume that $\Gamma$ contains $\bigcup\{ \Delta(\psi_i, \psi_j) : i, j \leq m \}$. Moreover, by Condition (\ref{Eq:gjdpqlvplqeWWsdwd}), we have that
\[
f[\Psi_i^2] \sub \boldsymbol{\varphi_i}^{k\boldsymbol{A}}(a_1^1, \dots, a_1^k, \dots, a_n^1, \dots, a_n^k)\text{, for every }i \leq m.
\]
Accordingly, for every $i \leq m$ and $\alpha \in \Psi_i^2$, there exists a formula $\beta_\alpha \in \boldsymbol{\varphi_i}^{k}(z_1^1, \dots, z_1^k, \dots, z_n^1, \dots, z_n^k)$ such that
\begin{align*}
f(\alpha) &= \beta_\alpha(a_1^1, \dots, a_1^k, \dots, a_n^1, \dots, a_n^k)\\
&= \beta_\alpha(f(z_1^1), \dots, f(z_1^k), \dots, f(z_n^1), \dots, f(z_n^k))\\
& = f(\beta_\alpha)
\end{align*}
and, therefore, $f[\Delta(\beta_\alpha, \alpha)] \subseteq F$. Consequently, we can add the sets $\Delta(\beta_\alpha, \alpha)$ to $\Gamma$, thereby completing its definition.

To conclude the proof, it suffices to show that
\begin{equation}\label{Eq:cksckcsnkcsklnscnlkslnk}
\Gamma \cup \boldsymbol{\varphi_i}^k(z_1^1, \dots, z_1^k, \dots, z_n^1, \dots, z_n^k) \vdash \psi_1, \text{ for every }i \leq m. 
\end{equation}
For if this is the case, the assumption that $\vdash$ validates the rules in $\mathsf{R}_\vdash(\Phi)$ implies that $\Gamma \vdash \psi_1$. Moreover, since $f[\Gamma] \sub F$ and $F$ is a deductive filter, we deduce $a = f(\psi_1) \in F$ as desired. 

Accordingly, we turn to prove Condition (\ref{Eq:cksckcsnkcsklnscnlkslnk}). Consider $i \leq m$. First, observe that
\begin{equation}\label{Eq:final-lkjdokpppwwer1}
\Gamma \, {\cup} \, \boldsymbol{\varphi_i}^k(z_1^1, \dots, z_1^k, \dots, z_n^1, \dots, z_n^k) \vdash \Psi_i^1,
\end{equation}
because $\Psi_i^1 \sub \Gamma$. We will prove that
\begin{equation}\label{Eq:final-lkjdokpppwwer2}
\Gamma \cup \boldsymbol{\varphi_i}^k(z_1^1, \dots, z_1^k, \dots, z_n^1, \dots, z_n^k) \vdash \Psi_i^2.
\end{equation}
To this end, consider a formula $\alpha \in \Psi_i^2$. By the construction of $\Gamma$, we have
\[
\beta_\alpha \in \boldsymbol{\varphi_i}^k(z_1^1, \dots, z_1^k, \dots, z_n^1, \dots, z_n^k) \, \, \text{ and } \, \, \Delta(\beta_\alpha,\alpha) \sub \Gamma.
\]
As the definition of a protoalgebraic logic gives $\beta_\alpha, \Delta(\beta_\alpha, \alpha) \vdash \alpha$, the above display guarantees that $\Gamma \cup \boldsymbol{\varphi_i}^k(z_1^1, \dots, z_1^k, \dots, z_n^1, \dots, z_n^k) \vdash \alpha$, thereby establishing Condition (\ref{Eq:final-lkjdokpppwwer2}).

Now, recall that $\Psi_i = \Psi_i^1 \cup \Psi_i^2$ and $\Psi_i \vdash \psi_i$ (see Conditions (\ref{Eq:gjdpqlvplqeWW}) and (\ref{Eq:gjdpqlvplqeWWsdwd}), if necessary). Together with the Conditions (\ref{Eq:final-lkjdokpppwwer1}) and (\ref{Eq:final-lkjdokpppwwer2}), this yields
\[
\Gamma \cup \boldsymbol{\varphi_i}^k(z_1^1, \dots, z_1^k, \dots, z_n^1, \dots, z_n^k) \vdash \psi_i.
\]
Finally, since by the construction of $\Gamma$ we have $\Delta(\psi_i,\psi_1) \sub \Gamma$ and, by protoalgebraicity, $\psi_i, \Delta(\psi_i, \psi_1) \vdash \psi_1$, we conclude that $\Gamma \vdash \psi_1$.
\end{proof}

\begin{Remark}\label{Rem:DT-Conjunction-simplifications}
A logic $\vdash$ is said to have a \emph{conjunction} if it possesses a term-definable binary connective $\land$ such that
\[
x, y \vdash x \land y \qquad x \land y \vdash x \qquad x \land y \vdash y.
\]
In this case, for every algebra $\A$ and $a_1, \dots, a_n \in A$,
\[
\fg{a_1, \dots, a_n} = \fg{ a_1 \land \dots \land a_n }.
\]
Consequently, the members of $\Fic{\A}$ are precisely the \emph{principal} deductive filters of $\vdash$ on $\A$, i.e., the sets of the form $\fg{a}$ for some $a \in A$.

As a consequence, if the logic $\vdash$ in the statement of Proposition \ref{prop : transfer} has a conjunction, then the positive integer $k$ in the proof of the implication (\ref{one1})$\Rightarrow$(\ref{three3}) can be taken to be $1$. Accordingly, for logics $\vdash$ with a conjunction, Condition (\ref{one1}) of Proposition \ref{prop : transfer} can be replaced by the simpler demand that $\vdash$ validates the metarules of the form

\begin{prooftree}
	\AxiomC{$\Gamma, \boldsymbol{\varphi_1}^1(\gamma_1, \dots, \gamma_n) \rhd \psi$} 
	\AxiomC{$\dots$}
	\AxiomC{$\Gamma, \boldsymbol{\varphi_m}^1(\gamma_1, \dots, \gamma_n) \rhd \psi$}
	\TrinaryInfC{$\Gamma \rhd \psi$}
\end{prooftree}
where $\Gamma \cup \{ \psi \} \cup \{ \gamma_1, \dots, \gamma_n \}$ is a finite subset of $Fm(\vdash)$.

 A similar simplification is possible when the logic $\vdash$ has the DT or the PC, as we proceed to explain. Suppose first that $\vdash$ has the DT. Given two finite subsets $\Gamma = \{ \varphi_1, \dots, \varphi_n \}$ and $\Sigma = \{ \psi_1, \dots, \psi_m \}$ of $Fm(\vdash)$, we will write 
\[
\Gamma \Rightarrow \Sigma\,\, \text{ as a shorthand for }\,\, (\varphi_1, \dots, \varphi_n)\Rightarrow_{nm}(\psi_1, \dots, \psi_m),
\]
where $\Rightarrow_{nm}$ is one of the sets witnessing the DT for $\vdash$. In the presence of the DT, Condition (\ref{one1}) of Proposition \ref{prop : transfer} becomes equivalent to the simpler demand that
\begin{equation}\label{Eq:IPC-simplication-Rules-k-1}
((\boldsymbol{\varphi_1}^k \Rightarrow y) \cup \dots \cup (\boldsymbol{\varphi_m}^k \Rightarrow y)) \Rightarrow y
\end{equation}
is a set of theorems of $\vdash$ for every $k \in \mathbb{Z}^+$, where $y$ is a variable that does not occur in any $\boldsymbol{\varphi_i}^k = \boldsymbol{\varphi_i}^k(x_1^1, \dots, x_1^k, \dots, x_n^1, \dots, x_n^k)$. We leave the easy proof, which relies only on the basic properties of the DT, to the reader.

Lastly, we turn to the case where $\vdash$ has the PC. Given finite subsets $\Gamma_1 = \{ \varphi_1, \dots, \varphi_{n_1} \}, \dots, \Gamma_m = \{ \varphi_1, \dots, \varphi_{n_m} \}$ of $Fm(\vdash)$, we will define recursively a finite set
\[
\Gamma_1 \! \bigcurlyvee \dots \bigcurlyvee \Gamma_p
\]
of formulas, for every $2 \leq p \leq m$. First, if $p =2$, we let
\[
\Gamma_1 \! \bigcurlyvee \Gamma_2 \coloneqq (\varphi_1, \dots, \varphi_{n_1})\!\bigcurlyvee_{n_1 n_2} (\varphi_1, \dots, \varphi_{n_2}),
\]
where $\bigcurlyvee_{n_1n_2}$ is one of the sets witnessing the PC for $\vdash$. On the other hand, if $2 < p < m$ and $\Gamma_1 \! \bigcurlyvee \dots \bigcurlyvee \Gamma_p = \{ \gamma_1, \dots, \gamma_t \}$, we let
\[
\Gamma_1 \! \bigcurlyvee \dots \bigcurlyvee \Gamma_{p+1} \coloneqq (\gamma_1, \dots, \gamma_t)\!\!\bigcurlyvee_{t n_{k+1}} (\varphi_1, \dots, \varphi_{n_{p+1}}).
\]

In the presence of the PC, Condition (\ref{one1}) of Proposition \ref{prop : transfer} becomes equivalent to the simpler demand that
\begin{equation}\label{Eq:IPC-PC-simplication-Rules-k-1}
\boldsymbol{\varphi_1}^k \! \bigcurlyvee \dots \bigcurlyvee\boldsymbol{\varphi_m}^k
\end{equation}
is a set of theorems of $\vdash$ for every $k \in \mathbb{Z}^+$.\footnote{If the logic $\vdash$ has a conjunction, we can restrict to the case where $k = 1$ both in Conditions (\ref{Eq:IPC-simplication-Rules-k-1}) and (\ref{Eq:IPC-PC-simplication-Rules-k-1}).} Also in this case, we leave the easy proof, which relies only on the basic properties of the PC, to the reader.
\qed
\end{Remark}

Sahlqvist theory for protoalgebraic logics centers on the following notion.

\begin{Definition}
Let $\vdash$ be a logic and $\boldsymbol{A}$ an algebra. The \emph{spectrum of} $\A$ \emph{relative to} $\vdash$, in symbols $\mathsf{Spec}_{\vdash}(\A)$, is the poset of meet irreducible deductive filters of $\vdash$ on $\A$ ordered under the inclusion relation. When $\A = \boldsymbol{Fm}(\vdash)$, we write $\mathsf{Spec}(\vdash)$ as a shorthand for $\mathsf{Spec}_\vdash(\A)$.
\end{Definition}

\begin{Remark}
In view of Remark \ref{Rem:prime-filters} and Example \ref{exa-filers-ipc}, the spectrum of a Heyting algebra $\A$ relative to $\mathsf{IPC}$ is the poset of prime filters of $\A$.
\qed 
\end{Remark}

Our main result establishes a correspondence between the validity of the metarules of the form $\mathsf{R}_\vdash(\Phi)$ and the structure of spectra $\mathsf{Spec}_{\vdash}(\A)$.\footnote{While the Abstract Sahlqvist Theorem takes the form of a correspondence result, it can also be used to derive canonicity theorems, as shown in Theorem \ref{Thm:canonicity-arrow-IPC}.}

\begin{Abstract Sahlqvist Theorem}\label{thm : correspondence}
	The following conditions are equivalent for a Sahlqvist quasiequation $\Phi$ compatible with a protoalgebraic logic $\vdash$:
\benroman
		\item\label{1one} The logic $\vdash$ validates the metarules in $\mathsf{R}_\vdash(\Phi)$;
		\item\label{2two} $\mathsf{Spec}(\vdash) \vDash \mathsf{tr}(\Phi)$;
		\item\label{3three} $\mathsf{Spec}_\vdash(\boldsymbol{A}) \vDash \mathsf{tr}(\Phi)$, for every algebra $\boldsymbol{A}$.
\eroman
\end{Abstract Sahlqvist Theorem}

\begin{proof}
(\ref{1one})$\Rightarrow$(\ref{3three}): Let $\A$ be an algebra. By applying Proposition \ref{prop : transfer} to the assumption that $\vdash$ validates the metarules in $\mathsf{R}_\vdash(\Phi)$, we obtain that the semilattice $\Fic{\boldsymbol{A}}$ validates $\Phi$. 

Then let $\mathcal{L}_\land$ be the sublanguage of $\mathcal{L}$ consisting of the connectives of $\mathsf{IPC}$ that occur in $\Phi$ with the addition of $\land$. As $\Phi$ is compatible with $\vdash$, from Theorems \ref{thm : il}, \ref{thm : dt}, and \ref{thm : pc} it follows that $\Fic{\boldsymbol{A}}$ is an $\mathcal{L}_\land$-subreduct of a Heyting algebra. Furthermore, observe that $\Phi$ is a Sahlqvist quasiequation in $\mathcal{L}_\land$. Therefore, we can apply Corollary \ref{corollary : corollary of enhanced sahlqvist} obtaining that $\Fic{\boldsymbol{A}}_\ast \vDash \mathsf{tr}(\Phi)$. 

Now, recall from Proposition \ref{prop : finitarity implies algebraicity} that the lattice $\Fi{\boldsymbol{A}}$ is algebraic. Therefore, from Theorem \ref{thm : algebraic lattices and filters of compact} we deduce that $\Fi{\boldsymbol{A}}$ is isomorphic to the lattice of filters of the semilattice $\Fic{\boldsymbol{A}}$. Thus, the poset of meet irreducible elements of $\Fi{\boldsymbol{A}}$, namely $\mathsf{Spec}_\vdash(\boldsymbol{A})$, is isomorphic to the poset of meet irreducible filters of $\Fic{\boldsymbol{A}}$, namely $\Fic{\boldsymbol{A}}_\ast$. Consequently, from $\Fic{\boldsymbol{A}}_\ast \vDash \mathsf{tr}(\Phi)$ it follows that $\mathsf{Spec}_\vdash(\boldsymbol{A}) \vDash \mathsf{tr}(\Phi)$ as desired.

(\ref{3three})$\Rightarrow$(\ref{2two}): Straightforward.

(\ref{2two})$\Rightarrow$(\ref{1one}): Assume $\mathsf{Spec}(\vdash) \vDash \mathsf{tr}(\Phi)$.\ As in the proof of the implication (\ref{1one})$\Rightarrow$(\ref{3three}), we have $\mathsf{Spec}(\vdash) \cong \Thc{\vdash}_\ast$. Consequently, we obtain $\Thc{\vdash}_\ast\vDash \mathsf{tr}(\Phi)$. Now, let $\mathcal{L}_\land$ be the language defined in the proof of the implication (\ref{1one})$\Rightarrow$(\ref{3three}).\ The same argument shows that $\Thc{\vdash}$ is an $\mathcal{L}_\land$-subreduct of a Heyting algebra and that $\Phi$ is a Sahlqvist quasiequation in $\mathcal{L}_\land$. Therefore, we can apply Corollary \ref{corollary : corollary of enhanced sahlqvist} to $\Thc{\vdash}_\ast\vDash \mathsf{tr}(\Phi)$, obtaining that $\Thc{\vdash} \vDash \Phi$. Lastly, by Proposition \ref{prop : transfer} we conclude that $\vdash$ validates the metarules in $\mathsf{R}_\vdash(\Phi)$ as desired. 
\end{proof}

Under additional assumptions, the Abstract Sahlqvist Theorem can be formulated in a more algebraic fashion. Given a quasivariety $\mathsf{K}$ and an algebra $\A$, we say that a congruence $\theta$ of $\A$ is a $\mathsf{K}$-\emph{congruence} of $\A$ when $\A / \theta \in \mathsf{K}$. When ordered under the inclusion relation, the set of $\mathsf{K}$-congruences of $\A$ forms an algebraic lattice, which we denote by $\mathsf{Con}_{\mathsf{K}}(\A)$. The poset of meet irreducible elements of $\mathsf{Con}_{\mathsf{K}}(\A)$ will then be denoted by $\mathsf{Spec}_{\mathsf{K}}(\A)$.

A logic $\vdash$ is \emph{algebraized} \cite{BP89} by a quasivariety $\mathsf{K}$ when there exist finite sets $\Delta(x, y)$ and $\tau(x)$ of formulas and equations, respectively, such that
\[
\mathsf{K} \vDash x \thickapprox y \, \,\text{ iff } \, \,\{ \epsilon(\varphi) \thickapprox \delta(\varphi) : \epsilon \thickapprox \delta \in \tau \text{ and } \varphi \in \Delta \}
\]
and, for every finite $\Gamma \cup \{ \varphi \} \subseteq Fm(\vdash)$,
\[
\Gamma \vdash  \varphi \, \, \text{ iff } \, \,  \mathsf{K} \vDash \foo \{ \epsilon(\gamma) \thickapprox \delta(\gamma) : \gamma \in \Gamma, \epsilon \thickapprox \delta \in \tau\} \Longrightarrow \epsilon'(\varphi) \thickapprox \delta'(\varphi)\text{, for all }\epsilon' \thickapprox \delta' \in \tau.
\]
In this case, for every algebra $\A$, the lattices $\Fi{\A}$ and $\mathsf{Con}_\mathsf{K}(\A)$ are isomorphic (see, e.g., \cite[Thm.\ 3.58]{AAL-AIT-f}) and, therefore, so are $\mathsf{Spec}_\vdash(\A)$ and $\mathsf{Spec}_{\mathsf{K}}(\A)$. Furthermore, the set of formulas $\Delta(x, y)$ witnesses the protoalgebraicity of $\vdash$.

\begin{exa}
The intuitionistic propositional calculus $\mathsf{IPC}$ is algebraized by the variety $\mathsf{HA}$ of Heyting algebras, as witnessed by the sets $\Delta = \{ x \to y, y \to x \}$ and $\tau = \{ x \thickapprox 1 \}$.
\qed
\end{exa}

\begin{Corollary}\label{cor : correspondence - algebraizable}
Let $\Phi$ be a Sahlqvist quasiequation compatible with a logic $\vdash$ that is algebraized by a quasivariety $\mathsf{K}$. Then $\vdash$ validates the metarules in $\mathsf{R}_\vdash(\Phi)$ iff $\mathsf{Spec}_\mathsf{K}(\boldsymbol{A}) \vDash \mathsf{tr}(\Phi)$, for every $\boldsymbol{A} \in \mathsf{K}$.
\end{Corollary} 

\begin{proof}
The ``only if'' part follows from the implication (\ref{1one})$\Rightarrow$(\ref{3three}) of the Abstract Sahlqvist Theorem and the observation that $\mathsf{Spec}_\vdash(\A) \cong\mathsf{Spec}_{\mathsf{K}}(\A)$, for every algebra $\A$. To prove the ``if'' part, suppose that $\mathsf{Spec}_\mathsf{K}(\boldsymbol{A}) \vDash \mathsf{tr}(\Phi)$, for every $\boldsymbol{A} \in \mathsf{K}$. Then consider an algebra $\A$, not necessarily in $\mathsf{K}$. By the Correspondence Theorem for quasivarieties, there exists $\B \in \mathsf{K}$ such that $\mathsf{Spec}_\mathsf{K}(\boldsymbol{A}) \cong \mathsf{Spec}_\mathsf{K}(\boldsymbol{B})$ (see, e.g., \cite[II.6.20]{BuSa00}). Together with the assumption, this implies that $\mathsf{Spec}_\mathsf{K}(\boldsymbol{A})\vDash \mathsf{tr}(\Phi)$, thus establishing Condition (\ref{3three}) of the Abstract Sahlqvist Theorem. By the implication (\ref{3three})$\Rightarrow$(\ref{1one}) of the same theorem, we conclude that $\vdash$ validates the metarules in $\mathsf{R}_\vdash(\Phi)$ as desired.
\end{proof}

\section{The excluded middle and the bounded top width laws}

We proceed to illustrate how the Abstract Sahlqvist Theorem can be used to obtain concrete correspondence results, some known and some new.

\begin{Definition}
A logic $\vdash$ is said to have the \emph{excluded middle law} (EML, for short) when for every $n \in \mathbb{Z}^+$ there exists a finite set $\sim_n\!\!(x_1, \dots, x_n) \subseteq Fm(\vdash)$ such that
\[
\{ x_1, \dots, x_n \} \, \cup \sim_n\!\!(x_1, \dots, x_n) \text{ is inconsistent}
\]
and the metarule
\begin{prooftree}
	\AxiomC{$\Gamma, \varphi_1, \dots, \varphi_n \rhd \psi$} 
	\AxiomC{$\Gamma, \sim_n\!\!(\varphi_1, \dots, \varphi_n) \rhd \psi$}
	\BinaryInfC{$\Gamma \rhd \psi$}
	\end{prooftree}
is valid in $\vdash$, for every finite $\Gamma \cup \{ \varphi_1, \dots, \varphi_n, \psi \} \subseteq Fm(\vdash)$. 
\end{Definition}
\begin{Remark}\label{Rem:EML-implies-IL}
Every logic with the EML has the IL, as witnessed by the sets $\sim_n\!\!(x_1, \dots, x_n)$.
\qed
\end{Remark}

In the presence of the IL, the semantic counterpart of the EML is the following property:

\begin{Definition} A logic $\vdash$ is said to be \emph{semisimple} when the order of $\mathsf{Spec}_\vdash(\A)$ is the identity relation, for every algebra $\A$.
\end{Definition}

\begin{Theorem}[\protect{\cite{PrenLav20}}]\label{Thm:Pre-Lav-thm}
A protoalgebraic logic has the excluded middle law iff it has the inconsistent lemma and is semisimple. 
\end{Theorem}

\begin{proof}
In view of Remark \ref{Rem:EML-implies-IL}, it suffices to prove that a protoalgebraic logic $\vdash$ with the IL has the EML iff it is semisimple. Accordingly, let $\{ \sim_n\!\!(x_1, \dots, x_n) : n \in \mathbb{Z}^+ \}$ be a family of sets witnessing the IL for $\vdash$. Moreover, observe that the Sahlqvist quasiequation
\[
\Phi = x \land y \leq z \, \& \, \lnot x \land y \leq z \Longrightarrow y \leq z
\]
corresponding to the excluded middle axiom $x \lor \lnot x$ is compatible with $\vdash$, because $\vdash$ has the IL.

Now, recall from Example \ref{Exa:BTWL-correspondence-easy} that a poset validates $\mathsf{tr}(\Phi)$ iff its order is the identity relation. Consequently, $\vdash$ is semisimple iff $\mathsf{Spec}_{\vdash}(\A) \vDash \mathsf{tr}(\Phi)$, for every algebra $\A$. By the Abstract Sahlqvist Theorem, the latter condition is equivalent to the demand that $\vdash$ validates the metarules in $\mathsf{R}_\vdash
(\Phi)$, namely,
\begin{prooftree}
	\AxiomC{$\Gamma, \varphi_1, \dots, \varphi_n \rhd \psi$} 
	\AxiomC{$\Gamma, \sim_n\!\!(\varphi_1, \dots, \varphi_n) \rhd \psi$}
	\BinaryInfC{$\Gamma \rhd \psi$}
	\end{prooftree}
for every finite $\Gamma \cup \{ \varphi_1, \dots, \varphi_n, \psi \} \subseteq Fm(\vdash)$. But, since the IL guarantees that the sets of the form $\{ x_1, \dots, x_n \} \, \cup \sim_n\!\!(x_1, \dots, x_n)$ are inconsistent, this amounts to the demand that $\vdash$ has the EML.
\end{proof}

In order to derive a similar result for the bounded top width axioms, we adopt the following convention: if a family $\{ \sim_n\!\!(x_1, \dots, x_n) : n \in \mathbb{Z}^+ \}$ of sets of formulas witnesses the IL for a logic $\vdash$, then for every finite set $\Gamma = \{ \gamma_1, \dots, \gamma_n \}$ of formulas we will write 
\[
\sim\!\Gamma \,\,\text{ as a shorthand for }\,\, \sim_n\!\!(\gamma_1, \dots, \gamma_n).
\]

\begin{Definition}
Let $\vdash$ be a logic with the IL witnessed by a family $\{ \sim_m\!\!(x_1, \dots, x_m) : m \in \mathbb{Z}^+ \}$ and let $n \in \mathbb{Z}^+$. The logic $\vdash$ has the \emph{bounded top width $n$ law} ($\textup{BTWL}_n$, for short) if it validates the metarule
\begin{prooftree}
	\AxiomC{$\Gamma, \sim\!(\sim\!(\gamma_i^1, \dots, \gamma_i^k)  \cup \{ \gamma_{j}^t : j < i, t \leq k \})  \rhd \psi$} 
	\AxiomC{for every $i \leq n+1$}
	\BinaryInfC{$\Gamma \rhd \psi$}
	\end{prooftree}
for every finite $\Gamma \cup \{ \gamma_1^1, \dots, \gamma_1^k, \dots, \gamma_{n+1}^1, \dots, \gamma_{n+1}^k, \psi \} \subseteq Fm(\vdash)$.
\end{Definition}

In the presence of the IL, the semantic counterpart of the $\textup{BTWL}_n$ can de described as follows:

\begin{Theorem}\label{Thm:BTWL-correspondence-AAL}
A protoalgebraic logic $\vdash$ with the inconsistency lemma has the bounded top width $n$ law iff for every algebra $\A$ and every $F \in \mathsf{Spec}_\vdash(\A)$, there are a positive integer $m \leq n$ and maximal elements $G_1, \dots, G_m$ of $\mathsf{Spec}_\vdash(\A)$ such that every $H \in \mathsf{Spec}_\vdash(\A)$ extending $F$ is contained in some $G_i$.
\end{Theorem}

\begin{proof}
Notice that $\vdash$ has the $\textup{BTWL}_n$ precisely when it validates the metarules in $\mathsf{R}_\vdash(\Phi_n)$ induced by the Sahlqvist quasiequation $\Phi_n$ corresponding to the axiom $\mathsf{btw}_n$, defined in Example \ref{Exa:BTWn}. Furthermore, $\Phi_n$ is compatible with $\vdash$, because $\vdash$ has the IL by assumption. Therefore, we can apply the Abstract Sahlqvist Theorem, obtaining that $\vdash$ has the $\textup{BTWL}_n$ iff $\mathsf{Spec}_\vdash(\A) \vDash \mathsf{tr}(\Phi_{n})$, for every algebra $\A$. In view of Example \ref{Exa:BTWL-correspondence-easy}, the latter amounts to the demand that 
for every algebra $\A$ and every $F, H_1, \dots, H_{n+1} \in \mathsf{Spec}_\vdash(\A)$ such that $F$ is contained in each $H_i$, there are $G_1, \dots, G_n \in \mathsf{Spec}_\vdash(\A)$ extending $F$ such that each $H_i$ is contained in at least one $G_j$. Therefore, it only remains to prove that this condition is equivalent to that in the right hand side of the statement. 

The fact that the condition in the statement implies the one above is clear. To prove the converse, consider an algebra $\A$ satisfying the condition above. Then let $M$ be the set of maximal proper deductive filters of $\vdash$ on $\A$. 

\begin{Claim}
Every element of $\mathsf{Spec}_\vdash(\A)$ is contained in some element of $M$.
\end{Claim}

\begin{proof}[Proof of the Claim]
Suppose, with a view to contradiction, that there exists some $F \in \mathsf{Spec}_\vdash(\A)$ that cannot be extended to an element of $M$. Then consider the subposet of $\Fi{\A}$ with universe
\[
Y \coloneqq \{ G \in \Fi{\A} : G \notin {\downarrow}M \text{ and }F \subseteq G \text{ and }G \ne A \}.
\]
Since $F$ is meet irreducible, it is different from $A$ and, therefore, it belongs to $Y$. Consequently, the poset $Y$ is nonempty and we can apply Zorn's Lemma to deduce that there exists a maximal chain $C$ in $Y$.

We will prove that the join of $C$ in $\Fi{\A}$ is $A$. Suppose, with a view to contradiction, that $\bigvee C \subsetneq A$. Observe that the maximality of $C$ guarantees that $F \in C$, whence $F \subseteq \bigvee C$.
Together with the assumption that $F \notin {\downarrow}M$, this implies that $\bigvee C \notin {\downarrow}M$. Since by assumption $\bigvee C \ne A$, there exists a proper $G \in \Fi{\A}$ such that $\bigvee C \subsetneq G$. As a consequence $G \notin C$, which, by the maximality of $C$, yields $G \notin Y$. Since $F \subseteq \bigvee C \subseteq G$ and $F \notin {\downarrow}M$, this means that $G = A$, a contradiction with the assumption that $G$ is proper. Hence, we conclude that $\bigvee C = A$.

Now, consider one of the finite sets $\sim_n\!\!(x_1, \dots, x_n)$ witnessing the IL for $\vdash$. Since the IL guarantees that the finite set
\[
\{ x_1, \dots, x_n \} \, \cup \sim_n\!\!(x_1, \dots, x_n)
\]
is inconsistent, we obtain that
\[
\fg{\{a_1, \dots, a_n\} \, \cup \sim_n\!\!(a_1, \dots, a_n)} = A,
\]
for every and $a_1, \dots, a_n \in A$. Therefore, the deductive filter $A$ is finitely generated.

By Proposition \ref{prop : finitarity implies algebraicity}, this implies that $A$ is a compact element of $\Fi{\A}$. As a consequence, from $\bigvee C = A$ it follows that there exists a finite $C' \subseteq C$ such that $\bigvee C' = A$. Since $F \in C$, we may assume that $C'$ contains $F$ and, therefore, is nonempty. As $C'$ is a finite nonempty chain, we have $\bigvee C' \in C'$, whence $A = \bigvee C' \in C' \subseteq C \subseteq Y$. But this contradicts the definition of $Y$, according to which $A \notin Y$. 
\end{proof}

Now, consider an element $F \in \mathsf{Spec}_\vdash(\A)$ and let
\[
M_{F} \coloneqq \{ G \in M : F \subseteq G \}.
\]
Clearly, $M_F$ is a set of maximal elements of $\mathsf{Spec}_\vdash(\A)$. Furthermore, in view of the Claim, every element of $\mathsf{Spec}_\vdash(\A)$ extending $F$ is contained in some element of $M_F$. Therefore, to conclude the proof, it suffices to show that $\vert M_F \vert \leq n$. Suppose, with a view to contradiction, that there are distinct $H_1, \dots, H_{n+1} \in M_F$. As $M_F \subseteq \mathsf{Spec}_\vdash(\A)$, we can apply the assumption obtaining that there are $G_1, \dots, G_n \in \mathsf{Spec}_\vdash(\A)$ such that each $H_i$ is contained into some $G_j$. Therefore, there are $m < k \leq n+1$ and $j \leq n$ such that $H_m, H_k \subseteq G_j$. Since $G_j$ is proper (because it belongs to $\mathsf{Spec}_\vdash(\A)$), the maximality of $H_m$ and $H_k$ implies that $H_m = G = H_j$. But this contradicts the assumption that $H_1, \dots, H_{n+1}$ are all different.
\end{proof}

It is easy to see that a logic $\vdash$ with the IL has the $\textup{BTWL}_1$ iff it has the \emph{weak excluded middle law} (WEML, for short) in the sense that it validates the metarule
\begin{prooftree}
	\AxiomC{$\Gamma, \sim\!(\varphi_1, \dots, \varphi_n)    \rhd \psi$} 
	\AxiomC{$\Gamma, \sim\sim\!(\varphi_1, \dots, \varphi_n)    \rhd \psi$}
	\BinaryInfC{$\Gamma \rhd \psi$}
	\end{prooftree}
for every finite $\Gamma \cup \{ \varphi_1, \dots, \varphi_n, \psi \} \subseteq Fm(\vdash)$. Bearing this in mind, from Theorem \ref{Thm:BTWL-correspondence-AAL} we deduce:

\begin{Corollary}[\protect{\cite[Thm.\ 6.3]{LMR22}}]
A protoalgebraic logic $\vdash$ with the inconsistency lemma has the weak excluded middle law iff for every algebra $\A$ and every $F \in \mathsf{Spec}_\vdash(\A)$, there exists the greatest element of $\mathsf{Spec}_\vdash(\A)$ extending $F$.
\end{Corollary}

\section{Sahlqvist theory for fragments of $\mathsf{IPC}$ with implication}

The Abstract Sahlqvist Theorem can be also employed to derive Sahlqvist theorems for concrete deductive systems. In this section, we will do this for fragments of $\mathsf{IPC}$ including the connective $\to$. To this end, it is convenient to recall some basic concepts. Let $\A$ be a subreduct of a Heyting algebra in a language $\mathcal{L}_\to$ containing $\to$. Then, the formula $x \to x$ induces a constant term function on $\A$, whose constant value will be denoted by $1$. Accordingly, a formula $\varphi$ of $\mathcal{L}_\to$ is \emph{valid} in $\A$, in symbols $\A \vDash \varphi$, when $\A$ satisfies the equation $\varphi \thickapprox 1$. Furthermore, a subset $F$ of $A$ is said to be an \emph{implicative filter} of $\A$ if it contains $1$ and, for every $a, b \in A$,
\[
\text{if }\, a, a\to b \in F, \,\text{ then } \, b \in F.
\]
When ordered under the inclusion relation, the set of implicative filters of $\A$ forms a lattice. We denote its subposet of meet irreducible elements by $\A_\ast$.
\begin{Remark}
When the language of $\A$ contains $\land$, Condition (\ref{Eq:res-law}) guarantees that the implicative and semilattice filters of $\A$ coincide. Therefore, there is not clash in our usage of the notation $\A_\ast$ both for the posets of meet irreducible semilattice and implicative filters.
\qed
\end{Remark}

The importance of implicative filters is made apparent by the following observation.

\begin{Proposition}\label{Porp:sqdlkjak}
Let $\mathsf{L}$ be a fragment of $\mathsf{IPC}$ containing $\to$. Then, for every subreduct $\A$ in the language of $\mathsf{L}$ of a Heyting algebra, the deductive filters of $\mathsf{L}$ on $\A$ coincide with the implicative filters of $\A$.
\end{Proposition}

\begin{proof}
 Since $\mathsf{L}$ is an implicative logic in the sense of \cite{Ra74}, the result follows from \cite[Prop.\ 2.28]{AAL-AIT-f}.
\end{proof}

Given a finite set $\Gamma \cup \{ \varphi\}$ of formulas of $\mathsf{IPC}$, with $\Gamma = \{ \gamma_1, \dots, \gamma_n \}$, we write
\[
\Gamma \to \varphi \, \, \text{ as a shorthand for the singleton } \, \, \{\gamma_1 \to (\gamma_2 \to (\dots (\gamma_{n} \to \varphi)\dots ))\}.
\]
Recall from Remark \ref{Rem:IPC-all-is-compatible} that every formula $\varphi$ of $\mathsf{IPC}$ is compatible with $\mathsf{IPC}$. Accordingly, given $k \in \mathbb{Z}^+$, we denote by $\boldsymbol{\varphi}^k$ the finite set of formulas of $\mathsf{IPC}$ associated with $\varphi$. Moreover, with every Sahlqvist quasiequation
\[
\Phi = \varphi_1 \land y \leq z \, \& \dots \& \, \varphi_m \land y \leq z \Longrightarrow y \leq z,
\]
we associate the set of formulas
\[
\mathsf{A}(\Phi) \coloneqq \bigcup \{ ((\boldsymbol{\varphi_1}^k \to y) \cup \dots \cup (\boldsymbol{\varphi_m}^k \to y)) \to y : k \in \mathbb{Z}^+ \},
\]
where $y$ is a variable that does not occur in $\boldsymbol{\varphi_1}^k, \dots, \boldsymbol{\varphi_m}^k$.

\begin{exa}\label{Exa:GD-axiom}
Consider the Sahlqvist quasiequation
\[
\Phi = (x_1 \to x_2) \land y \leq z \, \& \, (x_2 \to x_1) \land y \leq z \Longrightarrow y \leq z
\]
corresponding to the G\"odel-Dummett axiom (see Example \ref{Exa:BTWn}). In view of Example \ref{Exa:IPC-implications-translated}, for every $k \in \mathbb{Z}^+$, the sets $(\boldsymbol{x_1 \to x_2})^k \to y$ and $(\boldsymbol{x_2 \to x_1})^k \to y$ are the singletons containing, respectively, the formulas
\begin{align*}
\psi_1^k \coloneqq (x_1^1 \to ( \dots ( x_1^k \to x_2^1)\dots)) &\to ( (x_1^1 \to ( \dots ( x_1^k \to x_2^2)\dots))\\
& \to
(\dots ((x_2^1 \to ( \dots ( x_2^k \to x_1^k)\dots))\to y) \dots ))
\end{align*}
and 
\begin{align*}
\psi_2^k \coloneqq (x_2^1 \to ( \dots ( x_2^k \to x_1^1)\dots)) &\to ( (x_2^1 \to ( \dots ( x_2^k \to x_1^2)\dots))\\
& \to
(\dots ((x_2^1 \to ( \dots ( x_2^k \to x_1^k)\dots))\to y) \dots )).
\end{align*}
Consequently, $\mathsf{A}(\Phi)$ is the set $\{ \psi_1^k \to (\psi_2^k \to y) : k \in \mathbb{Z}^+\}$.
\qed
\end{exa}

Given an algebra $\A$, we denote by $\VVV(\A)$ the variety generated by $\A$. We rely on the following observation.

\begin{Lemma}\label{Lem:canonicity-arrow-IPC-trick}
Let $\Phi$ be a Sahlqvist quasiequation in a sublanguage $\mathcal{L}_\to$ of $\mathcal{L}$ containing $\to$. For every $\mathcal{L}_\to$-subreduct $\A$ of a Heyting algebra,
\[
\A \vDash \mathsf{A}(\Phi) \, \,  \Longleftrightarrow \, \, \B_\ast \vDash \mathsf{tr}(\Phi) \text{, for every } \B \in \VVV(\A).
\]
\end{Lemma}

\begin{proof}
Throughout the proof we will assume that
\[
\Phi = \varphi_1 \land y \leq z \, \& \dots \& \, \varphi_m \land y \leq z \Longrightarrow y \leq z.
\]
Furthermore, let $\mathsf{L}$ be the $\mathcal{L}_\to$-fragment of $\mathsf{IPC}$ and let $\mathsf{L}(\A)$ be the extension of $\mathsf{L}$ axiomatized, relatively to $\mathsf{L}$, by the formulas valid in $\A$. It is well known that $\mathsf{L}(\A)$ is algebraized by $\VVV(\A)$. 

We begin by proving the following equivalences:
\begin{align*}
\A \vDash \mathsf{A}(\Phi) \, \,& \Longleftrightarrow \, \,\emptyset \vdash_{\mathsf{L}(\A)} \mathsf{A}(\Phi)\\
\, \,& \Longleftrightarrow \, \,\mathsf{L}(\A) \text{ validates the metarules in }\mathsf{R}_{\mathsf{L}(\A)}(\Phi)\\
\, \,& \Longleftrightarrow \, \,\B_\ast \vDash \mathsf{tr}(\Phi) \text{, for every } \B \in \VVV(\A).
\end{align*}
The first equivalence follows from the definition of $\mathsf{L}$. For the second, observe that, being a fragment of $\mathsf{IPC}$ with $\to$, the logic $\mathsf{L}$ inherits the DT of $\mathsf{IPC}$. Since the DT persists in axiomatic extensions, the DT of $\mathsf{IPC}$ holds also in $\mathsf{L}(\A)$. Consequently, the second equivalence follows from the part of Remark \ref{Rem:DT-Conjunction-simplifications} devoted to the DT. To prove the third one, we begin by showing that $\Phi$ is compatible with $\mathsf{L}(\A)$. Suppose, for instance, that the connective $\lor$ occurs in some $\varphi_i$. Then $\mathcal{L}_\to$ contains $\lor$. Since the PC persists in axiomatic extensions of fragments of $\mathsf{IPC}$ with $\lor$, we conclude that $\mathsf{L}(\A)$ has the PC as desired. A similar argument applies to the cases where $0, \lnot$, or $\to$ occur in some $\varphi_i$, thereby yielding that $\Phi$ is compatible with $\mathsf{L}(\A)$. Furthermore, recall that $\VVV(\A)$ algebraizes $\mathsf{L}(\A)$. Lastly, Proposition \ref{Porp:sqdlkjak} guarantees that $\mathsf{Spec}_{\mathsf{L}(\A)}(\B) = \B_\ast$, for every $\B \in \VVV(\A)$. Therefore, we can apply Corollary \ref{cor : correspondence - algebraizable}, thus establishing the third equivalence.
\end{proof}

Notably, one can view $\mathsf{A}(\Phi)$ as the equational version of the quasiequation $\Phi$, as made precise by the next observation.

\begin{Proposition}\label{Prop:last-protosition-IPC-arrow}
A Heyting algebra $\A$ validates a Sahlqvist quasiequation $\Phi$ iff it validates the formulas in $\mathsf{A}(\Phi)$.
\end{Proposition}

\begin{proof}
In view of Lemma \ref{Lem:canonicity-arrow-IPC-trick}, it suffices to establish the following equivalences:
\begin{align*}
 \B_\ast \vDash \mathsf{tr}(\Phi) \text{, for every } \B \in \VVV(\A)\, \,& \Longleftrightarrow \, \,\B \vDash \Phi \text{, for every } \B \in \VVV(\A)\\
\, \,& \Longleftrightarrow \, \, \A \vDash \Phi.
\end{align*}
To this end, we will assume that
\[
\Phi = \varphi_1 \land y \leq z \, \& \dots \& \, \varphi_m \land y \leq z \Longrightarrow y \leq z.
\]

The first equivalence holds by Corollary \ref{corollary : corollary of enhanced sahlqvist}. To prove the nontrivial part of the second, suppose that $\A \vDash \Phi$. In view of Corollary \ref{Cor:meaning-Sahlqvist-IPC}, the equation $\varphi_1 \lor \dots \lor \varphi_n \thickapprox 1$ is valid in $\A$. Therefore, it is also valid $\VVV(\A)$. With another application of Corollary \ref{Cor:meaning-Sahlqvist-IPC}, we conclude that $\Phi$ is valid in all the members of $\VVV(\A)$ as desired.
\end{proof}

As in the case of Theorem \ref{thm : enhanced intuitionistic sahlqvist}, Sahlqvist Theorem from fragments of $\mathsf{IPC}$ with $\to$ takes the form of a canonicity result.

\begin{Theorem}\label{Thm:canonicity-arrow-IPC}
Let $\Phi$ be a Sahlqvist quasiequation in a sublanguage $\mathcal{L}_\to$ of $\mathcal{L}$ containing $\to$. If an $\mathcal{L}_\to$-subreduct $\A$ of a Heyting algebra validates $\mathsf{A}(\Phi)$, then also $\mathsf{Up}(\A_\ast)$ validates $\mathsf{A}(\Phi)$.
\end{Theorem}

\begin{proof}
Suppose that $\A \vDash \mathsf{A}(\Phi)$. In view of Lemma \ref{Lem:canonicity-arrow-IPC-trick}, this yields $\A_\ast \vDash \mathsf{tr}(\Phi)$. Therefore, we can apply the correspondence part of the Intuitionistic Sahlqvist Theorem, obtaining $\mathsf{Up}(\A_\ast) \vDash \Phi$. Lastly, by Proposition \ref{Prop:last-protosition-IPC-arrow}, this amounts to $\mathsf{Up}(\A_\ast) \vDash \mathsf{A}(\Phi)$.
\end{proof}

Bearing in mind that $\Phi$ and $\mathsf{A}(\Phi)$ axiomatize the same class of Heyting algebras (Proposition \ref{Prop:last-protosition-IPC-arrow}), a straightforward adaptation of the proof of Corollary \ref{corollary : corollary of enhanced sahlqvist} yields the following:

\begin{Corollary}\label{Cor:last-corollary-end}
Let $\Phi$ be a Sahlqvist quasiequation in a sublanguage $\mathcal{L}_\to$ of $\mathcal{L}$ containing $\to$. For every $\mathcal{L}_\to$-subreduct $\A$ of a Heyting algebra, it holds that $\boldsymbol{A} \vDash \mathsf{A}(\Phi) \text{ iff } \boldsymbol{A}_\ast \vDash \mathsf{tr}(\Phi)$.
\end{Corollary}

\begin{exa}
The $\langle \to \rangle$-subreduct of Heyting algebras are called \emph{Hilbert algebras} \cite{Di65,Di66}. Let $\Phi$ be the Sahlqvist quasiequation corresponding to the G\"odel-Dummett axiom. Since a Hilbert algebra validates $\mathsf{A}(\Phi)$ iff it validates the single formula $\varphi = ((x \to y) \to z) \to (  ((y \to z) \to z) \to z$, in view of Corollary \ref{Cor:last-corollary-end} and Example \ref{Exa:BTWL-correspondence-easy}, we obtain that
\[
\A \vDash \varphi \, \, \Longleftrightarrow \, \, \A_\ast \text{ is a root system},
\]
for every Hilbert algebra \cite[Thm.\ 4.5]{MontHil96}.
\qed
\end{exa}

\section{A correspondence theorem for intuitionistic linear logic}

We close this paper by deriving a correspondence theorem for intuitionistic linear logic \cite{Gi87} from the Abstract Sahlqvist Theorem (cf.\ \cite{Suzukia,Suzuki}). To this end, recall that a \textit{commutative FL-algebra} is a structure $\A = \langle A; \land, \lor, \cdot, \to, 0, 1 \rangle$ comprising a commutative monoid $\langle A; \cdot , 1 \rangle$ and a lattice $\langle A; \land, \lor \rangle$ such that for every $a, b, c \in A$,
\begin{equation}\label{Eq:FL-residuation-law}
a \cdot b \leq c \, \,\Longleftrightarrow\, \, a \leq b \to  c.
\end{equation}
The class of commutative FL-algebras forms a variety which we denote by $\mathsf{FL}_{\mathsf e}$ \cite{GaJiKoOn07}. \textit{Intuitionistic linear logic} $\mathsf{ILL}$ is the logic formulated in the language of commutative FL-algebras defined, for every set $\Gamma \cup \{ \varphi \}$ of formulas, as follows:
\[
\Gamma \vdash_{\mathsf{ILL}} \varphi \, \, \text{ iff } \, \, \text{ there exists a finite }\Sigma \subseteq \Gamma \text{ such that }\mathsf{FL}_{\mathsf{e}} \vDash \foo_{\gamma \in \Sigma} \gamma \geq 1 \Longrightarrow \varphi \geq 1.
\]
It is well known that every axiomatic extension $\vdash$ of $\mathsf{ILL}$ is algebraized by the variety $\mathsf{K}_\vdash$ of commutative FL-algebras axiomatized by the set of equations $\{ \varphi \geq 1 : \emptyset \vdash \varphi\}$, as witnessed by the sets $\tau \coloneqq \{ x \geq 1 \}$ and $\Delta(x, y) \coloneqq \{ x \to y, y \to x \}$ \cite[Thm.\ 3.3]{GaOn06}. In particular, $\mathsf{ILL}$ is algebraized by $\mathsf{FL}_{\mathsf{e}}$.

Given an algebra $\A$ in the language of $\mathsf{ILL}$, an element $a \in A$, and $n \in \mathbb{Z}^+$, we define an element $a^n$ of $A$ by setting
\[
a^1 \coloneqq a \, \, \text{ and } \, \, a^{m+1} \coloneqq a^m \cdot a, \text{ for every }m \geq 1.
\]
We will rely on the following property of $\mathsf{ILL}$:

\begin{Proposition}[\protect{\cite[Thm.\ 4.9]{GaOn06}}]\label{Prop:LDDT-FLe}
For every algebra $\A$ and $X \cup \{ a, b \} \subseteq A$,
\[
a \in \fgL{X \cup \{ b \}} \, \, \text{ iff } \, \, (1 \land b)^n \to a \in \fgL{X}\text{, for some }n \in \mathbb{Z}^+.
\]
\end{Proposition}

When $\A$ is the algebra of formulas $\boldsymbol{Fm}(\mathsf{ILL})$, this specializes to the following:

\begin{Corollary}\label{Cor:LDDT-FLe}
For every set $\Gamma \cup \{ \psi, \varphi \}$ of formulas of $\mathsf{ILL}$, we have
\[
\Gamma, \psi \vdash_{\mathsf{ILL}} \varphi \, \, \text{ iff } \, \, \Gamma \vdash_{\mathsf{ILL}} (1 \land \psi)^n \to \varphi\text{, for some }n \in \mathbb{Z}^+.
\]
\end{Corollary}

In order to obtain a correspondence theorem for $\mathsf{ILL}$, it is convenient to identify the axiomatic extensions of $\mathsf{ILL}$ with the IL, the DT, and the PC. For the DT and the PC we have the following:

\begin{Proposition}[\protect{\cite[Prop.\ 3.15]{Ga-PhD}}]
An axiomatic extension of $\mathsf{ILL}$ has the deduction theorem iff there exists some $k \in \mathbb{Z}^+$ such that the theorems of $\vdash$ include the formula $(1 \land x)^k \to (1 \land x)^{k+1}$. In this case, the DT is witnessed by the sets of the form
\[ 
(x_1, \dots, x_n)\!\Rightarrow_{nm} \!\!(y_1, \dots, y_m)= \{ (1 \land x_1 \land \dots \land x_n)^k \to (y_1 \land \dots \land y_m) \}.
\]
\end{Proposition}

\begin{Proposition}\label{Prop:PC-for-FLe}
Every axiomatic extension $\vdash$ of $\mathsf{ILL}$ has the proof by cases, as witnessed by the sets of the form
\[ 
(x_1, \dots, x_n)\!\bigcurlyvee_{nm} (y_1, \dots, y_m)= \{ (1 \land x_1 \land \dots \land x_n) \lor (1 \land y_1 \land \dots \land y_m) \}.
\]
\end{Proposition}

\begin{proof}
This is essentially \cite[Example 4.10.4]{CN20xx-book}, where the result is stated for the natural expansion $\mathsf{SL}_{\mathsf{aE}}$ of $\mathsf{ILL}$ with bounds. 
\end{proof}

In order to address the case of the IL, it is convenient to introduce the following shorthand for every formula $\varphi$ of $\mathsf{ILL}$:
\[
\bot \coloneqq 1 \land (1 \to 0) \land (0 \to 1) \land (1 \to (1 \to 1)) \land ((1 \to 1) \to 1) \, \, \text{ and } \, \, \lnot \varphi \coloneqq \varphi \to \bot.
\]
\begin{Proposition}\label{Prop:FLe-IL-description}
An axiomatic extension $\vdash$ of $\mathsf{ILL}$ has the inconsistency lemma iff there exist some $k \in \mathbb{Z}^+$ and a function $f \colon \mathbb{Z}^+ \to \mathbb{Z}^+$ such that the theorems of $\vdash$ include the formulas
\[
\bot^k \to x \, \, \text{ and } \, \, (1 \land \lnot (x \land 1)^m )^{f(m)} \to \lnot (1 \land x)^k,
\]
for every $m \in \mathbb{Z}^+$. In this case, the IL is witnessed by the sets of the form
\[
\sim_n\!\!(x_1, \dots, x_n) \coloneqq \{ \lnot (1 \land x_1 \land \dots \land x_n)^k \}.
\]
\end{Proposition}

	\begin{proof}
We will often use the fact that $\vdash$ is algebraized by $\mathsf{K}_\vdash$, as witnessed by the sets $\tau \coloneqq \{ x \geq 1 \}$ and $\Delta(x, y) \coloneqq \{ x \to y, y \to x \}$. Similarly, we will repeatedly appeal to the fact that for every $\A \in \mathsf{K}_\vdash$ and $a, b \in A$,
\[
a \leq b \, \, \Longleftrightarrow \, \, 1 \leq a \to b
\]
and
\[
\text{if } a \leq 1 \text{, then } a^{n+1} \leq a^n\text{, for every }n \in \mathbb{Z}^+.
\]
The first property follows from Condition (\ref{Eq:FL-residuation-law}) and the assumption that $\langle A; \cdot, 1 \rangle$ is a monoid. The second holds because the operation $\cdot$ is order preserving in both coordinates and, therefore, the assumption that $a \leq 1$ guarantees that $a^{n+1} = a^n \cdot a \leq a^n \cdot 1 = a^n$. These facts will be used in the proof without further notice.

We begin by proving the implication from left to right in the statement. Accordingly, suppose that $\vdash$ has the IL. We rely on the next observations:
		
\begin{Claim}\label{Claim : bot is inconsistent}
The formula $\bot$ is inconsistent in $\vdash$.
\end{Claim}		
	\begin{proof}[Proof of the Claim]
		First, we will prove that
				\begin{equation}\label{Eq:FL:IL_characterization-trick-1}
		\bot \vdash \varphi \to \psi\text{, for every pair $\varphi$ and $\psi$ of formulas in which no variable occurs.}
				\end{equation}
			 Accordingly, consider two such formulas $\varphi$ and $\psi$. It suffices to show that
		\[
\mathsf{K}_\vdash \vDash		\bot \geq 1 \Longrightarrow \varphi \to \psi \geq 1.	
\]
To this end, consider an algebra $\A \in \mathsf{K}_\vdash$ such that $1^\A \leq \bot^\A$. By the definition of $\bot$, this yields $0^\A = 1^\A = 1^\A \to 1^\A$. As
\[
\mathsf{FL}_{\mathsf{e}}\vDash 1 \thickapprox1 \land 1 \thickapprox 1 \lor 1 \thickapprox 1 \cdot 1,
\]
this implies that $\{ 1^\A \}$ is the universe of a subalgebra of $\A$. Consequently, $\varphi^\A, \psi^\A \in \{ 1^\A \}$, because $\varphi$ and $\psi$ have no variables. Therefore, $\varphi^\A = \psi^\A$. In particular, $\varphi^\A \leq \psi^\A$, which amounts to $1^\A \leq \varphi^\A \to \psi^\A$, as desired. This establishes Condition (\ref{Eq:FL:IL_characterization-trick-1}).

Now, recall that the IL guarantees that the set $\{1\} \?\cup\! \sim_1 \!\!(1)$ is inconsistent.\ In particular, $1 \land \bigwedge\! \sim_1\!\!(1) \vdash y$. Furthermore, from Condition (\ref{Eq:FL:IL_characterization-trick-1}) it follows that $\bot \vdash 1 \to (1 \land \bigwedge \!\sim_1\!\!(1))$. Since $1$ is a theorem of $\vdash$ and $x, x\to y \vdash y$, this yields that $\bot \vdash 1 \land \bigwedge \!\sim_1\!\!(1)$. Together with $1 \land \bigwedge\! \sim_1\!\!(1) \vdash y$, this implies that $\bot \vdash y$. By substitution invariance, we conclude that $\bot$ is inconsistent.
	\end{proof}	
	
	\begin{Claim}\label{cor : IL implies something}
		For every $\boldsymbol{A} \in \mathsf{K_\vdash}$ and $a \in A$, we have that 
		$
		1 \leq \bigwedge\!\thicksim_1\!\!(a)
		$
		iff there exists some $n \in \mathbb{Z}^+$ such that 
		$(1 \land a)^n \leq \bot$.
	\end{Claim}
	\begin{proof}[Proof of the Claim]
		Recall from Theorem \ref{thm : il} that the semilattice $\Fic{\A}$ is pseudocomplemented. Therefore, $\fg{a} = A$ iff the pseudocomplement $\fg{\thicksim_1\!\!(a)}$ of $\fg{a}$ in $\Fic{\A}$ is included in $\fg{\emptyset}$, in symbols,
		\[
		\fg{a} = A \, \, \Longleftrightarrow \, \, \thicksim_1\!\!(a) \subseteq \fg{\emptyset}.
		\]
		 On the other hand, by applying in succession Claim \ref{Claim : bot is inconsistent} and Proposition \ref{Prop:LDDT-FLe}, we obtain
 \[
 \fg{a} = A \, \, \Longleftrightarrow \, \, \bot \in \fg{a} \, \, \Longleftrightarrow \, \, \lnot (1 \land a)^n \in \fg{\emptyset}\text{, for some }n \in \mathbb{Z}^+.
 \] 
 From the two displays above and the fact that $\fg{\emptyset} = \{ c \in A : c \geq 1 \}$ it follows that
	\[
	1 \leq \bigwedge \!\sim_1 \!\! (a) \, \, \Longleftrightarrow \, \, 1 \leq \lnot (1 \land a)^n \text{, for some }n \in \mathbb{Z}^+,
	\]	 
 where the condition on the right hand side is equivalent to the demand that there exists some $n \in \mathbb{Z}^n$ such that $(1 \land a)^n \leq \bot$.
	\end{proof}

In virtue of Claim \ref{Claim : bot is inconsistent} and Corollary \ref{Cor:LDDT-FLe}, there exists some $t \in \mathbb{Z}^+$ such that $\emptyset \vdash (1 \land \bot)^{t} \to x$. Therefore, $\mathsf{K}_\vdash \vDash  (1 \land \bot)^t \to x \geq 1$, which amounts to $\mathsf{K}_\vdash \vDash (1 \land \bot)^t \leq x$. By the definition of $\bot$, we have $\mathsf{K}_\vdash \vDash \bot \leq 1$. As a consequence, $\mathsf{K}_\vdash \vDash \bot^s \leq x$, for every positive integer $s \geq t$. This, in turn, yields that
\begin{equation}\label{Eq:FL-decreasing-bot}
\emptyset \vdash \bot^{s} \to x \text{, for every positive intenger }s \geq t.
\end{equation}

On the other hand, in view of Claim \ref{cor : IL implies something}, we have that
		\[
		\mathsf{Th}(\mathsf{K_\vdash}) \cup \{(1 \land x)^n  \nleq \bot \colon n \in \mathbb{Z}^+\} \Vdash_{\mathsf{FOL}} 1  \nleq \bigwedge \! \thicksim_1\!\!(x),
		\]
		where $\mathsf{Th}(\mathsf{K_\vdash})$ is the elementary theory of $\mathsf{K}_\vdash$ and $\Vdash_{\mathsf{FOL}}$ is the deducibility relation of first order logic. 
				By the Compactness Theorem of first order logic, the previous display implies that there are some positive integers $n_1,\dots,n_m \in \mathbb{Z}^+$ such that 
		\[
		\mathsf{Th}(\mathsf{K_\vdash}) \cup \{ (1 \land x)^{n_1} \nleq \bot, \dots, (1 \land x)^{n_m} \nleq \bot \}  \Vdash_{\mathsf{FOL}} 1  \nleq \bigwedge \! \thicksim_1\!\!(x).
		\] 
Letting $k \coloneqq \max\{n_1,\dots,n_m, t\}$, the above display implies
		\[
		\mathsf{Th}(\mathsf{K_\vdash}) \cup \{ (1 \land x)^{k} \nleq \bot \}  \Vdash_{\mathsf{FOL}} 1  \nleq \bigwedge \! \thicksim_1\!\!(x).
		\]
		This, in turn, amounts to the following:
		\[
		\mathsf{K}_\vdash \vDash \bigwedge \! \thicksim_1\!\!(x) \geq 1 \Longrightarrow (1 \land x)^{k} \leq \bot.
		\]
		In view of Claim \ref{cor : IL implies something}, this yields that for every $m \in \mathbb{Z}^+$,
		\[
		\mathsf{K}_\vdash \vDash (1 \land x )^m \leq \bot \Longrightarrow  (1 \land x)^{k}  \leq \bot,
		\]
where the condition above can be equivalently phrased as 
		\[
	\mathsf{K}_\vdash \vDash \lnot (1 \land x )^m \geq 1 \Longrightarrow  \lnot (1 \land x)^{k} \geq 1.
		\]
Consequently, $\lnot (1 \land x )^m \vdash \lnot (1 \land x)^{k}$, for every $m \in \mathbb{Z}^+$. In view of Proposition \ref{Prop:LDDT-FLe}, for every $m \in \mathbb{Z}$ there exists some $f(m) \in \mathbb{Z}^+$ such that
\[
\emptyset \vdash (1 \land \lnot (1 \land x )^m)^{f(m)} \to \lnot (1 \land x)^{k}.
\]

Lastly, since the definition of $k$ guarantees that $k \geq t$, from Condition (\ref{Eq:FL-decreasing-bot}) it follows that $\emptyset \vdash \bot^k \to x$.

Then we turn to prove the implication from right to left in the statement. We will show that the sets of the form
\[
\thicksim_n\!\!(x_1,\dots,x_n) \coloneqq \{ \lnot(1 \land x_1 \land  \dots \land x_n)^k \}
\]
witness the IL for $\vdash$, i.e., that for every finite $\Gamma \cup \{\varphi_1,\dots,\varphi_n\} \sub Fm(\vdash)$, 
\[
				\Gamma \cup \{\varphi_1,\dots,\varphi_n\} \text{ is inconsistent } \, \, \text{ iff } \, \,  \Gamma \vdash \lnot(1 \land \varphi_1 \land  \dots \land \varphi_n)^k.
\]		
		Suppose first that $\Gamma \cup \{\varphi_1,\dots,\varphi_n\}$ is inconsistent. Then $\Gamma,   \varphi_1 \land \dots \land \varphi_n \vdash \bot$. In view of Corollary \ref{Cor:LDDT-FLe}, there exists some $m \in \mathbb{Z}^+$ such that $\Gamma \vdash \lnot(1 \land \varphi_1 \land \dots \land \varphi_n)^m$. As $x \vdash (1 \land x)^{f(m)}$, this yields $\Gamma \vdash (1 \land \lnot(1 \land \varphi_1 \land \dots \land \varphi_n)^m)^{f(m)}$.		
		Since by assumption $\emptyset \vdash (1 \land \lnot (1 \land x )^m)^{f(m)} \to \lnot (1 \land x)^{k}
$, with an application of modus ponens, we obtain that $\Gamma \vdash \lnot(1 \land \varphi_1 \land \dots \land \varphi_n)^k$, as desired.
		
To prove the converse, suppose that 	$\Gamma \vdash \lnot(1 \land \varphi_1 \land  \dots \land \varphi_n)^k$. By Corollary \ref{Cor:LDDT-FLe}, this implies that $\Gamma \cup \{ \varphi_1, \dots, \varphi_n \} \vdash \bot$. Furthermore, as $x \vdash x^k$, we get 	$\Gamma \cup \{ \varphi_1, \dots, \varphi_n \} \vdash \bot^k$. Since the set $\Gamma \cup \{ \varphi_1, \dots, \varphi_n \}$ is finite, there exists a variable $x$ that does not occur in any of its members. As by assumption $\emptyset \vdash \bot^k \to x$, by modus ponens we obtain $\Gamma \cup \{ \varphi_1, \dots, \varphi_n \} \vdash x$. Since, $x$ does not occur in the formulas of $\Gamma \cup \{ \varphi_1, \dots, \varphi_n \}$, by substitution invariance, we obtain that $\Gamma \cup \{ \varphi_1, \dots, \varphi_n \} \vdash \psi$ for every formula $\psi$. Hence, we conclude that $\Gamma \cup \{ \varphi_1, \dots, \varphi_n \}$ is inconsistent.
\end{proof}

Let $\vdash$ be an axiomatic extension of $\mathsf{ILL}$. Notice that, if $\varphi$ is a formula of $\mathsf{IPC}$ compatible with $\vdash$, then the finite set of formulas $\boldsymbol{\varphi}^k$ is interderivable in $\vdash$ with the conjunction $\bigwedge \boldsymbol{\varphi}^k$, for every $k \in \mathbb{Z}^+$. Accordingly, from now on we will assume that the expressions of the form $\boldsymbol{\varphi}^k$ stand for formulas $\bigwedge \boldsymbol{\varphi}^k$ of $\vdash$, as opposed to a sets of formulas of $\vdash$. 

Furthermore, recall that $\mathsf{K}_\vdash$ is a variety, whence it is closed under $\HHH$. Consequently, $\mathsf{Spec}_{\mathsf{K}_\vdash}(\A)$ coincides with the poset of meet irreducible congruences of $\A$, for every $\A \in \mathsf{K}_\vdash$. Because of this, when $\A \in \mathsf{K}_\vdash$, we will write  $\mathsf{Spec}(\A)$ as a shorthand for $\mathsf{Spec}_{\mathsf{K}_\vdash}(\A)$. Bearing this in mind, we obtain the desired correspondence theorem:

\begin{Theorem}\label{Thm:linear-main}
Let $\Phi = \varphi_1 \land y \leq z \, \& \dots \& \, \varphi_m \land y \leq z \Longrightarrow y \leq z$ be a Sahlqvist quasiequation compatible with an axiomatic extension $\vdash$ of $\mathsf{ILL}$. Then the theorems of $\vdash$ include the formula $(1 \land \boldsymbol{\varphi_1}^1) \lor \dots \lor (1 \land \boldsymbol{\varphi_m}^1)$ iff $\mathsf{Spec}(\A) \vDash \mathsf{tr}(\Phi)$, for every algebra $\A \in \mathsf{K}_\vdash$.
\end{Theorem}

\begin{proof}
Observe that $\land$ is a conjunction for $\vdash$ and that $\vdash$ has the PC, as witnessed by sets in Proposition \ref{Prop:PC-for-FLe}. Therefore, from Remark \ref{Rem:DT-Conjunction-simplifications} it follows that the theorems of $\vdash$ include the formula $(1 \land \boldsymbol{\varphi_1}^1) \lor \dots \lor (1 \land \boldsymbol{\varphi_m}^1)$ iff the logic $\vdash$ validates the metarules in $\mathsf{R}_\vdash(\Phi)$. But, since $\vdash$ is algebraized by $\mathsf{K}_\vdash$, we can apply Corollary \ref{cor : correspondence - algebraizable} obtaining that the latter condition is equivalent to the demand that $\mathsf{Spec}(\A) \vDash \mathsf{tr}(\Phi)$, for every $\A \in \mathsf{K}_\vdash$.
\end{proof}

\begin{exa}
Let $\vdash$ be an axiomatic extension of $\mathsf{ILL}$ with the IL. Then there are some $k\in \mathbb{Z}^+$ and a function $f \colon \mathbb{Z}^+ \to \mathbb{Z}^+$ witnessing the property in the statement of Proposition \ref{Prop:FLe-IL-description}. We will prove that the following conditions are equivalent for every $n \in \mathbb{Z}^+$:
\benroman
\item\label{item:FLe-IL-1} The logic $\vdash$ has the $\textup{BTWL}_n$;
\item\label{item:FLe-IL-2} The theorems of $\vdash$ include the formula
\[
\bigvee_{1 \leq i \leq n+1}\Big( 1 \land \lnot ( 1 \land \bigwedge_{1 \leq j < i} x_j \land  \lnot(1 \land x_i)^k)^k  \Big);
\]
\item\label{item:FLe-IL-3} For every $\A \in \mathsf{K}_\vdash$ and $\theta \in \mathsf{Spec}(\A)$, there are a positive integer $m \leq n$ and maximal elements $\phi_1, \dots, \phi_m$ of $\mathsf{Spec}(\A)$ such that every $\eta \in \mathsf{Spec}(\A)$ extending $\theta$ is contained in some $\phi_i$.
\eroman

First, recall from Example \ref{Exa:BTWn} that the Sahlqvist quasiequation
\[
\Phi_n = \varphi_1 \land y \leq z \, \& \dots \& \, \varphi_{n+1} \land y \leq z \Longrightarrow y \leq z,
\]
 corresponding to the $\mathsf{btw}_n$ axiom is defined setting, for every $i \leq n+1$,
\[
\varphi_i \coloneqq \lnot(\lnot x_i \land \bigwedge_{0 < j< i}x_j).
\]

To prove the equivalence between Conditions (\ref{item:FLe-IL-1}) and (\ref{item:FLe-IL-2}), recall that the logic $\vdash$ has the $\textup{BTWL}_n$ precisely when it validates the metarules in $\mathsf{R}_\vdash(\Phi_n)$. As observed in the proof of Theorem \ref{Thm:linear-main}, this happens iff the theorems of $\vdash$ include the formula $(1 \land \boldsymbol{\varphi_1}^1) \lor \dots \lor (1 \land \boldsymbol{\varphi_{n+1}}^1)$. But the latter coincides with the formula in Condition (\ref{item:FLe-IL-2}), because the IL for $\vdash$ is witnessed by the sets of formulas in Proposition \ref{Prop:FLe-IL-description}.

Lastly, recall that $\mathsf{Spec}(\A) \cong \mathsf{Spec}_\vdash(\A)$ for every $\A \in \mathsf{K}_\vdash$, because $\vdash$ is algebraized by $\mathsf{K}_\vdash$. Therefore, the implication (\ref{item:FLe-IL-1})$\Rightarrow$(\ref{item:FLe-IL-3}) follows from Theorem \ref{Thm:BTWL-correspondence-AAL}. To prove the converse, suppose that Condition (\ref{item:FLe-IL-3}) holds. Then it is easy to check that $\mathsf{Spec}(\A) \vDash \mathsf{tr}(\Phi_n)$ for every $\A \in \mathsf{K}_\vdash$, where $\mathsf{tr}(\Phi_n)$ is the first order sentence mentioned in Example \ref{Exa:BTWL-correspondence-easy}. Consequently, we can apply Corollary \ref{cor : correspondence - algebraizable}, obtaining that $\vdash$ validates the metarules in $\mathsf{R}_\vdash(\Phi_n)$, which means that $\vdash$ has the $\textup{BTWL}_n$.
\qed
\end{exa}

\begin{exa}
The following conditions are equivalent for an axiomatic extension $\vdash$ of $\mathsf{ILL}$:\footnote{For a description of semisimple axiomatic extensions of an expansion of $\mathsf{ILL}$ with bounds, see \cite[Thm.\ 3.45]{PrenLav20}.}
\benroman
\item\label{item:FLe-IL-1b} The logic $\vdash$ has the EML;
\item\label{item:FLe-IL-2b} There exist some $k \in \mathbb{Z}^+$ and a function $f \colon \mathbb{Z}^+ \to \mathbb{Z}^+$ such that the theorems of $\vdash$ include the formulas
\[
\bot^k \to x, \, \,\text{ }  (1 \land \lnot (x \land 1)^m )^{f(m)} \to \lnot (1 \land x)^k, \, \, \text{ and } \, \, (1 \land x) \lor ( 1 \land \lnot (x \land 1)^k),
\]
for every $m \in \mathbb{Z}^+$;
\item\label{item:FLe-IL-3b} The logic $\vdash$ is semisimple and has the IL.
\eroman

In view of Remark \ref{Rem:EML-implies-IL}, the logic $\vdash$ has the EML iff it has the IL and validates the metarules in $\mathsf{R}_\vdash(\Phi)$, where
\[
\Phi = x \land y \leq z \, \&  \, \lnot x \land y \leq z \Longrightarrow y \leq z.
\]
As observed in the proof of Theorem \ref{Thm:linear-main}, the logic $\vdash$ validates the rules in $\mathsf{R}_\vdash(\Phi)$ precisely when its theorems contain $(1 \land \boldsymbol{x}^1) \lor (1 \land (\boldsymbol{\lnot x})^1)$. Consequently, Condition (\ref{item:FLe-IL-1b}) can be rephrased as the demand that $\vdash$ has the IL and its theorems include the formula $(1 \land \boldsymbol{x}^1) \lor (1 \land (\boldsymbol{\lnot x})^1)$.

On the other hand, in view of Proposition \ref{Prop:FLe-IL-description}, the demand that the theorems of $\vdash$ contain the first two formulas in Condition (\ref{item:FLe-IL-2b}) amounts to the assumption that $\vdash$ has the IL. The third formula in Condition (\ref{item:FLe-IL-2b}) is precisely $(1 \land \boldsymbol{x}^1)\lor (1 \land (\boldsymbol{\lnot x})^1)$, whence this condition can be equivalently phrased as the requirement that $\vdash$ has the IL and that $\emptyset \vdash (1 \land \boldsymbol{x}^1) \lor (1 \land (\boldsymbol{\lnot x})^1)$. It follows that Conditions (\ref{item:FLe-IL-1b}) and (\ref{item:FLe-IL-2b}) are equivalent.

Lastly, the equivalence between Conditions (\ref{item:FLe-IL-1b}) and (\ref{item:FLe-IL-3b}) follows from Theorem \ref{Thm:Pre-Lav-thm}.
\qed
\end{exa}

\paragraph{\bfseries Acknowledgements.}
Thanks are due to G. Bezhanishvili, N. Bezhanishvili, J. Czelakowski, R. Jansana, and J.G. Raftery for many useful suggestions which helped to improve the presentation of the paper.
The first author was supported by the FPI scholarship PRE$2020$-$093696$, associated with the I+D+i research project PID$2019$-$110843$GA-100 I+D \textit{La geometría de las lógicas no clásicas} funded by the Ministry of Science and Innovation of Spain.\ The second author was supported by the I+D+i research project PID$2019$-$110843$GA-I$00$ \emph{La geometria de las logicas no-clasicas} funded by the Ministry of Science and Innovation of Spain, by the \emph{Beatriz Galindo} grant BEAGAL\-$18$/$00040$ funded by the Ministry of Science and Innovation of Spain, and by the MSCA-RISE-Marie Skłodowska-Curie Research and Innovation Staff Exchange (RISE) project MOSAIC $101007627$ funded by Horizon $2020$ of the European Union.

\bibliographystyle{plain}
	
\end{document}